\documentclass[12pt]{article}

\usepackage{ifthen}
\newboolean{draftmode}
\setboolean{draftmode}{false}

\usepackage{profile-wang}

\newcommand{\PaperTitle}{
    Linear Convergence of Nesterov's Acceleration for All Momentum Parameter
}

\author{
    Hongda Li
    \thanks{
        Department of Mathematics, I.K. Barber Faculty of Science,
        The University of British Columbia, Kelowna, BC Canada V1V 1V7.
        E-mail: \texttt{alto@student.ubc.ca}.
    } and
    Xianfu Wang
    \thanks{
        Department of Mathematics, I.K. Barber Faculty of Science,
        The University of British Columbia, Kelowna, BC Canada V1V 1V7.
        E-mail: \texttt{shawn.wang@ubc.ca}.
    }
}

\title{{\fontfamily{ptm}\selectfont \PaperTitle}}

\begin{document}

\date{\today}
\maketitle
\todoinline{Haec charta nunc in modo schedarum est. Fontem inspice ut optiones mutes.}

\begin{abstract}
    \noindent
    We study the convergence rate of the Nesetrov's accelerated proximal gradient method for a smooth, strongly convex composite minimization objective across a range of value of the momentum extrapolation. 
    The literature proves a convergence rate of the method for a specific momentum extrapolation parameter choice based on the knowledge of strong convexity and Lipschitz smoothness constants.
    For the heavy ball method it is established that every choice of the momentum parameter either fails to accelerate or fails to converge on the same function class.
    We prove in this paper that a range of Nesterov style constant momentum extrapolation parameters produces a linear convergence rate.
    We corroborate these results numerically via performance estimation which computes a method's exact worst case rate as a semidefinite program.
\end{abstract}

\todoinline{\noindent
    \begin{enumerate}[nosep]
        \item\label{title:cand1} \textbf{[Recommended]} Linear Convergence of Nesterov's Acceleration at Every Admissible Momentum Parameter
        \item On the Global Optimality of Nesterov's Momentum Coefficient
        \item Robustness of Nesterov's Momentum: Linear Convergence at Every Extrapolation Parameter
        \item Closing the Momentum Gap: Linear Convergence of Accelerated Proximal Gradient Methods at Every Admissible $\beta$
        \item Nesterov's Acceleration Converges Linearly at Every Momentum Parameter
        \item Global Optimality and Linear Convergence of Nesterov's Momentum for Strongly Convex Composite Minimization
        \item The Momentum Coefficient Does Not Matter for Convergence: A Lyapunov Analysis of Nesterov's Acceleration
    \end{enumerate}
}


\section{Introduction}\label{sec:introduction}
    In this paper we study the composite minimization problem $\min_{x \in \RR^n} F(x) := f(x) + g(x)$ where $f$ is $\mu$-strongly convex and $L$-Lipschitz smooth, and $g$ is proper, closed, and convex.
    Nesterov's acceleration \cite{nesterov_method_1983} improves the $\mathcal O\left(\tfrac{1}{k}\right)$ rate of gradient descent to $\mathcal O\left(\tfrac{1}{k^2}\right)$ for smooth convex minimization using an increasing extrapolation sequence $(\lambda_k)$.
    The sequence satisfies the recursive relation $\lambda_{k - 1}^2 = (\lambda_k - 1)\lambda_k$ for all $k \in \N$.
    Under the additional assumption of strong convexity with modulus $\mu > 0$ the literature proposes using a specific constant momentum parameter instead to achieve a faster linear convergence rate \cite{beck_first-order_2017, nesterov_lectures_2018}.
    We study this constant momentum scheme.
    \par
    We write $T_L^{(f, g)}(x) := \argmin_{z \in \RR^n}\left\{ g(z) + \langle \nabla f(x), z\rangle + \tfrac{L}{2}\Vert z - x\Vert^2 \right\}$,
    which is proposed by Beck and Teboulle \cite[Theorem 10.42]{beck_fast_2009} for the convex composite objective.
    Below is our abstraction of Nesterov's constant momentum extrapolation relation for the convex composite objective.
    At each iteration $k$ the iterates satisfy:
    \begin{align}
        x_{k + 1} &= x_k^+ + \beta\left(x_k^+ - x_{k - 1}^+\right),
        \notag
        \\
        x_{k + 1}^+ &= T_L^{(f, g)}(x_{k + 1}).
        \label{expr:const-momentum-recursion}
    \end{align}
    In Nesterov's original work $g \equiv 0$ and $T_L^{(f, 0)}(x_{k + 1}) = x_{k + 1} - \tfrac{1}{L}\nabla f(x_{k + 1})$ is a gradient step.
    Beck and Teboulle \cite[Theorem 10.42]{beck_fast_2009} and Nesterov \cite{nesterov_lectures_2018} fix $\beta$ to be $\beta^\star := \tfrac{1 - \sqrt{\mu/L}}{1 + \sqrt{\mu/L}}$.
    This choice gives a fast linear convergence rate of $\mathcal O((1 - \sqrt{\mu/L})^k)$.
    \par
    The goal of the paper is to study the convergence rate of \eqref{expr:const-momentum-recursion} for all $\beta \in [0, 1)$ assuming $F$ is strongly convex.
    We prove that the method converges linearly for all $\beta \in [0, 1)$.
    \par
    To the best of our knowledge the literature gives a convergence rate at $\beta \in \{0, \beta^\star\}$ and, by the substitution described below, for $\beta \in [\beta^\star, 1)$ but gives no rate for $\beta \in (0, \beta^\star)$ \cite[Theorem 10.42, 10.21]{beck_first-order_2017}.
    We are interested in resolving this gap because $\mu$ is rarely known exactly in practice.
    In addition, the literature has recently established that every choice of the momentum parameter of the heavy ball method either fails to accelerate or fails to converge.
    Lessard et al. \cite{lessard_analysis_2016} and Goujaud et al. \cite{goujaud_provable_2025} prove that every such choice either fails to accelerate or fails to converge on the same function class.
    \subsection{Notations}\label{ssec:intro-notations}
        Throughout $[\cdot]_+ : \RR \to \RR_+ : x \mapsto \max(x, 0)$ denotes the positive part and $\Vert \cdot \Vert$ denotes the Euclidean norm of $\RR^n$ unless stated otherwise.
        We write $\N := \{1, 2, \ldots\}$ for the positive integers and $\Z_+ := \{0, 1, 2, \ldots\}$ for the nonnegative integers.
    \subsection{Literature reviews}\label{ssec:intro-lit-review}
        Several decades ago Nesterov \cite{nesterov_method_1983} introduced acceleration through an estimating sequence.
        Guler \cite{guler_new_1992} extended acceleration to the proximal point algorithm.
        Beck and Teboulle \cite{beck_fast_2009, beck_first-order_2017} later extended it to composite objectives with a nonsmooth term.
        The literature then branches into several directions.
        Restart schemes improve the practical performance of the algorithm \cite{odonoghue_adaptive_2015, alamo_gradient_2019, alamo_restart_2019, aujol_parameter-free_2024}.
        Dynamical systems give a continuous time interpretation of the algorithm \cite{su_differential_2016, aujol_optimal_2019, kim_unifying_2023, attouch_first-order_2022}.
        Lin et al. \cite{lin_catalyst_2018} implement Guler's method in the Catalyst framework.
        Other inexact variants then treat practical problems outside the Catalyst framework \cite{villa_accelerated_2013, schmidt_convergence_2011, bello-cruz_inexact_2020}.
        Lee et al. \cite{lee_geometric_2021} construct an optimal method that minimizes the squared gradient norm.
        Stochastic acceleration methods form a further direction \cite{shalev-shwartz_accelerated_2014, allen-zhu_katyusha_2018, zhou_simple_2018, zhou_direct_2019}.
        Linear convergence under conditions weaker than strong convexity holds for first order methods \cite{necoara_linear_2019} and for Nesterov's acceleration \cite{aujol_heavy_2025}.
        \par
        Next, we review the convergence results of the constant momentum methods in \eqref{expr:const-momentum-recursion} and the closely related heavy ball method in the literature.
        In the literature the convergence rates at $\beta = 0$ and at $\beta = \beta^\star$ are well established by Nesterov \cite{nesterov_lectures_2018} and Beck \cite{beck_first-order_2017}.
        The former value gives gradient descent or proximal gradient descent and it converges at a rate of $\mathcal O((1 - \frac{\mu}{L})^k)$ via cocoercivity \cite{beck_first-order_2017}.
        The latter is Nesterov's constant momentum method and it converges at a rate of $\mathcal O((1 - \sqrt{\tfrac{\mu}{L}})^k)$ \cite{nesterov_method_1983, nesterov_lectures_2018}.
        \par
        The rate at $\beta^\star$ extends to every $\beta \in [\beta^\star, 1)$ by an ad hoc argument.
        Write $\hat\mu := L\left(\tfrac{1 - \beta}{1 + \beta}\right)^2$ so that $\beta = \tfrac{1 - \sqrt{\hat\mu/L}}{1 + \sqrt{\hat\mu/L}}$ is the conventional momentum parameter for modulus $\hat\mu$ in place of $\mu$.
        For every $\beta \in [\beta^\star, 1)$ this $\hat\mu$ satisfies $\hat\mu \le \mu$ and it is itself a valid strong convexity modulus of $F$ whenever $\mu$ is.
        The V-FISTA convergence result \cite[Chapter 10]{beck_first-order_2017} applied to $F$ with $\hat\mu$ in place of $\mu$ then gives a convergence rate at $\beta$.
        The same substitution fails for $\beta < \beta^\star$ where it would require $\hat\mu > \mu$ which is not a valid modulus of strong convexity.
        The literature gives no convergence rate below $\beta^\star$ by this or any other argument.
        This is the gap in the literature which we study in this paper.
        \par
        Using a different approach, we derive the convergence rate as a function of $\beta$ on the entire range $\beta \in (0, 1)$ where $L$ and $\mu$ have nothing to do with $\beta$.
        We do not choose a specific $\beta$ and then verify our convergence rate.
        We prove the convergence rate for all $\beta \in (0, 1)$.
        Our argument applies for $\beta < \beta^\star$ where the substitution argument above gives nothing.
        It also applies for all $\beta \in [\beta^\star, 1)$ where our rate is tighter than the one given by that substitution.
        \par
        We note that numerical methods exist for showing the convergence rate of numerical optimization algorithms.
        Lessard et al. \cite{lessard_analysis_2016} propose the integral quadratic constraint framework.
        Drori and Teboulle \cite{drori_performance_2014} and Taylor et al. \cite{taylor_smooth_2016} propose the performance estimation problem.
        Both methods certify the exact worst case convergence of a fixed instance of an algorithm by solving a semidefinite program.
        Here ``a fixed instance'' means that every parameter of the algorithm must be specified.
        These methods therefore give no convergence rate for a generic family of algorithms.
        Nonetheless we use the same performance estimation problem in Section \ref{sec:pep-verification} to corroborate our theoretical results for specific choices of $\beta$ in \eqref{expr:const-momentum-recursion}.
        \par
        We now switch our attention to the closely related heavy ball method.
        Polyak \cite{polyak_methods_1964} introduced the heavy ball method which is given by $x_{k + 1}^+ = T_L^{(f, g)}(x_k^+) + \beta(x_k^+ - x_{k - 1}^+)$.
        It evaluates the proximal gradient at the current iterate $x_k^+$ instead of the extrapolated point $x_{k + 1}$ of \eqref{expr:const-momentum-recursion}.
        An adjacent result in the literature shows that every choice of the momentum parameter of the heavy ball method either fails to accelerate or fails to converge.
        Lessard et al. \cite{lessard_analysis_2016} identified a function on which it cycles.
        Goujaud et al. \cite{goujaud_provable_2025} proved that no choice of parameters makes the method both accelerate and converge on the general class of convex functions.
        \par
        Returning to Nesterov's extrapolation relation, Goujaud et al. \cite{goujaud_counter-examples_2023} searched numerically for cycling trajectories of \eqref{expr:const-momentum-recursion} in the smooth case where $g \equiv 0$.
        Their search covered all $\beta \in [0, 1)$ and all step sizes less than $\tfrac{2}{L}\tfrac{1 + \beta}{1 + 2\beta}$.
        They found no cycles for lengths $2$ to $25$.
        They concluded that Nesterov's method is more robust than the heavy ball method.
        The step size $\tfrac{1}{L}$ studied in our paper is less than $\tfrac{2}{L}\tfrac{1 + \beta}{1 + 2\beta}$ for all $\beta \in [0, 1)$.
        Their numerical finding for the smooth case therefore does not contradict our claim and it supports our theoretical findings instead.
        \par
        Therefore whether an analytical convergence rate exists for Nesterov's method for all $\beta \in [0, 1)$ is a gap in the literature.
        We resolve it positively and establish a linear convergence rate for all $\beta \in [0, 1)$.
        We also verify our theoretical findings using numerical experiments.
    
    \subsection{Our contributions}\label{ssec:intro-contributions}
        Our paper makes three contributions.
        \begin{enumerate}[nosep]
            \item We propose a theoretical convergence framework based on the abstract formulation of Nesterov's extrapolated momentum in Assumption \ref{ass:abst-nes-extrp}.
            Under this framework we develop a Lyapunov inequality in Lemma \ref{lemma:apg-lyapunov-inequality}, and then strengthen it under strong convexity in Lemma \ref{lemma:apg-lya-ineq-mu}, which is essential to the proofs that follow.
            \item We identify the constant momentum recursion \eqref{expr:const-momentum-recursion} as a special case of the accelerated proximal gradient algorithm in Theorem \ref{thm:apg-momentum-equiv}.
            We derive the contraction factor $\overline M_\sigma(\beta)$ of \eqref{def:overline-m} which equals $M_\sigma^*(\beta)$ for $\beta \in [\beta^\star, 1)$ by Proposition \ref{prop:mstar-beta-ge-betastar} and $1 > \overline M_\sigma(\beta) \ge M_\sigma^*(\beta)$ for $\beta \in (0, \beta^\star]$ by Proposition \ref{prop:mstar-beta-le-bstar}.
            We prove the method converges linearly at every $\beta \in (0, 1)$ on the composite class in Theorem \ref{thm:lin-cnvg-constn-mmntm}.
            This is our main contribution which addresses the gap in the literature. 
            \item
            We verify numerically using the performance estimation problem (PEP) \cite{drori_performance_2014, taylor_smooth_2016} in Section \ref{sec:pep-verification} for certain choices of $\beta, \sigma$.
        \end{enumerate}
        \textbf{The paper is organized as follows.}
        Section \ref{sec:preliminaries} collects the standing assumptions and the algebraic identities.
        Section \ref{sec:aux-results} identifies the accelerated proximal gradient algorithm and develops its Lyapunov inequality.
        Section \ref{sec:major-results} presents two results of our paper. 
        The first shows the equivalence between the constant momentum algorithm and the accelerated proximal gradient algorithm proposed by us, enabling the use of our Lyapunov inequality.
        The second main result uses this equivalence to prove the convergence rate of the constant momentum algorithm for all $\beta \in (0, 1)$.
        Section \ref{sec:pep-verification} certifies our convergence rate using the exact worst case value computed via PEP on many discrete instances of $\beta, \sigma$.

\section{Preliminaries}\label{sec:preliminaries}
    This section collects the standing assumptions, notation, and elementary facts used throughout the paper.
    \subsection{Assumptions and the overall setup}\label{ssec:assumptions-setup}
        We introduce the Bregman divergence in the definition below.
        \begin{definition}[Bregman divergence]\;\label{def:bregman-divergence}\\
            Let $f : \RR^n \to \RR$ be a differentiable function. 
            The Bregman divergence is the mapping:
            \begin{align*}
                D_f : \RR^n \times \RR^n \to \RR : (x, y) \mapsto f(x) - f(y) - \langle \nabla f(y), x - y\rangle.
            \end{align*}
        \end{definition}
        We use the Bregman divergence to define strong convexity of a differentiable function. 
        \begin{definition}[Strong convexity]\;\label{def:strong-convexity}\\
            Let $f : \RR^n \to \RR$ be a differentiable function and let $\mu \ge 0$. 
            We say $f$ is strongly convex with modulus $\mu$ if for all $x \in \RR^n$ and $y \in\RR^n$ it satisfies:
            \begin{align*}
                D_f(x, y) \ge \frac{\mu}{2}\Vert x - y\Vert^2.
            \end{align*}
        \end{definition}
        The following is the assumption we make throughout the paper unless stated otherwise.
        \begin{assumption}[Global assumptions]\;\label{ass:global-setup}\\
            Assume that $(F, f, g, L, \mu)$ satisfy the following.
            \begin{enumerate}[nosep]
                \item\label{ass:global-setup:item1} $f : \RR^n \to \RR$ is convex and Lipschitz smooth with constant $L$, i.e., $D_f(x, y) \le \tfrac{L}{2}\Vert x - y\Vert^2$ for all $x, y \in \RR^n$ where $D_f(x, y)$ is given by Definition \ref{def:bregman-divergence}.
                \item\label{ass:global-setup:item2} $f$ is strongly convex with modulus $\mu \ge 0$ as defined in Definition \ref{def:strong-convexity}.
                \item\label{ass:global-setup:item3} $g : \RR^n \to \overline\RR$ is proper, closed, and convex.
                \item\label{ass:global-setup:item4} The overall objective is $F = f + g$.
            \end{enumerate}
        \end{assumption}
        Let $(F, f, g, L, \mu)$ satisfy Assumption \ref{ass:global-setup}.
        We consider the following optimization problem in this paper: $\min_{x \in \RR^n}\left\lbrace f(x) + g(x)\right\rbrace$.
    \subsection{Euclidean identities}\label{ssec:euclidean-identities}
        We record an elementary Euclidean identity in the fact below.
        \begin{fact}[{\cite[Lemma 1.7]{bauschke_convex_2017}}]\;\label{fact:ecld-id}\\
            Let $\lambda \in \RR$ and $a, b \in \RR^n$ be arbitrary. Then the following identity holds:
            \begin{align*}
                (1 - \lambda)\Vert a\Vert^2 + \lambda\Vert b\Vert^2
                = \Vert (1 - \lambda)a + \lambda b\Vert^2 + \lambda(1 - \lambda)\Vert a - b\Vert^2.
            \end{align*}
        \end{fact}
        The following lemma extends Fact \ref{fact:ecld-id}.
        \begin{lemma}\;\label{lemma:ecld-id-extd}\\
            Let $\lambda \in \RR$ and $\bar x, x \in \RR^n$ be arbitrary. 
            Then for all $w \in \RR^n$ the following holds:
            \begin{align*}
                (1 - \lambda)\Vert w - x\Vert^2 + \lambda\Vert w - \bar x\Vert^2
                = \Vert w - (1 - \lambda)x - \lambda \bar x\Vert^2 + \lambda(1 - \lambda)\Vert x - \bar x\Vert^2.
            \end{align*}
        \end{lemma}
        \begin{proof}
            \todoinline{
                Lemma \ref{lemma:ecld-id-extd} proof depends on:
                \begin{itemize}[nosep]
                    \item Fact \ref{fact:ecld-id}.
                \end{itemize}
            }
            Apply Fact \ref{fact:ecld-id} with $a = w - x$ and $b = w - \bar x$.
        \end{proof}
    \subsection{Abstract Nesterov's extrapolation form}\label{ssec:abstract-momentum-form}
        We propose the following abstract relation satisfied by the iterates of Nesterov's accelerated gradient method.
        The name ``abstract'' refers to two characteristics of the assumption.
        First, it only provides the relationship between the iterates without a proximal gradient step or a gradient descent step.
        Second, there are two indexed sequences $\alpha_k, q_k$ that control the recurrence relations instead of the fixed constants used in the literature.
        \begin{assumption}[Abstract Nesterov's extrapolation sequences]\;\label{ass:abst-nes-extrp}\\
            Let $x_{-1}^+, x_{-1}^\dagger$ be fixed.  
            The abstract Nesterov's extrapolation sequences are three sequences $(x_k)_{k \in \Z_+}, (x_k^+)_{k \in \Z_+}, (x_k^\dagger)_{k \in \Z_+}$ in $\RR^n$ such that there exist $(\alpha_k)_{k \in \Z_+}, (q_k)_{k \in \Z_+}$ and the following are true for all $k \in \Z_+$.
            \begin{enumerate}[nosep, label=(\roman*)]
                \item\label{ass:abst-nes-extrp-alphaq} $\alpha_k \in (0, 1], q_k \in [0, \alpha_k]\setminus \{1\}$.
                \item\label{ass:abst-nes-extrp-lambdak} $\lambda_k = \frac{\alpha_k - q_k}{1 - q_k}$.
                \item\label{ass:abst-nes-extrp-x} $x_k = \lambda_kx_{k - 1}^\dagger + \left(1 - \lambda_k\right)x_{k - 1}^+$.
                \item\label{ass:abst-nes-extrp-xd} $x_k^\dagger = x_{k - 1}^+ + \alpha_k^{-1}(x_k^+ - x_{k - 1}^+).$
            \end{enumerate}
        \end{assumption}
        \begin{remark}
            This assumption is a generalization of the SC-FGM form from Lee et al. \cite{lee_geometric_2021}.
            Setting $\alpha_k \equiv \sqrt{\frac{\mu}{L}}$ and $q_k \equiv \frac{\mu}{L}$ with initial conditions $x_{-1}^\dagger = x_{-1}^+$ recovers that form. 
        \end{remark}
        The assumption above highlights the algebraic properties of the iterates in Definition \ref{def:apg-algorithm}.
        It is used in the derivation of Lemma \ref{lemma:apg-lyapunov-inequality}.

\section{Auxiliary Results}\label{sec:aux-results}
    This section develops the inequality of Lemma \ref{lemma:apg-lyapunov-inequality} which is used to derive the convergence rate later in the paper.
    To start, the definition below introduces the proximal gradient operator for a composite objective consisting of a smooth convex function and a nonsmooth function.
    \begin{definition}[The proximal gradient operator]\;\label{def:prox-grad-operator}\\
        Let $(F, f, g, L, \mu)$ satisfy Assumption \ref{ass:global-setup}. 
        We define the proximal gradient operator for all $\gamma > 0$ by
        \begin{align*}
            T_\gamma^{(f, g)}(x) := \argmin_{z \in \RR^n}\left\{ g(z) + \langle \nabla f(x), z\rangle + \frac{\gamma}{2}\Vert z - x\Vert^2 \right\}.
        \end{align*}
    \end{definition}
    The proximal gradient operator provides us with the following inequality.
    \begin{fact}[The proximal gradient inequality {\cite[(10.37)]{beck_fast_2009}}]\;\label{fact:prox-grad-inequality}\\
        Let $(F, f, g, L, \mu)$ satisfy Assumption \ref{ass:global-setup}. 
        Then the following is true.
        \begin{enumerate}[nosep, label=(\roman*)]
            \item $T_\gamma^{(f, g)}(x)$ is a singleton which we can write $x^+ = T_\gamma^{(f, g)}(x)$ for short.
            \item for all $z \in \RR^n$ the following holds:
            \begin{align*}
                F(z) - F(x^+)
                \ge
                \frac{\gamma}{2}\Vert z - x^+\Vert^2
                - \frac{\gamma - \mu}{2}\Vert z - x\Vert^2.
            \end{align*}
        \end{enumerate}
    \end{fact}
    We propose the accelerated proximal gradient (APG) algorithm below.
    It satisfies the abstract extrapolation form of Assumption \ref{ass:abst-nes-extrp}.
    \begin{definition}[The accelerated proximal gradient algorithm]\;\label{def:apg-algorithm}\\
        Suppose that $(F, f, g, L, \mu)$ satisfy Assumption \ref{ass:global-setup}. 
        The accelerated proximal gradient algorithm generates three sequences $(x_k)_{k \in \Z_+}$, $(x_k^+)_{k \in \Z_+}$, $(x_k^\dagger)_{k \in \Z_+}$ along with parameters $(\alpha_k)_{k \in \Z_+}$, $(q_k)_{k \in \Z_+}$, $(\lambda_k)_{k \in \Z_+}$ such that they satisfy Assumption \ref{ass:abst-nes-extrp} and in addition there exists a line search sequence $(\gamma_k)_{k \in \Z_+}$ such that for all $k \in \Z_+$:
        \begin{align}
            &x_k^+ := T_{\gamma_k}^{(f, g)}(x_k),
            \label{def:apg-algorithm-xp}
            \\
            &D_f(x_k^+, x_k) \le \frac{\gamma_k}{2}\Vert x_k^+ - x_k\Vert^2,
            \label{def:apg-algorithm-ls}
        \end{align}
        where $T_{\gamma_k}^{(f, g)}$ is given by Definition \ref{def:prox-grad-operator}.
    \end{definition}
    \begin{remark}
        The constraint $q_k \le \alpha_k$ of Assumption \ref{ass:abst-nes-extrp}\ref{ass:abst-nes-extrp-alphaq} makes $x_k$ a convex combination of $x_{k - 1}^+$ and $x_{k - 1}^\dagger$.
        Assumption \ref{ass:global-setup} implies $D_f(x_k^+, x_k) \le \tfrac{L}{2}\Vert x_k^+ - x_k\Vert^2$ because of the Lipschitz continuity of $\nabla f$.
        More precisely it implies that the condition is automatically true for all $\gamma_k \ge L$.
        Furthermore, strong convexity implies $D_f(x_k^+, x_k) \ge \tfrac{\mu}{2}\Vert x_k^+ - x_k\Vert^2$ and hence $\gamma_k \ge \mu$.
    \end{remark}
    The lemma below links the abstract formulation of Definition \ref{def:apg-algorithm} to the constant momentum formulation of \eqref{expr:const-momentum-recursion}.
    \begin{lemma}[Abstract momentum form]\;\label{lemma:abstract-momentum-form}\\
        Let $(x_k)_{k \in \Z_+}, (x_k^+)_{k \in \Z_+}, (x_k^\dagger)_{k \in \Z_+}$ satisfy Assumption \ref{ass:abst-nes-extrp}. Then for all $k \in \Z_+$ the following holds:
        \begin{align*}
            x_{k + 1} = x_k^+ + 
            \lambda_{k + 1}\frac{1 - \alpha_k}{\alpha_k}\left(x_k^+ - x_{k - 1}^+\right).
        \end{align*}
        Here, $\lambda_{k + 1}$ is given by Assumption \ref{ass:abst-nes-extrp}\ref{ass:abst-nes-extrp-lambdak}.
    \end{lemma}
    \begin{proof}
        \todoinline{
            Lemma \ref{lemma:abstract-momentum-form} proof depends on:
            \begin{itemize}[nosep]
                \item Assumption \ref{ass:abst-nes-extrp}\ref{ass:abst-nes-extrp-x}.
                \item Assumption \ref{ass:abst-nes-extrp}\ref{ass:abst-nes-extrp-xd}.
            \end{itemize}
        }
        From Assumption \ref{ass:abst-nes-extrp}\ref{ass:abst-nes-extrp-x} the following holds:
        \begin{align*}
            x_{k + 1} &= \lambda_{k + 1}x_k^\dagger + (1 - \lambda_{k + 1})x_k^+
            \\
            &\underset{(1)}{=} x_k^+ + \lambda_{k + 1}\left(x_k^\dagger - x_k^+\right)
            \\
            &\underset{(2)}{=} x_k^+ + \lambda_{k + 1}\frac{1 - \alpha_k}{\alpha_k}\left(x_k^+ - x_{k - 1}^+\right).
        \end{align*}
        At (1) we added and subtracted $x_k^+$ isolating the momentum term as a multiple of $x_k^\dagger - x_k^+$.
        At (2) we substituted the recursion for $x_k^\dagger$ at index $k$ given by Assumption \ref{ass:abst-nes-extrp}\ref{ass:abst-nes-extrp-xd}. It gives $x_k^\dagger - x_k^+ = (x_{k - 1}^+ - x_k^+) + \alpha_k^{-1}(x_k^+ - x_{k - 1}^+) = \tfrac{1 - \alpha_k}{\alpha_k}(x_k^+ - x_{k - 1}^+)$.
        Substituting $\lambda_{k + 1}$ back completes the proof.
    \end{proof}
    \begin{remark}
        The map $(\alpha_k, \alpha_{k + 1}, q_k, q_{k + 1}) \mapsto \lambda_{k + 1}\tfrac{1 - \alpha_k}{\alpha_k}$ is not injective.
        Many choices of $\alpha_k, \alpha_{k + 1}, q_k$ and $q_{k + 1}$ produce a $\beta_{k + 1}$ such that $x_{k + 1} = x_k^+ + \beta_{k + 1}(x_k^+ - x_{k - 1}^+)$.
        An algorithm fitting the extrapolated momentum update by $\beta_{k + 1}$ therefore maps back to a whole family of instances of Assumption \ref{ass:abst-nes-extrp}.
    \end{remark}
    The following three lemmas are used in Section \ref{ssec:apg-lya-ineq}. 
    \begin{lemma}[Abstract deviation form]\;\label{lemma:abstract-dvt}\\
        Let $\bar x \in \RR^n$ be arbitrary and let $(x_k)_{k \in \Z_+}, (x_k^+)_{k \in \Z_+}, (x_k^\dagger)_{k \in \Z_+}$ and $(\lambda_k)_{k\in\Z_+}$ satisfy Assumption \ref{ass:abst-nes-extrp}. 
        Then for all $k \in \Z_+$ the following are satisfied:
        \begin{align*}
            \alpha_k(x_{k - 1}^\dagger - \bar x) &= x_k - z_k + q_k(x_{k - 1}^\dagger - x_k),
            \\
            \alpha_k(x_k^\dagger - \bar x) &= x_k^+ - z_k.
        \end{align*}
    \end{lemma}
    \begin{proof}
        See Section \ref{ssec:abstract-dvt-pf}. 
    \end{proof}
    \begin{lemma}[Abstract distance splitting form]\;\label{lemma:abstract-dist-splt}\\
        Let $\bar x \in \RR^n$ be arbitrary and let $(x_k)_{k \in \Z_+}, (x_k^+)_{k \in \Z_+}, (x_k^\dagger)_{k \in \Z_+}$ satisfy Assumption \ref{ass:abst-nes-extrp}.
        Here, $(\lambda_k)_{k \in \Z_+}$ is given by Assumption \ref{ass:abst-nes-extrp}\ref{ass:abst-nes-extrp-lambdak}.
        Then for all $k \in \Z_+$ the following holds:
        \begin{align*}
            \Vert \bar x - x_k\Vert^2
            = (1 - \lambda_k)\Vert \bar x - x_{k - 1}^+\Vert^2
            + \lambda_k\left\Vert \bar x - x_{k - 1}^\dagger\right\Vert^2
            - \lambda_k(1 - \lambda_k)\left\Vert x_{k - 1}^\dagger - x_{k - 1}^+\right\Vert^2.
        \end{align*}
    \end{lemma}
    \begin{proof}
        See Section \ref{ssec:abstract-dist-splt-pf}.
    \end{proof}
    \begin{lemma}[Abstract gap splitting form]\;\label{lemma:abstract-gp-splt}\\
        Let $\bar x \in \RR^n$ be arbitrary and let $(x_k)_{k \in \Z_+}, (x_k^+)_{k \in \Z_+}, (x_k^\dagger)_{k \in \Z_+}$ satisfy Assumption \ref{ass:abst-nes-extrp}.
        Here, $(\lambda_k)_{k \in \Z_+}$ is given by Assumption \ref{ass:abst-nes-extrp}\ref{ass:abst-nes-extrp-lambdak}. 
        Let $z_k := \alpha_k \bar x + (1 - \alpha_k)x_{k - 1}^+$. 
        Then for all $k \in \Z_+$ the following holds:
        \begin{align*}
            &\frac{1}{2}\Vert x_k - z_k\Vert^2
            \\
            &= \frac{\alpha_k \lambda_k}{2}\left\Vert x_{k - 1}^\dagger - \bar x\right\Vert^2
            + \frac{\alpha_k q_k(1 - \lambda_k)}{2}\left\Vert x_{k - 1}^+ - \bar x\right\Vert^2
            - \frac{q_k \lambda_k(1 - \lambda_k)}{2}\left\Vert x_{k - 1}^\dagger - x_{k - 1}^+\right\Vert^2.
        \end{align*}
    \end{lemma}
    \begin{proof}
        See Section \ref{ssec:abstract-gp-splt-pf}.
    \end{proof}

    \subsection{APG Lyapunov Inequality}\label{ssec:apg-lya-ineq}
        Choose initial conditions $x_{-1}^+, x_{-1}^\dagger \in \RR^n$.
        Let $(F, f, g, L, \mu)$ be given by Assumption \ref{ass:global-setup}. 
        Let $(x_k^+)_{k \in \Z_+}, (x_k^\dagger)_{k \in \Z_+}$ be sequences in $\RR^n$ and let $(\alpha_k)_{k \in \Z_+}$, $(q_k)_{k \in \Z_+}$, $(\gamma_k)_{k \in \Z_+}$, $(\lambda_k)_{k \in \Z_+}$ be sequences in $\RR$ such that they satisfy Definition \ref{def:apg-algorithm}. 
        For all $\bar x \in \RR^n$, $k \in \N$ we define:
        \begin{align}
            \label{def:delta}
            \Delta_k(\bar x) &:= F(x_k^+) - F(\bar x),
            \\
            \label{def:upsilon}
            \Upsilon_k(\bar x) &:= \Delta_k(\bar x) + \frac{\gamma_k \alpha_k^2}{2}\Vert x_k^\dagger - \bar x\Vert^2,
            \\
            \label{def:ck}
            C_k &:= 
            \frac{(\gamma_k - \mu)\alpha_k \lambda_k}{\gamma_{k - 1}\alpha_{k - 1}^2},
            \\
            \label{def:xi-k}
            \Xi_k(\bar x) 
            &:= \frac{\alpha_k(1 - \lambda_k)(\gamma_k q_k - \mu)}{2}\left\Vert x_{k - 1}^+ - \bar x \right\Vert^2.
        \end{align}
        The following lemma derives a descent inequality for the family of algorithms that satisfy Definition \ref{def:apg-algorithm} under Assumption \ref{ass:global-setup}.
        \begin{lemma}[APG Lyapunov inequality]\;\label{lemma:apg-lyapunov-inequality}\\
            Let $x_{-1}^+, x_{-1}^\dagger \in \RR^n$ be any initial conditions.
            Let $(x_k)_{k \in \Z_+}, (x_k^+)_{k \in \Z_+}, (x_k^\dagger)_{k \in \Z_+}$, $(\alpha_k)_{k \in \Z_+}$, $(q_k)_{k \in \Z_+}$, $(\lambda_k)_{k \in \Z_+}$ and $(\gamma_k)_{k \in \Z_+}$ be given by Definition \ref{def:apg-algorithm}. 
            Let $\Delta_k$, $\Upsilon_k$, $C_k$ and $\Xi_k$ be given by \eqref{def:delta}, \eqref{def:upsilon}, \eqref{def:ck}, and \eqref{def:xi-k} respectively.
            Then the following statements are true for all minimizers $\bar x \in \RR^n$.
            \begin{enumerate}[label=(\roman*)]
                \item For all $k \in \N$ the following inequality holds:
                \begin{align}
                    \Upsilon_k(\bar x)
                    \le
                    \max\left(1 - \alpha_k,\, C_k\right)
                    \Upsilon_{k - 1}(\bar x) + \Xi_k(\bar x).
                    \label{eqn:apg-lyapunov-main}
                \end{align}
                \item The following inequality holds
                \begin{align}
                    \Upsilon_0(\bar x) \le
                    \frac{\gamma_0}{2}\Vert \bar x - x_0\Vert^2
                    - \frac{\gamma_0}{2}\left\Vert \bar x - x_0^+ \right\Vert^2
                    + \frac{\gamma_0 \alpha_0^2}{2}\Vert x_0^\dagger - \bar x\Vert^2.
                \label{eqn:apg-lyapunov-base}
                \end{align}
                \item If additionally $x_{-1}^+ = x_{-1}^\dagger$ and $\alpha_0 = 1$ then $x_0 = x_{-1}^+ = x_{-1}^\dagger$ and:
                \begin{align}
                    \label{eqn:apg-lyapunov-base-special}
                    \Upsilon_0(\bar x) \le \frac{\gamma_0}{2}\Vert \bar x - x_0\Vert^2.
                \end{align}
            \end{enumerate}
        \end{lemma}
        \begin{proof}
            \todoinline{
                Lemma \ref{lemma:apg-lyapunov-inequality} proof depends on:
                \begin{itemize}[nosep]
                    \item Fact \ref{fact:prox-grad-inequality}.
                    \item Lemma \ref{lemma:ecld-id-extd}.
                    \item Lemma \ref{lemma:abstract-dvt}.
                    \item Lemma \ref{lemma:abstract-gp-splt}.
                    \item Assumption \ref{ass:abst-nes-extrp}.
                    \item Assumption \ref{ass:abst-nes-extrp}\ref{ass:abst-nes-extrp-lambdak}.
                    \item Assumption \ref{ass:abst-nes-extrp}\ref{ass:abst-nes-extrp-x}.
                    \item Assumption \ref{ass:abst-nes-extrp}\ref{ass:abst-nes-extrp-xd}.
                    \item Definition \ref{def:apg-algorithm}.
                    \item Definition \ref{def:strong-convexity}.
                    \item \eqref{def:delta}.
                    \item \eqref{def:upsilon}.
                    \item \eqref{def:xi-k}.
                \end{itemize}
            }
            The proof is a Lyapunov argument.
            First, we establish for all $k \in \Z_+$:
            \begin{align*}
                & (1 - \alpha_k)\left(F(x_{k - 1}^+) - F(\bar x)\right)
                - \left(F(x_k^+) - F(\bar x)\right)
                \\
                &= (1 - \alpha_k)F(x_{k - 1}^+) - (1 - \alpha_k)F(\bar x)
                - F(x_k^+) + F(\bar x)
                \\
                &= (1 - \alpha_k)F(x_{k - 1}^+) - F(x_k^+) + \alpha_k F(\bar x)
                \\
                &= (1 - \alpha_k)\left(F(x_{k - 1}^+) - F(x_k^+) \right)
                + \alpha_k\left(F(\bar x) - F(x_k^+)\right).
            \end{align*}
            Recall $\Delta_k(\bar x) := F(x_k^+) - F(\bar x)$ and $z_k := \alpha_k\bar x + (1 - \alpha_k)x_{k - 1}^+$ for all $k \in \Z_+$. 
            Using the above, we deduce
            \begin{align*}
                &(1 - \alpha_k)\Delta_{k - 1}(\bar x) - \Delta_{k}(\bar x)
                \\
                &= (1 - \alpha_k)\left(F(x_{k - 1}^+) - F(x_k^+)\right)
                + \alpha_k\left(F(\bar x) - F(x_k^+)\right)
                \\
                &\underset{(1)}\ge
                (1 - \alpha_k)\left(
                    \frac{\gamma_k}{2}
                    \left\Vert x_{k - 1}^+ - x_k^+ \right\Vert^2
                    -
                    \frac{\gamma_k - \mu}{2}\left\Vert x_{k - 1}^+ - x_k \right\Vert^2
                \right) 
                    \\&\quad
                    + \alpha_k\left(
                        \frac{\gamma_k}{2}\left\Vert \bar x - x_k^+ \right\Vert^2
                        - \frac{\gamma_k - \mu}{2}\Vert \bar x - x_k \Vert^2
                    \right).
            \end{align*}
            At (1) we applied Fact \ref{fact:prox-grad-inequality} with $x = x_k$, $x^+ = x_k^+$, $\gamma = \gamma_k$ for two values of $z$: $z = \bar x$ and  $z = x_{k - 1}^+$.
            Continuing and rearranging the RHS of (1):
            \begin{align*}
                &\frac{\gamma_k}{2}\left(
                    (1 - \alpha_k)\left\Vert x_{k - 1}^+ - x_k^+\right\Vert^2
                    + \alpha_k\left\Vert \bar x - x_k^+ \right\Vert^2
                \right)
                    \\&\quad
                    - \frac{\gamma_k - \mu}{2}
                    \left(
                        (1 - \alpha_k)\left\Vert x_{k - 1}^+ - x_k \right\Vert^2
                        + \alpha_k\Vert \bar x - x_k \Vert^2
                    \right)
                \\
                &\underset{(2)}{=}
                \frac{\gamma_k}{2}\left(
                    \left\Vert x_k^+ - z_k\right\Vert^2
                    + \alpha_k(1 - \alpha_k)\left\Vert x_{k - 1}^+ - \bar x\right\Vert^2
                \right)
                    \\&\quad
                    - \frac{\gamma_k - \mu}{2}\left(
                        \Vert x_k - z_k\Vert^2
                        + \alpha_k(1 - \alpha_k)\left\Vert x_{k - 1}^+ - \bar x \right\Vert^2
                    \right)
                \\
                &\underset{(3)}{=}
                \frac{\gamma_k\alpha_k^2}{2}\left\Vert x_k^\dagger - \bar x\right\Vert^2
                - \frac{(\gamma_k - \mu)\alpha_k\lambda_k}{2}\left\Vert x_{k - 1}^\dagger - \bar x\right\Vert^2
                - \frac{\alpha_k(1 - \lambda_k)(\gamma_kq_k - \mu)}{2}\left\Vert x_{k - 1}^+ - \bar x \right\Vert^2
                    \\&\quad 
                    + \underbrace{\frac{(\gamma_k - \mu)q_k\lambda_k(1 - \lambda_k)}{2}\left\Vert 
                        x_{k - 1}^\dagger - x_{k - 1}^+ 
                    \right\Vert^2}_{\ge 0}.
            \end{align*}
            At (2) we applied Lemma \ref{lemma:ecld-id-extd} twice at $(x, \bar x, \lambda) = (x_{k - 1}^+, \bar x, \alpha_k)$ once with $w = x_k^+$ and coefficient $\gamma_k$, and once with $w = x_k$.
            At (3) we applied Lemma \ref{lemma:abstract-dvt} to $\Vert x_k^+ - z_k\Vert^2$ and Lemma \ref{lemma:abstract-gp-splt} to $\Vert x_k - z_k\Vert^2$ then collected the coefficients of $\Vert x_{k - 1}^+ - \bar x\Vert^2$. 
            Expand and group $\Vert x_{k - 1}^+ - \bar x \Vert^2$ on right hand side of (2) to produce:
            \begin{align*}
                &\frac{\gamma_k}{2}\left\Vert x_k^+ - z_k \right\Vert^2 
                - \frac{\gamma_k - \mu}{2}\Vert x_k - z_k\Vert^2
                + \frac{\mu\alpha_k(1 - \alpha_k)}{2}\left\Vert x_{k - 1}^+ - \bar x \right\Vert^2
                \\
                &\underset{(3)}{=}
                \frac{\gamma_k\alpha_k^2}{2}\left\Vert x_k^\dagger - \bar x \right\Vert^2
                - \frac{(\gamma_k - \mu)\alpha_k\lambda_k}{2}\Vert x_{k - 1}^\dagger - \bar x\Vert^2
                - \frac{(\gamma_k - \mu)\alpha_kq_k(1 - \lambda_k)}{2}\left\Vert x_{k - 1}^+ - \bar x \right\Vert^2
                \\&\quad
                + \frac{\mu\alpha_k(1 - \alpha_k)}{2}\left\Vert x_{k - 1}^+ - \bar x \right\Vert^2
                + \frac{(\gamma_k - \mu)q_k\lambda_k(1 - \lambda_k)}{2}\left\Vert x_{k - 1}^\dagger - x_{k - 1}^+ \right\Vert^2
                \\
                &=
                \frac{\gamma_k\alpha_k^2}{2}\Vert x_k^\dagger - \bar x\Vert^2
                - \frac{(\gamma_k - \mu)\alpha_k\lambda_k}{2}\Vert x_{k - 1}^\dagger - \bar x\Vert^2
                - \frac{\alpha_k(1 - \lambda_k)(\gamma_kq_k - \mu)}{2}\left\Vert x_{k - 1}^+ - \bar x \right\Vert^2
                \\&\quad
                + \frac{(\gamma_k - \mu)q_k\lambda_k(1 - \lambda_k)}{2}\left\Vert x_{k - 1}^\dagger - x_{k - 1}^+ \right\Vert^2. 
            \end{align*}
            The coefficient of $\Vert x_{k - 1}^+ - \bar x\Vert^2$ follows from substituting $1 - \alpha_k = (1 - q_k)(1 - \lambda_k)$ obtained from Assumption \ref{ass:abst-nes-extrp}\ref{ass:abst-nes-extrp-lambdak}:
            \begin{align*}
                - \frac{(\gamma_k - \mu)\alpha_kq_k(1 - \lambda_k)}{2} + \frac{\mu\alpha_k(1 - \alpha_k)}{2}
                &= - \frac{(\gamma_k - \mu)\alpha_kq_k(1 - \lambda_k)}{2} + \frac{\mu\alpha_k(1 - q_k)(1 - \lambda_k)}{2}
                \\
                &= \frac{\alpha_k(1 - \lambda_k)}{2}\Big[\mu(1 - q_k) - (\gamma_k - \mu)q_k\Big]
                \\
                &= \frac{\alpha_k(1 - \lambda_k)}{2}(\mu - \gamma_kq_k)
                \\
                &= -\frac{\alpha_k(1 - \lambda_k)(\gamma_kq_k - \mu)}{2}.
            \end{align*}
            Finally, the coefficient of $\Vert x_{k - 1}^\dagger - x^+_{k - 1}\Vert^2$ on the RHS of equality (3) is nonnegative because $\gamma_k \ge \mu$, $q_k \ge 0$, $\gamma_k \in [0, 1]$ by Definition \ref{def:apg-algorithm} and Definition \ref{def:strong-convexity}.
            For all $k \in \Z_+$ the chain of inequalities (1), (2) and (3) in this proof yields:
            \begin{align*}
                &(1 - \alpha_k)\Delta_{k - 1}(\bar x) - \Delta_{k}(\bar x)
                \\
                &\ge
                \frac{\gamma_k\alpha_k^2}{2}\left\Vert x_k^\dagger - \bar x\right\Vert^2
                - \frac{(\gamma_k - \mu)\alpha_k\lambda_k}{2}\left\Vert x_{k - 1}^\dagger - \bar x\right\Vert^2
                - \frac{\alpha_k(1 - \lambda_k)(\gamma_kq_k - \mu)}{2}\left\Vert x_{k - 1}^+ - \bar x \right\Vert^2.
            \end{align*}
            Rearranging the above and using the notations from \eqref{def:delta}, \eqref{def:upsilon} and \eqref{def:xi-k} yields for all $k \ge 1$:
            \begin{align*}
                \Upsilon_k(\bar x) &= \Delta_k(\bar x) + \frac{\gamma_k\alpha_k^2}{2}\left\Vert x_k^\dagger - \bar x\right\Vert^2
                \\
                &\le (1 - \alpha_k)\Delta_{k - 1}(\bar x) + \frac{(\gamma_k - \mu) \alpha_k \lambda_k}{2}\left\Vert
                    x_{k - 1}^\dagger - \bar x
                \right\Vert^2
                + \Xi_k(\bar x)
                \\
                &=
                (1 - \alpha_k)\Delta_{k - 1}(\bar x) + C_k\cdot\frac{\gamma_{k - 1}\alpha_{k - 1}^2}{2}\left\Vert
                    x_{k - 1}^\dagger - \bar x
                \right\Vert^2
                + \Xi_k(\bar x)
                \\
                &\underset{(4)}\le
                \max\left(
                    1 - \alpha_k, C_k
                \right)\Upsilon_{k - 1}(\bar x) + \Xi_k(\bar x).
            \end{align*}
            At (4) we used the fact that $\bar x$ is a minimizer of $F$. 
            This makes $\Delta_{k - 1}(\bar x) \ge 0$ and hence we can replace both coefficients $(1 - \alpha_k), C_k$ by their maximum
            which establishes \eqref{eqn:apg-lyapunov-main}.
            \par
            However, $C_k$ is not defined at $k = 0$ so we derive separately an upper bound for $\Upsilon_0(\bar x)$.
            We apply Fact \ref{fact:prox-grad-inequality} with $z = \bar x$ and $x = x_0$ which produces $x^+ = T_{\gamma_0}^{(f, g)}(x_0) = x_0^+$ and gives:
            \begin{align*}
                F(\bar x) - F(x_0^+)
                &\ge
                \frac{\gamma_0}{2}\left\Vert \bar x - x_0^+ \right\Vert^2
                - \frac{\gamma_0}{2}\Vert \bar x - x_0 \Vert^2
                + \frac{\mu}{2}\Vert x_0 - \bar x \Vert^2
                \\
                &\ge
                \frac{\gamma_0}{2}\left\Vert \bar x - x_0^+ \right\Vert^2
                - \frac{\gamma_0}{2}\Vert \bar x - x_0\Vert^2
                \\
                \iff
                F(x_0^+) - F(\bar x) &\le
                \frac{\gamma_0}{2}\Vert \bar x - x_0\Vert^2
                - \frac{\gamma_0}{2}\left\Vert \bar x - x_0^+ \right\Vert^2
                \\
                \iff
                \Upsilon_0(\bar x) &\le
                \frac{\gamma_0}{2}\Vert \bar x - x_0\Vert^2
                - \frac{\gamma_0}{2}\left\Vert \bar x - x_0^+ \right\Vert^2
                + \frac{\gamma_0 \alpha_0^2}{2}\Vert x_0^\dagger - \bar x\Vert^2,
            \end{align*}
            which establishes \eqref{eqn:apg-lyapunov-base}.
            \par
            Suppose additionally $x_{-1}^+ = x_{-1}^\dagger$ and $\alpha_0 = 1$.
            Then $x_0 = \left(\tfrac{\alpha_0 - q_0}{1 - q_0}\right)x_{-1}^\dagger + \left(\tfrac{1 - \alpha_0}{1 - q_0}\right)x_{-1}^+ = x_{-1}^\dagger = x_{-1}^+$ by Assumption \ref{ass:abst-nes-extrp}\ref{ass:abst-nes-extrp-lambdak} and Assumption \ref{ass:abst-nes-extrp}\ref{ass:abst-nes-extrp-x}.
            Also $x_0^\dagger = x_{-1}^+ + \alpha_0^{-1}(x_0^+ - x_{-1}^+) = x_0^+$ by Assumption \ref{ass:abst-nes-extrp}\ref{ass:abst-nes-extrp-xd} because $\alpha_0 = 1$.
            Substituting $\alpha_0 = 1$ and $x_0^\dagger = x_0^+$ into \eqref{eqn:apg-lyapunov-base}:
            \begin{align*}
                \Upsilon_0(\bar x)
                &\le
                \frac{\gamma_0}{2}\Vert \bar x - x_0\Vert^2
                - \frac{\gamma_0}{2}\Vert \bar x - x_0^+\Vert^2
                + \frac{\gamma_0}{2}\Vert x_0^+ - \bar x\Vert^2
                =
                \frac{\gamma_0}{2}\Vert \bar x - x_0\Vert^2,
            \end{align*}
            which establishes \eqref{eqn:apg-lyapunov-base-special}.
        \end{proof}
\section{Main Results}\label{sec:major-results}
    This section proves that the constant momentum algorithm described by \eqref{expr:const-momentum-recursion} converges linearly for all $\beta \in [0, 1)$.
    Unlike the ad hoc extension in the literature our rate depends on $\beta, L$ and $\mu$ and it does not require $\beta$ to be chosen in advance. 
    \par
    Section \ref{ssec:apg-const-momentum-alg} identifies the algorithm as a special case of the abstract Lyapunov framework of Section \ref{sec:aux-results}.
    Section \ref{ssec:lyp-ineqr-s-cnvx} strengthens the results of Section \ref{sec:aux-results} using strong convexity.
    Section \ref{ssec:numerical-illustrations} illustrates the resulting bound numerically.
    \par
    The proof of Theorem \ref{thm:lin-cnvg-constn-mmntm} has three stages.
    The first stage is Theorem \ref{thm:apg-momentum-equiv} in Section \ref{ssec:apg-const-momentum-alg}.
    It shows that the constant momentum algorithm is the APG algorithm with a specific set of constant parameters $(\alpha, q)$ that satisfies Assumption \ref{ass:const-param}.
    The second stage in Section \ref{ssec:lyp-ineqr-s-cnvx} gives the recurrence $\Upsilon_k(\bar x) \le M_\sigma(\alpha, q)\Upsilon_{k - 1}(\bar x)$ for these parameters under Assumption \ref{ass:strong-convexity}.
    The third stage is Propositions \ref{prop:mstar-beta-ge-betastar} and \ref{prop:mstar-beta-le-bstar} in Section \ref{ssec:optimal-contraction-factor}.
    They minimize $M_\sigma(\alpha, q)$ over all pairs realizing a given $\beta$.
    Proposition \ref{prop:mstar-beta-ge-betastar} evaluates the minimum for $\beta \ge \beta^\star$ and Proposition \ref{prop:mstar-beta-le-bstar} bounds it for $\beta \le \beta^\star$.
    Theorem \ref{thm:lin-cnvg-constn-mmntm} combines the three stages into a linear rate for every $\beta \in (0, 1)$.
    \par
    \todoinline{\noindent
        Figure \ref{fig:major-results-dependency} shows how the partial results of this section build up to Theorem \ref{thm:lin-cnvg-constn-mmntm}.
        An arrow from one result to another means that the proof of the second result uses the first.
        \begin{figure}[H]
            \centering
            \resizebox{0.98\textwidth}{!}{
            \begin{tikzpicture}[
                >=stealth, thick,
                res/.style={draw, rounded corners, align=center, text width=5.0cm, inner sep=4pt, font=\footnotesize, execute at begin node={\hyphenpenalty=10000 \relax}},
                up/.style={res, dashed},
                main/.style={res, very thick},
                asm/.style={res, densely dotted},
                fromasm/.style={->, densely dotted}
            ]
                \node[asm]  (a41)   at (2.2, 2.4)  {Assumption \ref{ass:const-param}\\ constant parameters};
                \node[asm]  (a23)   at (8.0, 2.4)  {Assumption \ref{ass:global-setup}\\ global assumptions};
                \node[up]   (amf)   at (-1.2, 0)   {Lemma \ref{lemma:abstract-momentum-form}\\ abstract momentum form};
                \node[up]   (lyap)  at (6, 0)      {Lemma \ref{lemma:apg-lyapunov-inequality}\\ APG Lyapunov inequality};
                \node[res]  (dbeta) at (12, 0)     {Lemma \ref{lemma:dbeta-equiv}\\ $(\alpha, q)$ for every $\beta$};
                \node[main] (equiv) at (-1.2, -2.6){Theorem \ref{thm:apg-momentum-equiv}\\ APG momentum equivalences};
                \node[res]  (lcst)  at (6, -2.6)   {Lemma \ref{lemma:apg-lya-const-param}\\ Lyapunov inequality under constant parameters};
                \node[asm, text width=3.6cm] (asc) at (-1.2, -5.2) {Assumption \ref{ass:strong-convexity}\\ strong convexity};
                \node[res]  (lmu)   at (6, -5.2)   {Lemma \ref{lemma:apg-lya-ineq-mu}\\ Lyapunov inequality under strong convexity};
                \node[res]  (min)   at (12, -2)    {Lemma \ref{lemma:m-in-alpha-beta}\\ $M_\sigma$ in terms of $\alpha, \beta$ only};
                \node[res]  (xi)    at (12, -4)    {Lemma \ref{lemma:xi-sigma-beta-nonpos}\\ when $\xi_{\sigma, \beta} \le 0$};
                \node[res]  (ge)    at (10.4, -6.8){Proposition \ref{prop:mstar-beta-ge-betastar}\\ $M_\sigma^*$ for $\beta \ge \beta^\star$};
                \node[res]  (le)    at (16.0, -6.8){Proposition \ref{prop:mstar-beta-le-bstar}\\ $M_\sigma^*$ bounds for $\beta \le \beta^\star$};
                \node[main] (thm)   at (6, -9.4)   {Theorem \ref{thm:lin-cnvg-constn-mmntm}\\ linear convergence for all $\beta \in (0, 1)$};
                \draw[->] (amf) -- (equiv);
                \draw[->] (lyap) -- (lcst);
                \draw[->] (lcst) -- (lmu);
                \draw[->] (dbeta) -- (min);
                \draw[->] (min) -- (xi);
                \draw[->] (xi) -- (ge);
                \draw[->] (xi) -- (le);
                \draw[->] (dbeta.east) -| ([xshift=9mm]le.east) -- ([yshift=2mm]le.east);
                \draw[->] (min.east) -| ([xshift=5mm]le.east) -- ([yshift=-2mm]le.east);
                \draw[->] (equiv) -- (thm);
                \draw[->] (lmu) -- (thm);
                \draw[->] (ge) -- (thm);
                \draw[->] (le) -- (thm);
                \draw[fromasm] (a41.south) -- (2.2, -1.3) -| ([xshift=15mm]equiv.north);
                \draw[fromasm] (a41.south) |- (lcst.west);
                \draw[fromasm] (a23.south) -- (lyap.north -| a23.south);
                \draw[fromasm] (asc.east) -- (lmu.west);
            \end{tikzpicture}
            }
            \caption{
                Dependencies among the results and assumptions of this section.
                Dotted boxes are assumptions.
                An arrow from $A$ to $B$ means that $B$ is derived using $A$, so $B$ shares the hypotheses of $A$.
            }
            \label{fig:major-results-dependency}
        \end{figure}
    }
    The following assumption gives a special case of Definition \ref{def:apg-algorithm}.
    \begin{assumption}[Constant parameters]\;\label{ass:const-param}\\
        Let $(\alpha_k)_{k \in \Z_+}, (q_k)_{k \in \Z_+}, (\gamma_k)_{k \in \Z_+}$ be given by Definition \ref{def:apg-algorithm}.
        We assume that there exist $\alpha$, $q$ and $\lambda$ such that:
        \begin{enumerate}[nosep]
            \item $\alpha \in (0, 1]$, $q \in [0, \alpha]\setminus \{1\}$.
            \item\label{ass:const-param:item1} $\gamma_k \equiv L$ for every $k \in \Z_+$ where $L$ is from Assumption \ref{ass:global-setup}\ref{ass:global-setup:item1}.
            \item\label{ass:const-param:item2} $\alpha_0 = 1$.
            \item\label{ass:const-param:item3} $\alpha_k = \alpha$ for all $k \in \N$ and $q_k = q$ for all $k \in \Z_+$.
            \item\label{ass:const-param:item4} $\lambda_k = \lambda = \tfrac{\alpha - q}{1 - q}$.
        \end{enumerate}
    \end{assumption}

    \subsection{APG and the constant momentum algorithm}\label{ssec:apg-const-momentum-alg}
        This subsection links the constant momentum algorithm in Definition \ref{def:const-momentum-alg} to a family of algorithms satisfying Definition \ref{def:apg-algorithm}.
        This result is presented in Theorem \ref{thm:apg-momentum-equiv} for the purpose of applying \eqref{eqn:apg-lyapunov-main} to Definition \ref{def:apg-algorithm}.
        \begin{definition}[The constant momentum algorithm]\;\label{def:const-momentum-alg}\\
            Let $x_{-1}^+, x_0^+ \in \RR^n$ be arbitrary and let $\beta \in [0, 1)$ be constant.
            Suppose that $(F, f, g, L, \mu)$ satisfy Assumption \ref{ass:global-setup}.
            The constant momentum algorithm generates two sequences of iterates $(x_k)_{k \in \N}, (x_k^+)_{k \in \N}$ such that for all $k \in \Z_+$:
            \begin{align}
                &x_{k + 1} = x_k^+ + \beta\left(x_k^+ - x_{k - 1}^+\right),
                \label{def:const-momentum-alg-x}
                \\
                &x_{k + 1}^+ = T_L^{(f, g)}(x_{k + 1}).
                \label{def:const-momentum-alg-xp}
            \end{align}
        \end{definition}
        Our first result is the following theorem which makes the correspondence between the APG algorithm and the constant momentum algorithm precise.
        \begin{theorem}[APG momentum equivalences]\;\label{thm:apg-momentum-equiv}\\
            Let $(x_k)_{k \in \Z_+}, (x_k^+)_{k \in \Z_+}, (x_k^\dagger)_{k \in \Z_+}$ satisfy Definition \ref{def:apg-algorithm}, and $(y_k)_{k \in \N}, (y_k^+)_{k \in \N}$ satisfy Definition \ref{def:const-momentum-alg}. 
            Let $\hat x\in \RR^n$ and assume that the following hold: 
            \begin{enumerate}[nosep]
                \item\label{thm:apg-momentum-equiv-hyp1} the initial conditions of Definition \ref{def:apg-algorithm} satisfy $x_{-1}^+ = x_{-1}^\dagger =: \hat x$ and the initial conditions of Definition \ref{def:const-momentum-alg} satisfy $y_{-1}^+ = y_0^+ := T_L^{(f,g)}(\hat x)$;
                \item\label{thm:apg-momentum-equiv-hyp2} the sequences $(\alpha_k)_{k \in \Z_+}, (q_k)_{k \in \Z_+}, (\gamma_k)_{k \in \Z_+}$ satisfy Assumption \ref{ass:const-param} with constant pair $(\alpha, q)$ and the constant $\beta$ of Definition \ref{def:const-momentum-alg} satisfies $\beta =\tfrac{\alpha - q}{1 - q}\cdot\tfrac{1 - \alpha}{\alpha}$.
            \end{enumerate}
            Then the iterates produced by these two algorithms are the same. That is, $y_k = x_k$ and $y_k^+ = x_k^+$ for all $k \in \N$.
        \end{theorem}
        \begin{proof}
            \todoinline{
                Theorem \ref{thm:apg-momentum-equiv} proof depends on:
                \begin{itemize}[nosep]
                    \item Definition \ref{def:apg-algorithm}.
                    \item Assumption \ref{ass:abst-nes-extrp}.
                    \item Assumption \ref{ass:abst-nes-extrp}\ref{ass:abst-nes-extrp-lambdak}.
                    \item Assumption \ref{ass:abst-nes-extrp}\ref{ass:abst-nes-extrp-x}.
                    \item Lemma \ref{lemma:abstract-momentum-form}.
                    \item Assumption \ref{ass:const-param}\ref{ass:const-param:item1}.
                    \item Assumption \ref{ass:const-param}\ref{ass:const-param:item3}.
                    \item Definition \ref{def:const-momentum-alg}.
                    \item \eqref{def:const-momentum-alg-x}.
                    \item \eqref{def:const-momentum-alg-xp}.
                    \item \eqref{def:apg-algorithm-xp}.
                \end{itemize}
            }
            The sequences $(x_k)_{k \in \Z_+}, (x_k^+)_{k \in \Z_+}, (x_k^\dagger)_{k \in \Z_+}$ satisfy Assumption \ref{ass:abst-nes-extrp} by Definition \ref{def:apg-algorithm} so Lemma \ref{lemma:abstract-momentum-form} applies to them.
            At $k = 0$ Assumption \ref{ass:const-param} gives $\alpha_0 = 1$ and $\gamma_0 = L$.
            Hence $\lambda_0 = \tfrac{\alpha_0 - q_0}{1 - q_0} = 1$ by Assumption \ref{ass:abst-nes-extrp}\ref{ass:abst-nes-extrp-lambdak}.
            It follows from Assumption \ref{ass:abst-nes-extrp}\ref{ass:abst-nes-extrp-x} and \eqref{def:apg-algorithm-xp} that:
            \begin{align*}
                x_0
                &= \lambda_0 x_{-1}^\dagger + (1 - \lambda_0)x_{-1}^+
                = x_{-1}^\dagger
                \underset{\text{\ref{thm:apg-momentum-equiv-hyp1}}} = \hat x
                \implies
                x_0^+ = T_{\gamma_0}^{(f,g)}(x_0) = T_L^{(f,g)}(\hat x) 
                \underset{\text{\ref{thm:apg-momentum-equiv-hyp1}}}= y_0^+.
            \end{align*}
            Using Lemma \ref{lemma:abstract-momentum-form} and Assumption \ref{ass:const-param} yields for all $k \in \Z_+$:
            \begin{align*}
                x_{k + 1} 
                &= 
                x_k^+ + \lambda_{k + 1}\tfrac{1 - \alpha_k}{\alpha_k}
                \left(x_k^+ - x_{k - 1}^+\right)
                \\
                &=
                \begin{cases}
                    x_k^+ + 
                    \underbrace{
                        \left(
                            \tfrac{\alpha - q}{1 - q}\tfrac{1 - \alpha}{\alpha}
                        \right)
                    }_{=\beta}
                    \left(x_k^+ - x_{k - 1}^+\right)
                    & k \in \N, 
                    \\
                    x_0^+  & k = 0,
                \end{cases}
                \\
                &\underset{\text{\ref{thm:apg-momentum-equiv-hyp2}}}=
                \begin{cases}
                    x_k^+ + \beta\left(x_k^+ - x_{k - 1}^+\right) & k \in \N, 
                    \\
                    x_0^+  & k = 0.
                \end{cases}
            \end{align*}
            The case $k \in \N$ holds because $\alpha_k = \alpha_{k + 1} = \alpha$ and $q_{k + 1} = q$ by Assumption \ref{ass:const-param}\ref{ass:const-param:item3}.
            The case $k = 0$ holds because $\alpha_0 = 1$ makes $\tfrac{1 - \alpha_0}{\alpha_0} = 0$.
            Therefore we established:
            \begin{align}
                x_{k + 1} =
                \begin{cases}
                    x_k^+ + \beta\left(x_k^+ - x_{k - 1}^+\right) & k \in \N, 
                    \\
                    x_0^+  & k = 0,
                \end{cases}
            \label{eqn:apg-momentum-equiv-recursion}
            \end{align}
            where $\beta = \tfrac{\alpha - q}{1 - q}\cdot\frac{1 - \alpha}{\alpha}$.
            \par
            Next we show $y_k = x_k$ and $y_k^+ = x_k^+$ for every $k \in \N$ by strong induction on $k$.
            For every $k \in \N$ denote by $\mathbb H_k$ the statement $y_{k - 1}^+ = x_{k - 1}^+$ and $y_k^+ = x_k^+$.
            Consider the base case $k = 1$:
            \begin{align*}
                &y_1
                \underset{(1)}{=} y_0^+ + \beta(y_0^+ - y_{-1}^+)
                \underset{\text{\ref{thm:apg-momentum-equiv-hyp1}}} = y_0^+
                \underset{\text{\ref{thm:apg-momentum-equiv-hyp1}}} = T_L^{(f, g)}(\hat x)
                \underset{\text{\eqref{eqn:apg-momentum-equiv-recursion}}}= x_1
                \\
                \implies
                &y_1^+
                \underset{(2)}{=} T_L^{(f, g)}(y_1)
                = T_L^{(f, g)}(x_1)
                \underset{(3)}{=} x_1^+.
            \end{align*}
            At (1) we use \eqref{def:const-momentum-alg-x} from Definition \ref{def:const-momentum-alg}.
            At (2) we use \eqref{def:const-momentum-alg-xp}.
            At (3) we use \eqref{def:apg-algorithm-xp} with $\gamma_1 = L$ from Assumption \ref{ass:const-param}\ref{ass:const-param:item1}.
            The last equality in the first row is \eqref{eqn:apg-momentum-equiv-recursion} at $k = 0$ together with $x_0^+ = T_L^{(f, g)}(\hat x)$.
            Since $y_0^+ = T_L^{(f, g)}(\hat x) = x_0^+$ by \ref{thm:apg-momentum-equiv-hyp1} this establishes the base case:
            \begin{align*}
                \mathbb H_1: y_0^+ = x_0^+, \quad y_1^+ = x_1^+.
            \end{align*}
            Consider the inductive step.
            Fix $k \ge 1$ and assume $\mathbb H_k$.
            The following holds:
            \begin{align*}
                &y_{k + 1}
                \underset{(4)}{=} y_k^+ + \beta(y_k^+ - y_{k - 1}^+)
                \underset{(\mathbb H_k)}{=} x_k^+ + \beta(x_k^+ - x_{k - 1}^+)
                \underset{(5)}{=} x_{k + 1}
                \\
                \implies 
                &y_{k + 1}^+
                \underset{(6)}{=} T_L^{(f, g)}(y_{k + 1})
                = T_L^{(f, g)}(x_{k + 1})
                \underset{(7)}{=} x_{k + 1}^+.
            \end{align*}
            At (4) and (6) we use \eqref{def:const-momentum-alg-x} and \eqref{def:const-momentum-alg-xp} in Definition \ref{def:const-momentum-alg} respectively.
            At (5) we use \eqref{eqn:apg-momentum-equiv-recursion} at $k$.
            At (7) we use \eqref{def:apg-algorithm-xp} with $\gamma_{k + 1} = L$ from Assumption \ref{ass:const-param}\ref{ass:const-param:item1}.
            Together with $y_k^+ = x_k^+$ from $\mathbb H_k$ this establishes $\mathbb H_{k + 1}$ and $y_{k + 1} = x_{k + 1}$.
            This completes the induction.
        \end{proof}
        \begin{remark}
            The momentum coefficient $\beta := \tfrac{\alpha - q}{1 - q}\cdot\tfrac{1 - \alpha}{\alpha}$ depends only on the constant pair $(\alpha, q)$ of Assumption \ref{ass:const-param}\ref{ass:const-param:item3}.
        \end{remark}
        To make the above equivalence useful, we continue to strengthen the Lyapunov inequality of Lemma \ref{lemma:apg-lyapunov-inequality} under strong convexity.
        Under Assumption \ref{ass:const-param}, we express $C_k, \xi_k, M_k$ defined in Lemma \ref{lemma:apg-lya-const-param} as functions of $\alpha, q$ parameterized by
        \begin{align}
            \label{def:sigma}
            \sigma := \frac{\mu}{L}.
        \end{align}
        First, we define these functions.
        For all $\alpha \in (0, 1]$, $q \in [0, \alpha] \setminus \{1\}$ we define:
        \begin{align}
            \label{def:c-sigma-alpha-q}
            \forall \sigma \in [0, 1]& \quad 
            C_\sigma(\alpha, q) := \tfrac{(1 - \sigma)(\alpha - q)}{\alpha(1 - q)},
            \\
            \label{def:xi-sigma-alpha-q}
            \forall \sigma \in (0, 1]& \quad 
            \xi_\sigma(\alpha, q) := \tfrac{\alpha(1 - \alpha)(q - \sigma)}{\sigma(1 - q)},
            \\
            \label{def:m-sigma-alpha-q}
            \forall \sigma \in (0, 1]& \quad 
            M_\sigma(\alpha, q) :=
            \max\left(
                1 - \alpha, C_\sigma(\alpha, q)
            \right)
            + \left[
                \xi_\sigma(\alpha, q)
            \right]_+.
        \end{align}
        Second, the following lemma writes the parameters $C_k, \xi_k, M_k$ in terms of these functions using Assumption \ref{ass:const-param} to adapt the results of Lemma \ref{lemma:apg-lyapunov-inequality}.
        \begin{lemma}[APG Lyapunov inequality under constant parameters]\;\label{lemma:apg-lya-const-param}\\
            Assume the same hypotheses as Lemma \ref{lemma:apg-lyapunov-inequality} and Assumption \ref{ass:const-param}.
            Let $\sigma$ be given by \eqref{def:sigma}.
            Let $\Xi_k(\bar x)$ be given by \eqref{def:xi-k}, $C_\sigma(\alpha, q)$ be given by \eqref{def:c-sigma-alpha-q}.
            Then for all minimizers $\bar x$ of $F$ the following statements are true.
            \begin{enumerate}[nosep]
                \item For all $k \in \N$ the following inequality holds:
                \begin{align}\label{lemma:apg-lya-const-param-eqn1}
                    \Upsilon_k(\bar x) &\le
                    \max(1 - \alpha, C_\sigma(\alpha, q))\Upsilon_{k - 1}(\bar x) + \Xi_k(\bar x),
                    \\
                    \Xi_k(\bar x) &= \frac{L\alpha(1 - \alpha)(q - \sigma)}{2(1 - q)}\left\Vert x_{k - 1}^+ - \bar x \right\Vert^2.
                \end{align}
                \item If in addition the initial conditions satisfy $x_{-1}^+ = x_{-1}^\dagger$ then $x_0 = x_{-1}^+ = x_{-1}^\dagger$ and we have the base case
                \begin{align}\label{lemma:apg-lya-const-param-eqn2}
                    \Upsilon_0(\bar x) \le \frac{L}{2}\Vert \bar x - x_0\Vert^2.
                \end{align}
            \end{enumerate}
        \end{lemma}
        \begin{proof}
            \todoinline{
                Lemma \ref{lemma:apg-lya-const-param} proof depends on:
                \begin{itemize}[nosep]
                    \item Lemma \ref{lemma:apg-lyapunov-inequality}.
                    \item Definition \ref{def:apg-algorithm}.
                    \item Assumption \ref{ass:abst-nes-extrp}\ref{ass:abst-nes-extrp-alphaq}.
                    \item Assumption \ref{ass:abst-nes-extrp}\ref{ass:abst-nes-extrp-lambdak}.
                    \item Assumption \ref{ass:const-param}.
                    \item Assumption \ref{ass:global-setup}.
                    \item \eqref{eqn:apg-lyapunov-main}.
                    \item \eqref{eqn:apg-lyapunov-base-special}.
                    \item \eqref{def:delta}.
                    \item \eqref{def:upsilon}.
                    \item \eqref{def:ck}.
                    \item \eqref{def:xi-k}.
                    \item \eqref{def:sigma}.
                    \item The proof implicitly depends on: \eqref{def:c-sigma-alpha-q}.
                \end{itemize}
            }
            Assumption \ref{ass:const-param} gives $\alpha_k = \alpha, q_k = q, \gamma_k = L$ for all $k \in \N$, with the base case $\alpha_0 = 1$, $q_0 = q$ and $\gamma_0 = L$.
            Assumption \ref{ass:global-setup} gives $\tfrac{\mu}{2}\Vert x - y\Vert^2 \le D_f(x, y) \le \tfrac{L}{2}\Vert x - y\Vert^2$ for all $x, y \in \RR^n$ so $0\le \mu \le L$ and therefore $\sigma \in [0, 1]$ by \eqref{def:sigma}.
            Recall $C_k = \tfrac{(\gamma_k - \mu)\alpha_k\lambda_k}{\gamma_{k - 1}\alpha_{k - 1}^2}$ from \eqref{def:ck} and $\lambda_k = \tfrac{\alpha_k - q_k}{1 - q_k}$ from Assumption \ref{ass:abst-nes-extrp}\ref{ass:abst-nes-extrp-lambdak}.
            For all $k \ge 2$ substituting $\lambda_k = \tfrac{\alpha - q}{1 - q}, \gamma_k = \gamma_{k - 1} = L, \alpha_k = \alpha_{k - 1} = \alpha$ and $q_k = q$ into $C_k$ and using $\tfrac{L - \mu}{L} = 1 - \sigma$ from \eqref{def:sigma} yields:
            \begin{align*}
                C_k 
                =
                \frac{(L - \mu)\alpha}{L\alpha^2}
                \cdot\frac{\alpha - q}{1 - q}
                = 
                \frac{(1 - \sigma)}{\alpha}
                \cdot\frac{\alpha - q}{1 - q}
                = C_\sigma(\alpha, q).
            \end{align*}
            We consider $C_1$ separately:
            \begin{align*}
                C_1 = \frac{\alpha(L - \mu)}{L} \cdot\frac{\alpha - q}{1 - q}
                = \frac{\alpha^2(1 - \sigma)}{\alpha} \cdot\frac{\alpha - q}{1 - q}
                = \alpha^2 C_\sigma(\alpha, q).
            \end{align*}
            Moreover $C_\sigma(\alpha, q) \ge 0$ because $\sigma \le 1$, $\alpha \ge q$ and $1 - q > 0$.
            Since $\alpha \in (0, 1]$, it follows that for all $k \ge 1$:
            \begin{align}\label{lemma:apg-lya-const-param-pitem1}
                C_1 \le C_k \le C_\sigma(\alpha, q).
            \end{align}
            Recall that $\Xi_k(\bar x) = \tfrac{\alpha_k(1 - \lambda_k)(\gamma_kq_k - \mu)}{2}\Vert x_{k - 1}^+ - \bar x\Vert^2$ by \eqref{def:xi-k}.
            Under Assumption \ref{ass:const-param}, for all $k \in \N$ we have $\alpha_k = \alpha, q_k = q, \gamma_k = L$ and $\lambda_k = \lambda$.
            Substituting them into the fraction $\tfrac{\alpha_k(1 - \lambda_k)(\gamma_kq_k - \mu)}{2}$ and using $\mu = \sigma L$ from \eqref{def:sigma} yields:
            \begin{align}
                \frac{\alpha(1 - \lambda)(Lq - \mu)}{2} 
                = \frac{\alpha(1 - \alpha)(Lq - \mu)}{2(1 - q)}
                = \frac{L\alpha(1 - \alpha)(q - \sigma)}{2(1 - q)}.
            \end{align}
            Therefore for all $k \in \N$ we obtain the following expression for $\Xi_k(\bar x)$: 
            \begin{align}\label{lemma:apg-lya-const-param-pitem2}
                \Xi_k(\bar x) &= \frac{L\alpha(1 - \alpha)(q - \sigma)}{2(1 - q)}
                \left\Vert
                    x_{k - 1}^+ - \bar x
                \right\Vert^2.
            \end{align}
            Using \eqref{lemma:apg-lya-const-param-pitem1} and \eqref{eqn:apg-lyapunov-main} from Lemma \ref{lemma:apg-lyapunov-inequality} yields for all $k \in \N$ the inequality:
            \begin{align*}
                \Upsilon_k(\bar x)
                &\le \max(1 - \alpha_k, C_k)\Upsilon_{k - 1}(\bar x) + \Xi_k(\bar x)
                \\
                &\le \max(1 - \alpha, C_\sigma(\alpha, q))\Upsilon_{k - 1}(\bar x) + \Xi_k(\bar x).
            \end{align*}
            The second inequality applies \eqref{lemma:apg-lya-const-param-pitem1} to the coefficient of $\Upsilon_{k - 1}(\bar x)$.
            The latter step needs $\Upsilon_{k - 1}(\bar x) \ge 0$.
            This holds because $\bar x$ minimizes $F$ so $\Delta_{k - 1}(\bar x) \ge 0$ by \eqref{def:delta} and \eqref{def:upsilon} adds a nonnegative term.
            Therefore, \eqref{lemma:apg-lya-const-param-eqn1} is verified.
            \par
            Finally, any sequence satisfying Assumption \ref{ass:const-param} automatically satisfies the constraints of $(\alpha_k)_{k \in \Z_+}$ and $(q_k)_{k \in \Z_+}$ in Assumption \ref{ass:abst-nes-extrp}\ref{ass:abst-nes-extrp-alphaq} that Definition \ref{def:apg-algorithm} requires.
            Assumption \ref{ass:const-param} also fixes $\alpha_0 = 1$.
            Therefore if $x_{-1}^+ = x_{-1}^\dagger$ then Lemma \ref{lemma:apg-lyapunov-inequality} gives $x_0 = x_{-1}^+ = x_{-1}^\dagger$ and \eqref{eqn:apg-lyapunov-base-special} reads $\Upsilon_0(\bar x) \le \tfrac{\gamma_0}{2}\Vert \bar x - x_0\Vert^2 = \tfrac{L}{2}\Vert \bar x - x_0\Vert^2$ because $\gamma_0 = L$.
        \end{proof}

    \subsection{The Lyapunov inequality under strong convexity}\label{ssec:lyp-ineqr-s-cnvx}
        The following definition introduces a condition weaker than strong convexity.
        \begin{definition}[Quadratic growth condition]\label{def:q-growth}
            Suppose that $F : \RR^n \to \overline\RR$ is proper and closed with set of minimizers $X$.
            Then $F$ satisfies the quadratic growth condition if there exists $\theta \ge 0$ such that for every $x \in \RR^n$ there exists $\bar x \in X$ satisfying:
            \begin{align*}
                F(x) - F(\bar x) &\ge \frac{\theta}{2}\Vert x - \bar x\Vert^2. 
            \end{align*}
        \end{definition}
        If $F$ satisfies Assumption \ref{ass:global-setup} then it satisfies Definition \ref{def:q-growth} with $\theta = \mu$.
        The lemma below builds on top of Lemma \ref{lemma:apg-lyapunov-inequality} using strong convexity and Assumption \ref{ass:const-param}.
        \begin{assumption}[Strong convexity]\;\label{ass:strong-convexity}\\
            Let $(F, f, g, L, \mu)$ satisfy Assumption \ref{ass:global-setup} and let $\sigma$ be given by \eqref{def:sigma}.
            Assume $(f, \mu)$ satisfy the following.
            \begin{enumerate}[nosep]
                \item\label{ass:strong-convexity:item1} $f$ is strongly convex with modulus $\mu > 0$ as defined in Definition \ref{def:strong-convexity}.
            \end{enumerate}
        \end{assumption}
        \begin{remark}
            Assumption \ref{ass:global-setup} gives $\tfrac{\mu}{2}\Vert x - y\Vert^2 \le D_f(x, y) \le \tfrac{L}{2}\Vert x - y\Vert^2$ for all $x, y \in \RR^n$ so $\mu \le L$.
            Hence Assumption \ref{ass:strong-convexity}\ref{ass:strong-convexity:item1} gives $\sigma \in (0, 1]$.
        \end{remark}
        \begin{lemma}[APG Lyapunov inequality under strong convexity]\;
        \label{lemma:apg-lya-ineq-mu}\\
            Assume the same hypotheses as Lemma \ref{lemma:apg-lya-const-param} and Assumption \ref{ass:strong-convexity}.
            Suppose that for all $k \in \Z_+$, $\Delta_k, \Upsilon_k$ are given by \eqref{def:delta}, \eqref{def:upsilon} respectively.
            Let $M_\sigma(\alpha, q)$ be given by \eqref{def:m-sigma-alpha-q}.
            Then the following holds for all $k \in \N$ and the unique minimizer $\bar x \in \RR^n$ of $F$:
            \begin{align}\label{eqn:apg-lya-ineq-mu-rslt}
                \Upsilon_k(\bar x) &\le M_\sigma(\alpha, q) \Upsilon_{k - 1}(\bar x).
            \end{align}
        \end{lemma}
        \begin{proof}
            \todoinline{
                Lemma \ref{lemma:apg-lya-ineq-mu} proof depends on:
                \begin{itemize}[nosep]
                    \item Lemma \ref{lemma:apg-lya-const-param}.
                    \item Definition \ref{def:q-growth}.
                    \item Assumption \ref{ass:global-setup}.
                    \item \eqref{def:delta}.
                    \item \eqref{def:upsilon}.
                    \item \eqref{def:xi-k}.
                    \item \eqref{def:sigma}.
                    \item \eqref{def:xi-sigma-alpha-q}.
                    \item \eqref{def:m-sigma-alpha-q}.
                \end{itemize}
                The proof implicitly depends on:
                \begin{itemize}[nosep]
                    \item Assumption \ref{ass:const-param}.
                    \item Assumption \ref{ass:strong-convexity}.
                \end{itemize}
            }
            Fix any $k \in \N$ and let $\Xi_k(\bar x)$ be given by \eqref{def:xi-k}.
            Lemma \ref{lemma:apg-lya-const-param} evaluates it under Assumption \ref{ass:const-param}.
            Let $\xi_\sigma(\alpha, q)$ be given by \eqref{def:xi-sigma-alpha-q}.
            Then, using $\mu = \sigma L$ from \eqref{def:sigma}, the following holds:
            \begin{align*}
                \Xi_k(\bar x)
                &= \tfrac{\mu}{2}\xi_\sigma(\alpha, q)\Vert x_{k - 1}^+ - \bar x\Vert^2
                \\
                &\underset{(1)}\le
                [\xi_\sigma(\alpha, q)]_+\tfrac{\mu}{2}\Vert x_{k - 1}^+ - \bar x\Vert^2
                \underset{(2)}\le
                [\xi_\sigma(\alpha, q)]_+\Delta_{k - 1}(\bar x)
                \underset{(3)}\le
                [\xi_\sigma(\alpha, q)]_+\Upsilon_{k - 1}(\bar x).
            \end{align*}
            At (1) the sign of $\xi_\sigma(\alpha, q)$ is unknown so we replace it by its positive part to obtain an upper bound.
            At (2) we invoke Definition \ref{def:q-growth} with $\theta = \mu$, which holds for $F$ by Assumption \ref{ass:global-setup}.
            Together with $\Delta_{k - 1}(\bar x) = F(x_{k - 1}^+) - F(\bar x)$ from \eqref{def:delta} it gives $\tfrac{\mu}{2}\Vert x_{k - 1}^+ - \bar x\Vert^2 \le \Delta_{k - 1}(\bar x)$.
            This upper bound applies because $[\xi_\sigma(\alpha, q)]_+ \ge 0$.
            At (3) we use \eqref{def:upsilon} and combine it with Assumption \ref{ass:const-param} which gives $\Upsilon_{k - 1}(\bar x) = \Delta_{k - 1}(\bar x) + \tfrac{L\alpha^2}{2}\Vert x_{k - 1}^\dagger - \bar x\Vert^2 \ge \Delta_{k - 1}(\bar x)$ and $[\xi_\sigma(\alpha, q)]_+ \ge 0$.
            Lemma \ref{lemma:apg-lya-const-param} and the bound on $\Xi_k(\bar x)$ above give:
            \begin{align*}
                \Upsilon_k(\bar x) 
                &\le
                \max(1 - \alpha, C_\sigma(\alpha, q))\Upsilon_{k - 1}(\bar x) + \Xi_k(\bar x)
                \\
                &\le
                \left(
                    \max(1 - \alpha, C_\sigma(\alpha, q)) + [\xi_\sigma(\alpha, q)]_+
                \right)\Upsilon_{k - 1}(\bar x).
            \end{align*}
            The coefficient of $\Upsilon_{k - 1}(\bar x)$ is $M_\sigma(\alpha, q)$ by \eqref{def:m-sigma-alpha-q}.
            This establishes \eqref{eqn:apg-lya-ineq-mu-rslt}.
        \end{proof}
    
    \subsection{Characterizing the linear convergence of constant momentum algorithm}\label{ssec:optimal-contraction-factor}
        Assumption \ref{ass:const-param} continues to hold in this subsection and we suppose that $\sigma$ given by \eqref{def:sigma} satisfies $\sigma \in (0, 1)$.
        The results from Theorem \ref{thm:apg-momentum-equiv} are therefore available.
        We then define the following set which characterizes the parameters $\alpha, q, \beta$ for which Definitions \ref{def:apg-algorithm} and \ref{def:const-momentum-alg} produce the same iterates:
        \begin{align}\label{def:dbeta-set}
            \mathcal D_\beta:= 
            \left\{
                (\alpha, q): \alpha \in (0, 1],\; 
                q \in [0, \alpha]\setminus \{1\},\;
                \beta = \left(
                    \tfrac{\alpha - q}{1 - q}
                \right)
                \tfrac{1 - \alpha}{\alpha}
            \right\}.
        \end{align}
        Under Assumption \ref{ass:const-param} and the equivalence between constant momentum and abstract APG we minimize $M_\sigma$ from Lemma \ref{lemma:apg-lya-ineq-mu} over every $(\alpha, q) \in \mathcal D_\beta$.
        That infimum is given by:
        \begin{align}\label{def:mstar-sigma-beta-bnd}
            M^*_\sigma(\beta) &:= \inf_{(\alpha, q) \in \mathcal D_\beta} M_\sigma(\alpha, q). 
        \end{align}
        This is an infimum because it is not yet clear if the infimum is attained by an element of $\mathcal D_\beta$.
        This subsection derives an algebraic form of an upper bound of $M_\sigma^*(\beta)$ for all $\sigma \in (0, 1)$ and $\beta \in (0, 1)$.
        To achieve that we present the following three lemmas which establish the necessary intermediate results with each proof given in the appendix.
        \par
        The first lemma characterizes an equivalent expression for $\mathcal D_\beta$ defined in \eqref{def:dbeta-set}.
        \begin{lemma}[$(\alpha, q)$ for every $\beta$]\label{lemma:dbeta-equiv}
            Let $\mathcal D_\beta$ be given by \eqref{def:dbeta-set}.
            Then it holds that for all $\beta \in (0, 1)$
            \begin{align*}
                \mathcal D_\beta =
                \left\lbrace
                    (\alpha, q): 
                    q = \tfrac{\alpha(1 - \alpha - \beta)}{1 - \alpha(1 + \beta)}, \;
                    \alpha \in (0, 1 - \beta]
                \right\rbrace. 
            \end{align*}
            When $\beta = 0$, $\mathcal D_0 = \{(\alpha, \alpha) : \alpha \in (0, 1)\} \cup (\{1\}\times [0, 1))$.
        \end{lemma}
        \begin{proof}
            See Section \ref{ssec:dbeta-equiv-pf}.
        \end{proof}
        The second lemma below expresses \eqref{def:c-sigma-alpha-q}, \eqref{def:xi-sigma-alpha-q} and \eqref{def:m-sigma-alpha-q} in terms of only $\alpha, \beta$ for all $(\alpha, q) \in \mathcal D_\beta$ assuming $\beta \in (0, 1)$ and $\sigma \in (0, 1)$.
        In this way $q$ is eliminated and $C_\sigma, \xi_\sigma$ and $M_\sigma$ become simpler.
        \begin{lemma}[$M_\sigma$ in terms of $\alpha, \beta$ only]\label{lemma:m-in-alpha-beta}
            For all $\beta\in (0, 1)$ and $\sigma \in (0, 1)$ the following are true for all $(\alpha, q) \in \mathcal D_\beta$:
            \begin{enumerate}[nosep]
                \item $C_\sigma(\alpha, q) = \tfrac{\beta(1 - \sigma)}{1 - \alpha}$,
                \item $\xi_\sigma(\alpha, q) = \tfrac{\alpha}{\sigma}\left((\alpha - \sigma) - \alpha\cdot \tfrac{\beta(1 - \sigma)}{1 - \alpha}\right)$. 
            \end{enumerate}
        \end{lemma}
        \begin{proof}
            See Section \ref{ssec:m-in-alpha-beta-pf}.
        \end{proof}
        From the results of Lemma \ref{lemma:m-in-alpha-beta} we define for all $\beta \in (0, 1)$, $\sigma \in (0, 1)$ and $\alpha\in (0, 1 - \beta]$ the following for short:
        \begin{align}
            C_{\sigma, \beta}(\alpha) 
            &:= \tfrac{\beta(1 - \sigma)}{1 - \alpha}, 
            \label{def:c-sigma-beta-in-alpha}
            \\
            \xi_{\sigma, \beta}(\alpha)
            &:= \tfrac{\alpha}{\sigma}\left((\alpha - \sigma) - \alpha\cdot C_{\sigma, \beta}(\alpha)\right),
            \label{def:xi-sigma-beta-in-alpha}
            \\
            M_{\sigma, \beta}(\alpha) &:=
            \max(1 - \alpha, C_{\sigma, \beta}(\alpha)) + [\xi_{\sigma, \beta}(\alpha)]_+. 
            \label{def:m-sigma-beta-in-alpha}
        \end{align}
        By the results of Lemma \ref{lemma:m-in-alpha-beta} the following holds for all $(\alpha, q) \in \mathcal D_\beta$:
        \begin{align*}
            C_{\sigma, \beta}(\alpha) = C_\sigma(\alpha, q),
            \quad
            \xi_{\sigma, \beta}(\alpha) = \xi_\sigma(\alpha, q),
            \quad 
            M_{\sigma, \beta}(\alpha) = M_\sigma(\alpha, q). 
        \end{align*}
        Consequently $M_\sigma^*$ from \eqref{def:mstar-sigma-beta-bnd} is equivalent to:
        \begin{align}\label{eqn:mstar-with-alpha-only}
            M_\sigma^*(\beta) &:= \inf_{(\alpha, q)\in \mathcal D_\beta} M_\sigma(\alpha, q)
            = \inf_{\alpha \in (0, 1 - \beta]} M_{\sigma, \beta}(\alpha).
        \end{align}
        The third lemma below characterizes a condition on $\beta$ such that $\xi_{\sigma, \beta}(\alpha) \le 0$.
        The goal here is to characterize the conditions under which evaluating $M_\sigma^*$ is simpler.
        \begin{lemma}[Characterizing when $\xi_{\sigma, \beta}(\alpha) \le 0$]\label{lemma:xi-sigma-beta-nonpos}
            Fix any $\sigma \in (0, 1)$. We write $\beta^\star := \frac{1 - \sqrt{\sigma}}{1 + \sqrt{\sigma}}$.
            Let $\xi_{\sigma, \beta}$ be given by \eqref{def:xi-sigma-beta-in-alpha}. 
            Then the following holds for all $\alpha \in (0, 1 - \beta]$:
            \begin{align*}
                \xi_{\sigma, \beta}(\alpha) \le 0 \iff 
                1 > \beta \ge \beta^\star. 
            \end{align*}
        \end{lemma}
        \begin{proof}
            See Section \ref{ssec:xi-sigma-beta-nonpos-pf}.
        \end{proof}
        We are now in a position to introduce two propositions that give an upper bound of $M_\sigma^*$ and show in addition that the infimum is attained.
        By Lemma \ref{lemma:xi-sigma-beta-nonpos}, $\xi_{\sigma, \beta}(\alpha) \le 0$ for all $\alpha \in (0, 1 - \beta]$ whenever $\beta \ge \beta^\star$ which reduces $M_{\sigma, \beta}(\alpha)$ from \eqref{def:m-sigma-beta-in-alpha} to $\max(1 - \alpha, C_{\sigma, \beta}(\alpha))$ on this range.
        The proposition below evaluates $M_\sigma^*(\beta)$ for all $\beta \in [\beta^\star, 1)$ using this reduction.
        \begin{proposition}[The value of $M_\sigma^*$ for all $\beta \ge \beta^\star$]\label{prop:mstar-beta-ge-betastar}
            Fix any $\sigma \in (0, 1)$ and write $\beta^\star := \tfrac{1 - \sqrt\sigma}{1 + \sqrt\sigma}$.
            Then for all $\beta \in [\beta^\star, 1)$ the following holds:
            \begin{align*}
                M_\sigma^*(\beta) =
                \begin{cases}
                    \sqrt{(1 - \sigma)\beta}, & \beta \in [\beta^\star, 1 - \sigma],
                    \\
                    \beta, & \beta \in (1 - \sigma, 1).
                \end{cases}
            \end{align*}
            Moreover the infimum defining $M_\sigma^*(\beta)$ in \eqref{eqn:mstar-with-alpha-only} is attained: at $\alpha^\star = 1 - \sqrt{(1 - \sigma)\beta}$ in the first case and at $\alpha^\star = 1 - \beta$ in the second case.
        \end{proposition}
        \begin{proof}
            \todoinline{
                Proposition \ref{prop:mstar-beta-ge-betastar} proof depends on:
                \begin{itemize}[nosep]
                    \item Lemma \ref{lemma:xi-sigma-beta-nonpos}.
                    \item \eqref{def:m-sigma-beta-in-alpha}.
                    \item \eqref{def:c-sigma-beta-in-alpha}.
                    \item \eqref{eqn:mstar-with-alpha-only}.
                \end{itemize}
            }
            Fix any $\beta \in [\beta^\star, 1)$.
            By Lemma \ref{lemma:xi-sigma-beta-nonpos} and \eqref{def:m-sigma-beta-in-alpha}, $M_{\sigma, \beta}(\alpha) = \max(1 - \alpha, C_{\sigma, \beta}(\alpha))$ for all $\alpha \in (0, 1 - \beta]$.
            The map $\alpha \mapsto 1 - \alpha$ is strictly decreasing and by \eqref{def:c-sigma-beta-in-alpha} the map $\alpha \mapsto C_{\sigma, \beta}(\alpha) = \tfrac{(1 - \sigma)\beta}{1 - \alpha}$ is strictly increasing on $(0, 1 - \beta]$.
            Consequently $\max(1 - \alpha, C_{\sigma, \beta}(\alpha))$ is minimized at the point $\alpha^\star$ where the two maps are equal provided $\alpha^\star$ lies in $(0, 1 - \beta]$.
            Solving $1 - \alpha^\star = C_{\sigma, \beta}(\alpha^\star)$ gives
            \begin{align*}
                1 - \alpha^\star = \frac{(1 - \sigma)\beta}{1 - \alpha^\star}
                \iff
                (1 - \alpha^\star)^2 = (1 - \sigma)\beta
                \iff
                \alpha^\star = 1 - \sqrt{(1 - \sigma)\beta}.
            \end{align*}
            \par
            We verify that $\alpha^\star \in (0, 1 - \beta]$.
            Because $(1 - \sigma)\beta \in (0, 1)$, $\alpha^\star \in (0, 1)$.
            Therefore it remains to check $\alpha^\star \le 1 - \beta$:
            \begin{align*}
                \alpha^\star \le 1 - \beta
                &\iff
                \beta \le \sqrt{(1 - \sigma)\beta}
                \underset{(1)}\iff
                \beta^2 \le (1 - \sigma)\beta
                \underset{(2)}\iff
                \beta \le 1 - \sigma.
            \end{align*}
            At (1) squaring both sides is valid because $\beta \ge 0$ and $\sqrt{(1 - \sigma)\beta} \ge 0$.
            At (2) we divided by $\beta > 0$.
            Therefore the product is nonpositive exactly when the second factor is nonpositive.
            \par
            Consider the case $\beta \in [\beta^\star, 1 - \sigma]$.
            Then $\alpha^\star \in (0, 1 - \beta]$.
            Therefore $\alpha^\star$ minimizes $M_{\sigma, \beta}(\alpha)$ over $\alpha \in (0, 1 - \beta]$.
            By \eqref{eqn:mstar-with-alpha-only} the minimum value is
            \begin{align*}
                M_\sigma^*(\beta) = 1 - \alpha^\star = \sqrt{(1 - \sigma)\beta}.
            \end{align*}
            \par
            Consider the case $\beta \in (1 - \sigma, 1)$.
            For all $\alpha \in (0, 1 - \beta]$ it satisfies $1 - \alpha \ge \beta$.
            From \eqref{def:c-sigma-beta-in-alpha} we also observe that $C_{\sigma, \beta}$ is increasing on $(0, 1 - \beta]$ and therefore: 
            \begin{align*}
                C_{\sigma, \beta}(\alpha) \le C_{\sigma, \beta}(1 - \beta) = \tfrac{(1 - \sigma)\beta}{\beta} = 1 - \sigma.
            \end{align*}
            Hence $1 - \alpha \ge \beta > 1 - \sigma \ge C_{\sigma, \beta}(\alpha)$ and $M_{\sigma, \beta}(\alpha) = 1 - \alpha$ on the whole domain $(0, 1 - \beta]$.
            This is minimized by taking $\alpha$ to be as large as possible, namely $\alpha = 1 - \beta$.
            By \eqref{eqn:mstar-with-alpha-only} the minimum value is $M_\sigma^*(\beta) = M_{\sigma, \beta}(1 - \beta) = 1 - (1 - \beta) = \beta$.
        \end{proof}
        The next proposition provides an upper bound of $M_\sigma^*(\beta)$ when $\beta \le \beta^\star$.
        \begin{proposition}[$M_\sigma^*$ bounds when $\beta\le\beta^\star$]
        \label{prop:mstar-beta-le-bstar}
            Fix any $\beta \in (0, 1)$ and $\sigma \in (0, 1)$.
            Let $S:= 1 + \sigma - \beta(1 - \sigma)$ and $\beta^\star = \frac{1 - \sqrt{\sigma}}{1 + \sqrt{\sigma}}$.
            Then $M_\sigma^*$ from \eqref{def:mstar-sigma-beta-bnd} satisfies for all $\beta \le \beta^\star$: 
            \begin{align*}
                M_\sigma^*(\beta) &\le 
                \frac{(1 - \sigma)(1 + \beta) + \sqrt{S^2 - 4\sigma}}{2}. 
            \end{align*}
            Additionally the upper bound is attained by the real number: $\alpha^\star := \frac{S - \sqrt{S^2 - 4\sigma}}{2} \in (\sigma, 1 - \beta)$. 
        \end{proposition}
        \begin{proof}
            \todoinline{
                Proposition \ref{prop:mstar-beta-le-bstar} proof depends on:
                \begin{itemize}[nosep]
                    \item Lemma \ref{lemma:dbeta-equiv}.
                    \item \eqref{eqn:xi-sigma-beta-nonpos-pitem1}.
                    \item Lemma \ref{lemma:xi-sigma-beta-nonpos}.
                    \item \eqref{def:m-sigma-beta-in-alpha}.
                    \item Lemma \ref{lemma:m-in-alpha-beta}.
                \end{itemize}
            }
            We identify an $\alpha^\star$ such that $(\alpha^\star, \sigma) \in \mathcal D_\beta$.
            From the definition of the set it follows that: 
            \begin{align*}
                \sigma = \tfrac{\alpha(1 - \alpha - \beta)}{1 - \alpha(1 + \beta)}
                \iff &
                \sigma(1 - \alpha(1 + \beta)) = \alpha(1 - \alpha - \beta). 
            \end{align*}
            Rearranging the equation gives: 
            \begin{align*}
                0 &= \alpha - \alpha^2 - \alpha \beta - \sigma + \alpha\sigma + \alpha\beta\sigma
                \\
                &= 
                - \alpha^2 + \alpha(1 - \beta + \sigma + \beta\sigma) - \sigma
                \\
                &= 
                - \alpha^2 + \alpha(1 + \sigma - \beta(1 - \sigma)) - \sigma=:p(\alpha).
            \end{align*}
            The equality can be written as $0 = - \alpha^2 + \alpha S - \sigma$ where $S:= 1 + \sigma - \beta(1 - \sigma)$ and observe $S \ge 2 \sigma \ge 0$. 
            The roots are $\alpha^\star = \tfrac{S \pm \sqrt{S^2 - 4\sigma}}{2}$. 
            We claim that the roots are real because $S^2 - 4\sigma \ge 0$ for all $\beta \in (0, \beta^\star]$. 
            To see that, note that $S^2 - 4 \sigma = (S - 2\sqrt{\sigma})(S + 2\sqrt{\sigma})$ so it agrees with the sign of $S - 2\sqrt{\sigma}$.
            Since $\beta \le \beta^\star:= \frac{1 - \sqrt{\sigma}}{1 + \sqrt{\sigma}}$ and $\sigma \in (0,1)$ assumed at the start of the section it follows that $S - 2 \sqrt{\sigma} \ge 0$ because:
            \begin{align*}
                S = 1 + \sigma - \beta(1 - \sigma)
                &\ge 1 + \sigma - \beta^\star(1 - \sigma)
                \\
                &= 1 + \sigma - (1 - \sqrt{\sigma})^2
                \\
                &= 1 + \sigma - (1 + \sigma - 2\sqrt{\sigma}) = -2\sqrt{\sigma}.
            \end{align*}
            Therefore the roots $\alpha^\star$ are always real.
            From Lemma \ref{lemma:dbeta-equiv}, $(\alpha^\star, \sigma) \in \mathcal D_\beta$ implies that $\alpha^\star \in [\sigma, 1 - \beta]$.
            This interval is nonempty because of the assumption $\beta \le \beta^\star$: 
            \begin{align*}
                1 - \beta 
                \ge \tfrac{2\sqrt{\sigma}}{1 + \sqrt{\sigma}} 
                \ge \sqrt{\sigma} > \sigma. 
            \end{align*}
            Recall $p(\alpha):= - \alpha^2 + \alpha S - \sigma$. 
            We show below that $p(\sigma) < 0$ and $p(1 - \beta) < 0$.
            Also $p(\sqrt\sigma) = \sqrt\sigma(S - 2\sqrt\sigma) \ge 0$ because $S - 2\sqrt\sigma \ge 0$ was shown above.
            Moreover $\sigma < \sqrt\sigma \le 1 - \beta$ because $\sigma \in (0, 1)$ and $\beta \le \beta^\star \le 1 - \sqrt\sigma$.
            A downward quadratic function is negative outside its roots and nonnegative between them.
            Hence $\sigma$ and $1 - \beta$ lie outside the roots while $\sqrt\sigma$ lies between them.
            Therefore both roots $\alpha^\star$ of $p(\alpha) = 0$ lie in the interval $(\sigma, 1 - \beta)$.
            $p$ is negative at $\alpha = \sigma, \alpha = 1 - \beta$ because: 
            \begin{align*}
                p(\sigma) &= - \sigma^2 + \sigma(1 + \sigma - \beta(1 - \sigma)) - \sigma
                \\
                &= - \sigma^2 + \sigma + \sigma^2 - \beta\sigma + \beta\sigma^2 - \sigma
                \\
                &= \beta\sigma(\sigma - 1) < 0, 
                \\
                p(1 - \beta) &=
                - (1 - \beta)^2 + (1 - \beta)(1 + \sigma - \beta(1 - \sigma)) - \sigma
                \\
                &= (1 - \beta)(
                    - 1 + \beta + 1 + \sigma - \beta + \beta\sigma
                ) - \sigma
                \\
                &= (1 - \beta)(\sigma + \beta\sigma) - \sigma
                \\
                &= \sigma((1 - \beta)(1 + \beta) - 1) 
                \\
                &= \sigma(1 - \beta^2 - 1)
                \\
                &= - \sigma\beta^2 < 0. 
            \end{align*}
            Therefore, $\left(\tfrac{S \pm \sqrt{S^2 - 4\sigma}}{2}, \sigma\right) \in \mathcal D_\beta$. 
            We will use this element to upper bound $M_{\sigma}^*(\beta)$.
            From \eqref{eqn:xi-sigma-beta-nonpos-pitem1} in the proof of Lemma \ref{lemma:xi-sigma-beta-nonpos} we make the observation that: 
            \begin{align*}
                \tfrac{\sigma}{\alpha}(1 - \alpha)\xi_{\sigma, \beta}(\alpha) 
                = p(\alpha)
                \implies 
                \xi_{\sigma, \beta}(\alpha^\star) = 0.
            \end{align*}
            Hence $M_{\sigma, \beta}(\alpha^\star) = \max(1 - \alpha^\star, C_{\sigma, \beta}(\alpha^\star))$. 
            Observe that when $q = \sigma$ combining \eqref{def:m-sigma-beta-in-alpha} and Lemma \ref{lemma:m-in-alpha-beta} yields:
            \begin{align*}
                C_{\sigma, \beta}(\alpha) = C_\sigma(\alpha, \sigma)
                = \tfrac{(1 - \sigma)(\alpha - \sigma)}{\alpha(1 - \sigma)} = 1 - \tfrac{\sigma}{\alpha}. 
            \end{align*}
            Denote by $\alpha_+^\star$ the larger root and by $\alpha_-^\star$ the smaller one.
            Then it follows that: 
            \begin{align*}
                \left(\alpha_-^\star\right)^2 
                \le \alpha_-^\star\alpha_+^\star = \sigma 
                \implies 
                \alpha_-^\star\le \tfrac{\sigma}{\alpha_-^\star}
                \iff 
                1 - \alpha_-^\star
                \ge 1 - \tfrac{\sigma}{\alpha_-^\star}. 
            \end{align*}
            Therefore we have:
            \begin{align*}
                M_{\sigma, \beta}(\alpha_-^\star) 
                &= \max\left(
                    1 - \alpha_-^\star, 1 - \tfrac{\sigma}{\alpha_-^\star}
                \right) 
                \\
                &= 
                1 - \frac{S - \sqrt{S^2 - 4\sigma}}{2}
                \\
                &= 
                \frac{1 - \sigma + \beta(1 - \sigma) + \sqrt{S^2 - 4\sigma}}{2}
                \\
                &= 
                \frac{(1 - \sigma)(1 + \beta) + \sqrt{S^2 - 4\sigma}}{2} \ge M_{\sigma}^*(\beta). 
            \end{align*}
        \end{proof}
        By Propositions \ref{prop:mstar-beta-ge-betastar} and \ref{prop:mstar-beta-le-bstar} we identify an upper bound of $M_\sigma^*(\beta)$:
        \begin{align}
            \overline M_\sigma(\beta) &:= 
            \begin{cases}
                \frac{(1 - \sigma)(1 + \beta) + \sqrt{S^2 - 4\sigma}}{2},
                & \beta \le \beta^\star, 
                \\
                \sqrt{(1 - \sigma)\beta}, & \beta \in (\beta^\star, 1 - \sigma),
                \\
                \beta, & \beta \in [1 - \sigma, 1),
            \end{cases}
            \label{def:overline-m}
            \\
            \beta^\star &:= \frac{1 - \sqrt{\sigma}}{1 + \sqrt{\sigma}}, 
            \label{def:beta-star}
            \\
            S &:= 1 + \sigma - \beta(1 - \sigma),
            \label{def:big-s}
        \end{align}
        In relation to \eqref{def:mstar-sigma-beta-bnd} it satisfies for all $\beta \in (0, 1)$ the inequality $M_\sigma^*(\beta) \le \overline M_\sigma(\beta)$ with equality when $\beta \in [\beta^\star, 1)$.
        Lemma \ref{lemma:dbeta-equiv} also implies that there exists $(\overline\alpha_\sigma(\beta), \overline q_\sigma(\beta)) \in \mathcal D_\beta$ such that $\overline M_\sigma(\beta)$ is attained for all $\beta \in (0, 1)$ and they are given by:
        \begin{align}
            \overline\alpha_\sigma(\beta) &:= 
            \begin{cases}
                \tfrac{S - \sqrt{S^2 - 4\sigma}}{2}, & \beta\in (0, \beta^\star), 
                \\
                1 - \sqrt{(1 - \sigma)\beta}, & \beta \in [\beta^\star, 1 - \sigma],
                \\
                1 - \beta, & \beta \in (1 - \sigma, 1),
            \end{cases}
            \label{def:overline-alpha}\\
            \overline q_\sigma(\beta) &:= 
            \frac{
                \overline\alpha_\sigma(\beta)(1 - \overline\alpha_\sigma(\beta) - \beta)
            }{
                1 - \overline\alpha_\sigma(\beta)(1 + \beta)
            }.
            \label{def:overline-beta}
        \end{align}
        Therefore for all $\beta\in(0, 1)$:
        \begin{align}\label{eqn:overline-m-realized}
            \overline M_\sigma(\beta) &= M_\sigma(\overline\alpha_\sigma(\beta), \overline q_\sigma(\beta)).
        \end{align}
        Our second main result of the paper is the theorem below.
        It presents the convergence rate of the constant momentum algorithm for all $\beta \in (0, 1)$.
        \begin{theorem}[Linear convergence for all $\beta \in (0, 1)$]\;\label{thm:lin-cnvg-constn-mmntm}\\
            Fix any $\beta \in (0, 1)$ and $\hat x \in \RR^n$ and let $\sigma$ be given by \eqref{def:sigma} such that $\sigma \in (0, 1)$.
            Suppose that:
            \begin{enumerate}[nosep]
                \item the sequences $(x_k)_{k \in \N}, (x_k^+)_{k \in \N}$ satisfy Definition \ref{def:const-momentum-alg} with initial condition $x_{-1}^+ = x_0^+ = T_{L}^{(f, g)}(\hat x)$.
                \item $(F, f, g, L, \mu)$ satisfy Assumption \ref{ass:strong-convexity}.
            \end{enumerate}
            Let $\bar x$ be the unique minimizer of $F = f + g$.
            Then the iterates $x_k^+$ satisfy for all $k\in \N$:
            \begin{align*}
                F(x_k^+) - F(\bar x) &\le \overline M_\sigma(\beta)^k\cdot \frac{L}{2}\Vert \hat x - \bar x\Vert^2,
            \end{align*}
            where $\overline M_\sigma(\beta)$ is given by \eqref{def:overline-m} and lies in the interval $(0, 1)$.
        \end{theorem}
        \begin{proof}
            \todoinline{
                Theorem \ref{thm:lin-cnvg-constn-mmntm} proof depends on:
                \begin{itemize}[nosep]
                    \item \eqref{eqn:overline-m-realized}.
                    \item Definition \ref{def:apg-algorithm}.
                    \item Theorem \ref{thm:apg-momentum-equiv}.
                    \item Assumption \ref{ass:const-param}.
                    \item Assumption \ref{ass:strong-convexity}.
                    
                    \item Lemma \ref{lemma:apg-lya-ineq-mu}.
                    \item \eqref{def:m-sigma-alpha-q}.
                    \item \eqref{def:upsilon}.
                    \item \eqref{lemma:apg-lya-const-param-eqn2}.
                    \item \eqref{def:overline-m}.
                    \item \eqref{def:overline-alpha}.
                    \item \eqref{def:overline-beta}.
                    \item \eqref{def:big-s}.
                    \item Proposition \ref{prop:mstar-beta-ge-betastar}.
                    \item Proposition \ref{prop:mstar-beta-le-bstar}.
                \end{itemize}
                The proof implicitly depends on:
                \begin{itemize}[nosep]
                    \item Lemma \ref{lemma:dbeta-equiv}.
                    \item \eqref{def:delta}.
                    \item \eqref{def:dbeta-set}.
                    \item \eqref{def:sigma}.
                \end{itemize}
            }
            Write for short $(\overline \alpha, \overline q) := (\overline\alpha_\sigma(\beta), \overline q_\sigma(\beta))$ given by \eqref{def:overline-alpha} and \eqref{def:overline-beta}.
            By \eqref{eqn:overline-m-realized}, $M_\sigma(\overline \alpha, \overline q) = \overline M_\sigma(\beta)$ where $\overline M_\sigma$ is given by \eqref{def:overline-m}.
            Propositions \ref{prop:mstar-beta-ge-betastar} and \ref{prop:mstar-beta-le-bstar} together prove that $(\overline \alpha, \overline q) \in \mathcal D_\beta$.
            \par
            Let $(z_k)_{k \in \Z_+}, (z_k^+)_{k \in \Z_+}, (z_k^\dagger)_{k \in \Z_+}$ satisfy Definition \ref{def:apg-algorithm} with initial conditions $z_{-1}^+ = z_{-1}^\dagger := \hat x$ and let $(\alpha_k)_{k \in \Z_+}, (q_k)_{k \in \Z_+}, (\gamma_k)_{k \in \Z_+}$ satisfy Assumption \ref{ass:const-param} with constant pair $(\overline \alpha, \overline q)$.
            The following holds for all $k\in \N$:
            \begin{align*}
                F(x_k^+) - F(\bar x)
                \underset{(1)}{=} F(z_k^+) - F(\bar x) = \Delta_k(\bar x)
                \underset{(2)}\le \Upsilon_k(\bar x)
                \underset{(3)}\le M_\sigma(\overline\alpha, \overline q)^{k}\, \Upsilon_0(\bar x)
                \underset{(4)}\le \tfrac{L\overline M_\sigma(\beta)^k}{2}\Vert \bar x - \hat x\Vert^2.
            \end{align*}
            \par
            At (1) Theorem \ref{thm:apg-momentum-equiv} applies to $(x_k)_{k \in \N}, (x_k^+)_{k \in \N}$ and $(z_k)_{k \in \Z_+}, (z_k^+)_{k \in \Z_+}, (z_k^\dagger)_{k \in \Z_+}$ with parameter $\beta$ because $(\overline \alpha, \overline q) \in \mathcal D_\beta$.
            It gives $x_k^+ = z_k^+$ for every $k \in \N$.
            \par
            At (2), \eqref{def:upsilon} gives $\Upsilon_k(\bar x) = \Delta_k(\bar x) + \tfrac{L \overline\alpha^2}{2}\Vert z_k^\dagger - \bar x\Vert^2 \ge \Delta_k(\bar x)$.
            \par
            At (3) Lemma \ref{lemma:apg-lya-ineq-mu} applies to $(z_j)_{j \in \Z_+}, (z_j^+)_{j \in \Z_+}, (z_j^\dagger)_{j \in \Z_+}$ because $(\alpha_j, q_j, \gamma_j)$ satisfy Assumption \ref{ass:const-param} with constant pair $(\overline\alpha, \overline q)$ and $(F, f, g, L, \mu)$ satisfy Assumption \ref{ass:strong-convexity}.
            It gives $\Upsilon_j(\bar x) \le M_\sigma(\overline\alpha, \overline q) \Upsilon_{j - 1}(\bar x)$ for every $j \ge 1$ where $M_\sigma$ is given by \eqref{def:m-sigma-alpha-q}.
            Since $M_\sigma(\overline\alpha, \overline q) \ge 1 - \overline\alpha \ge 0$ by \eqref{def:m-sigma-alpha-q} and $\overline\alpha \in (0, 1]$ multiplying an inequality by it preserves the inequality.
            Applying this for $j = k, k - 1, \ldots, 1$ by induction on $k$ gives $\Upsilon_k(\bar x) \le M_\sigma(\overline\alpha, \overline q)^k \Upsilon_0(\bar x)$.
            \par
            At (4) we use $M_\sigma(\overline\alpha, \overline q) = \overline M_\sigma(\beta)$.
            Since $z_{-1}^+ = z_{-1}^\dagger = \hat x$ the inequality \eqref{lemma:apg-lya-const-param-eqn2} gives $z_0 = \hat x$ and $\Upsilon_0(\bar x) \le \tfrac{L}{2}\Vert \bar x - z_0\Vert^2 = \tfrac{L}{2}\Vert \bar x - \hat x\Vert^2$.
            \par
            Finally we show that $\overline M_\sigma(\beta)$ from \eqref{def:overline-m} is in the interval $(0, 1)$ for all $\beta \in (0, 1)$.
            For all $\beta \in (\beta^\star, 1)$ the value $\overline M_\sigma(\beta)$ is either $\sqrt{(1 - \sigma)\beta}$ or $\beta$ by \eqref{def:overline-m}.
            Both cases lie in $(0, 1)$ because $\sigma, \beta \in (0, 1)$.
            For all $\beta \in (0, \beta^\star]$ observe that \eqref{def:overline-alpha}, \eqref{def:big-s} together yield the equalities:
            \begin{align*}
                1 - \overline \alpha_\sigma(\beta) &= 1 - \tfrac{S - \sqrt{S^2 - 4\sigma}}{2}
                = \tfrac{(1 - \sigma)(1 + \beta) + \sqrt{S^2 - 4\sigma}}{2}
                = \overline M_\sigma(\beta). 
            \end{align*}
            From Proposition \ref{prop:mstar-beta-le-bstar} it follows that $\overline \alpha_\sigma(\beta) \in (0, 1 - \beta)$ and hence for all $\beta \in (0, \beta^\star]$ the above says $\overline M_\sigma(\beta) = 1 - \overline \alpha_\sigma(\beta) \in (\beta, 1) \subset (0, 1)$.
        \end{proof}

    \subsection{Numerical illustrations of upper bounds}\label{ssec:numerical-illustrations}
        In this section we present plots to show that $\overline M_\sigma(\beta)$ is loose when $\beta \in (0, \beta^\star)$ by making the comparison between $\overline M_\sigma$ and $M^*_\sigma$. 
        For all values of $\beta \in (0, 1)$ we conduct a grid search to find the minimum and minimizers of $M_\sigma(\alpha, q)$ over $(\alpha, q) \in \mathcal D_\beta$.
        We then compare them to those attaining the upper bound $\overline M_\sigma(\beta)$.
        We conclude that whenever $\overline M_\sigma(\beta) > M_\sigma^*(\beta)$ the numerically obtained minimizer makes $\xi_{\sigma, \beta}(\alpha)$ as defined in \eqref{def:xi-sigma-beta-in-alpha} greater than zero.
        We also visualize the numerically obtained minimizer $(\alpha, q)\in \mathcal D_\beta$ of \eqref{eqn:mstar-with-alpha-only} against those from \eqref{def:overline-alpha}, \eqref{def:overline-beta} that attain the upper bound $\overline M_\sigma$ for $\sigma = 0.1$ and $\beta\in(0, 1)$.
        \par
        For discussion we fix $\beta \in (0, 1)$ and let $\sigma \in (0, 1)$ be given by \eqref{def:sigma}, and write the numerically obtained minimizer of a given $\beta$ satisfying $M_{\sigma}(\alpha^\star, q^\star) = M_\sigma^*(\beta)$ as $(\alpha^\star, q^\star)$.
        We write $(\overline \alpha, \overline q)$ for the pair satisfying $M_\sigma(\overline \alpha, \overline q) = \overline M_\sigma(\beta)$ as defined in \eqref{def:overline-alpha} and \eqref{def:overline-beta}.
        Propositions \ref{prop:mstar-beta-ge-betastar} and \ref{prop:mstar-beta-le-bstar} prove that this pair exists.
        \par
\begin{figure}[H]
            \centering
            \subfloat[$M_\sigma^*(\beta)$]{
                \label{fig:mstar-beta:a}\includegraphics[width=0.47\textwidth]{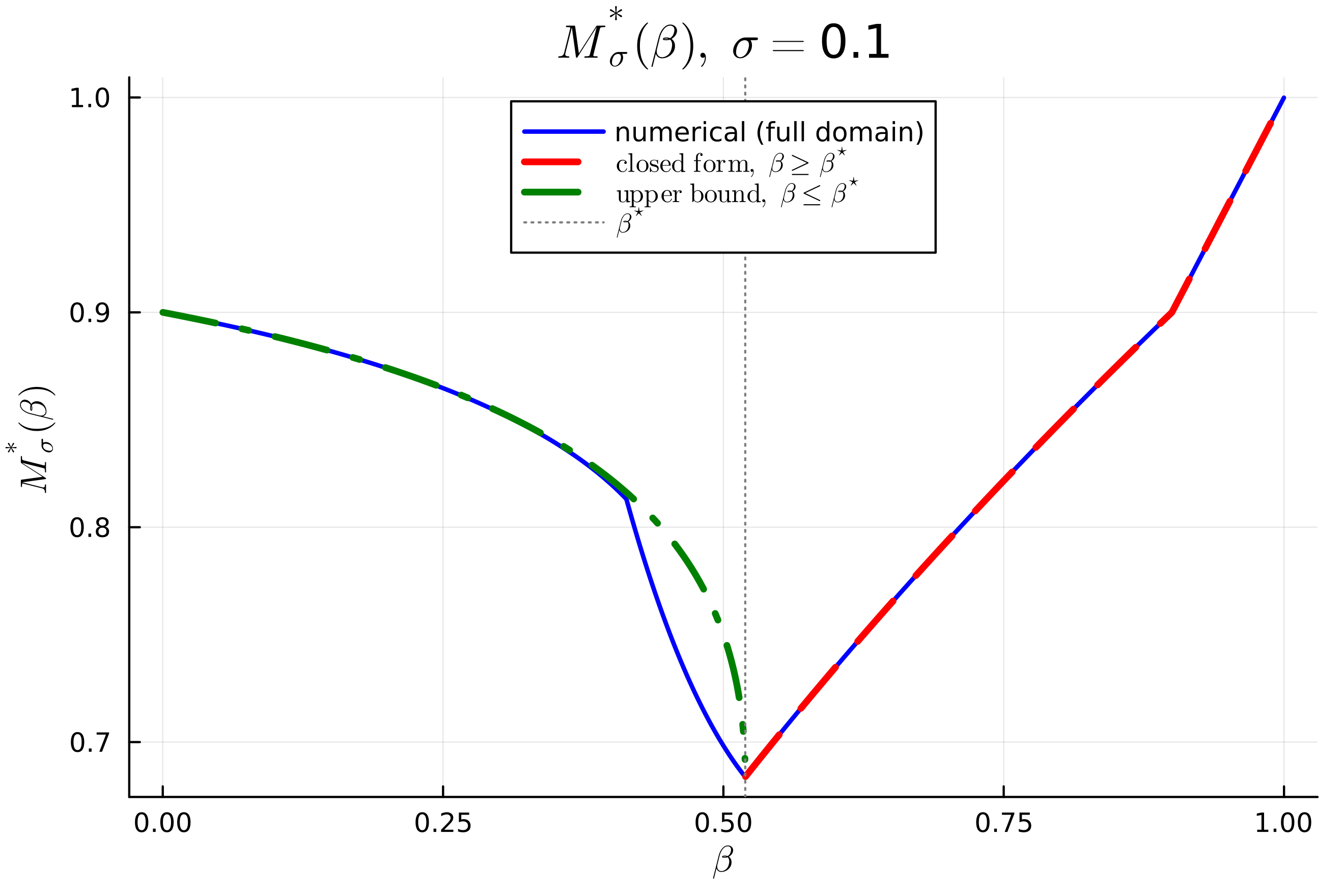}
            }
            \hfill
            \subfloat[relative optimality vs $\rho$]{
                \label{fig:mstar-beta:b}\includegraphics[width=0.47\textwidth]{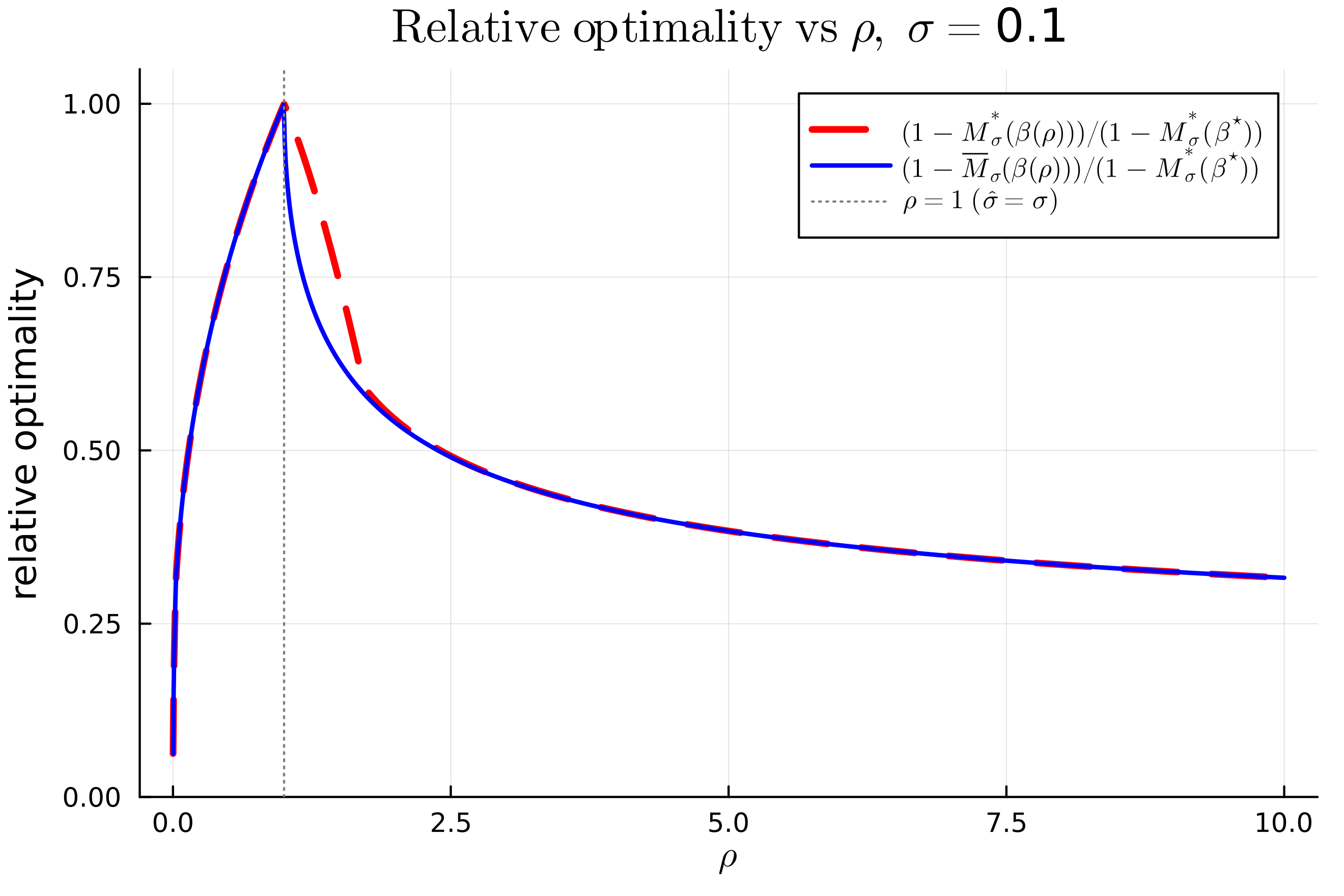}
            }
            \caption{
                (a) illustrates $M_\sigma^*(\beta)$ at $\sigma = 0.1$.
                It plots $M_{\sigma, \beta}(\beta), \overline M_{\sigma, \beta}(\beta)$ over $\beta\in(0, 1)$ with $[1 - \sqrt{\sigma}, 1]$ as the range of the vertical axis.
                (b) plots relative optimality for all $\rho \in (0, \sigma^{-1}]$.
                The dashed curve is the exact value using $M_\sigma^*$.
                The solid curve is the achievable value using $\overline M_\sigma$.
            }
            \label{fig:mstar-beta}
        \end{figure}

        Figure \ref{fig:mstar-beta} illustrates Propositions \ref{prop:mstar-beta-ge-betastar} and \ref{prop:mstar-beta-le-bstar} at $\sigma = 0.1$.
        Figure \ref{fig:mstar-beta:a} illustrates $M_\sigma^*(\beta)$.
        It compares $M_\sigma^*(\beta)$ against the upper bound from Proposition \ref{prop:mstar-beta-ge-betastar} and Proposition \ref{prop:mstar-beta-le-bstar}.
        \par
        Figure \ref{fig:mstar-beta:b} illustrates the effect of estimating $\sigma$ using the scaling factor $\rho \in (0, \sigma^{-1}]$.
        Suppose that under some practical settings we are given an estimate of $\sigma$ by $\hat\sigma := \rho\sigma$ for a $\rho \in (0, \sigma^{-1}]$.
        Suppose that we also use the formula $\frac{1 - \sqrt{\sigma}}{1 + \sqrt{\sigma}}$ for $\beta^\star$ with $\hat\sigma$ in place of $\sigma$ to produce the momentum constant: $\beta(\rho) := \tfrac{1 - \sqrt{\rho\sigma}}{1 + \sqrt{\rho\sigma}} \in [0, 1)$.
        This simulates a common situation in practice where $\sigma$ is imprecise but we still use the formula for the momentum parameter $\beta^\star$ given by V-FISTA.
        Observe that when $\rho = 1$ it follows that $\beta(1) = \beta^\star$.
        Additionally in the figure we visualize the relative optimality $\tfrac{1 - M_\sigma^*(\beta(\rho))}{1 - M_\sigma^*(\beta^\star)}$. 
        This is the ratio of the contraction gap $1 - M_\sigma^*(\beta(\rho))$ caused by the mistuned estimate $\rho\sigma$ of $\sigma$ to the best possible contraction gap $1 - M_\sigma^*(\beta^\star) = \sqrt\sigma$ produced when $\rho = 1$.
        Figure \ref{fig:mstar-beta:b} shows this quantity bounded below by $\sqrt\sigma$ for every $\rho \in [1, \sigma^{-1}]$.
        This happens when $\rho\sigma$ is an overestimate of $\sigma$ and makes $\beta(\rho) \le \beta^\star$.
        Otherwise, in the interval $\rho \in (0, 1]$ the ratio decreases to zero.
        This happens when $\rho\sigma$ is an underestimate of $\sigma$ and makes $\beta(\rho) > \beta^\star$.
        \par
\begin{figure}[H]
            \centering
            \subfloat[looseness]{
                \label{fig:mstar-argmin:a}\includegraphics[width=0.47\textwidth]{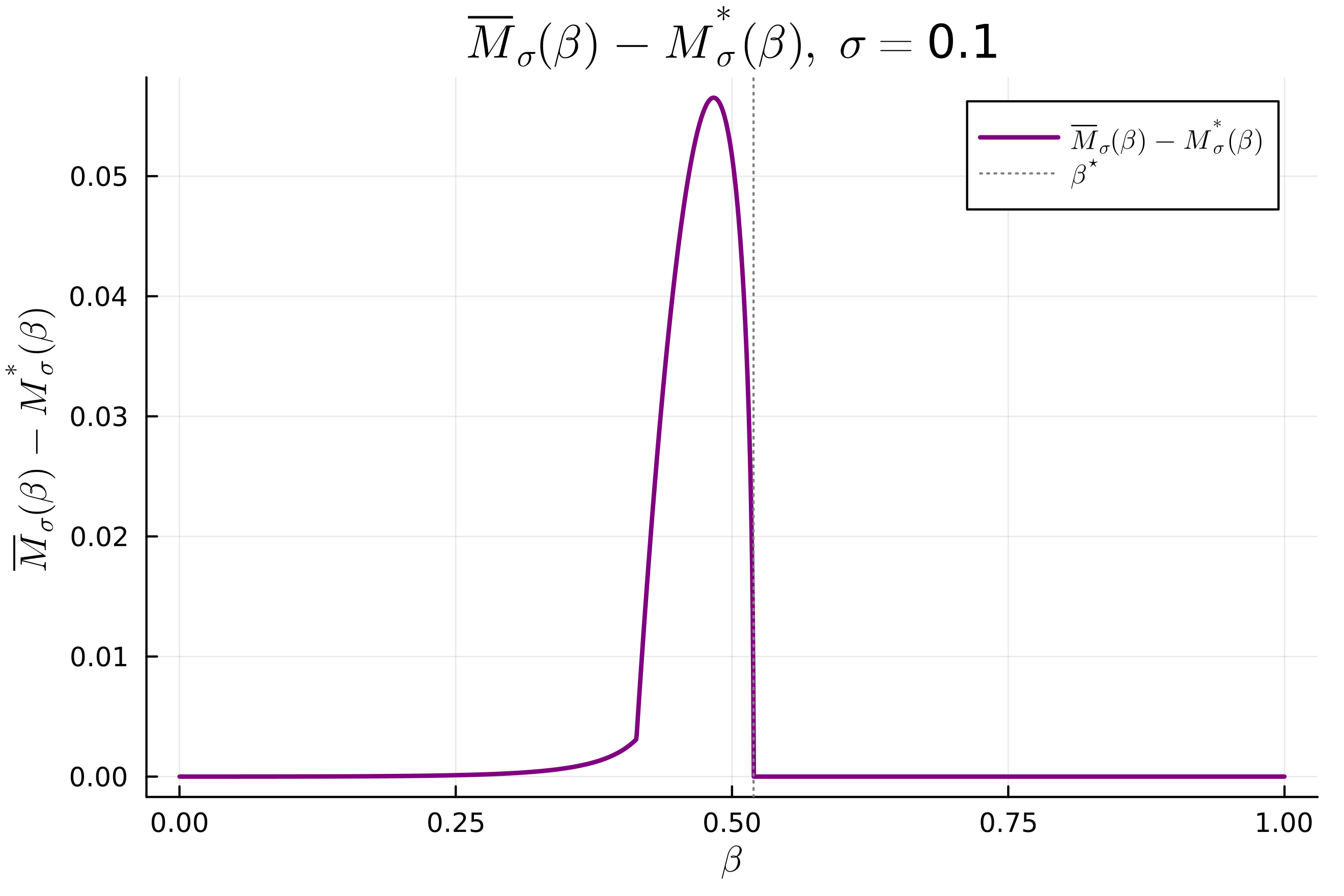}
            }
            \hfill
            \subfloat[$\xi$ at the minimizer]{
                \label{fig:mstar-argmin:b}\includegraphics[width=0.47\textwidth]{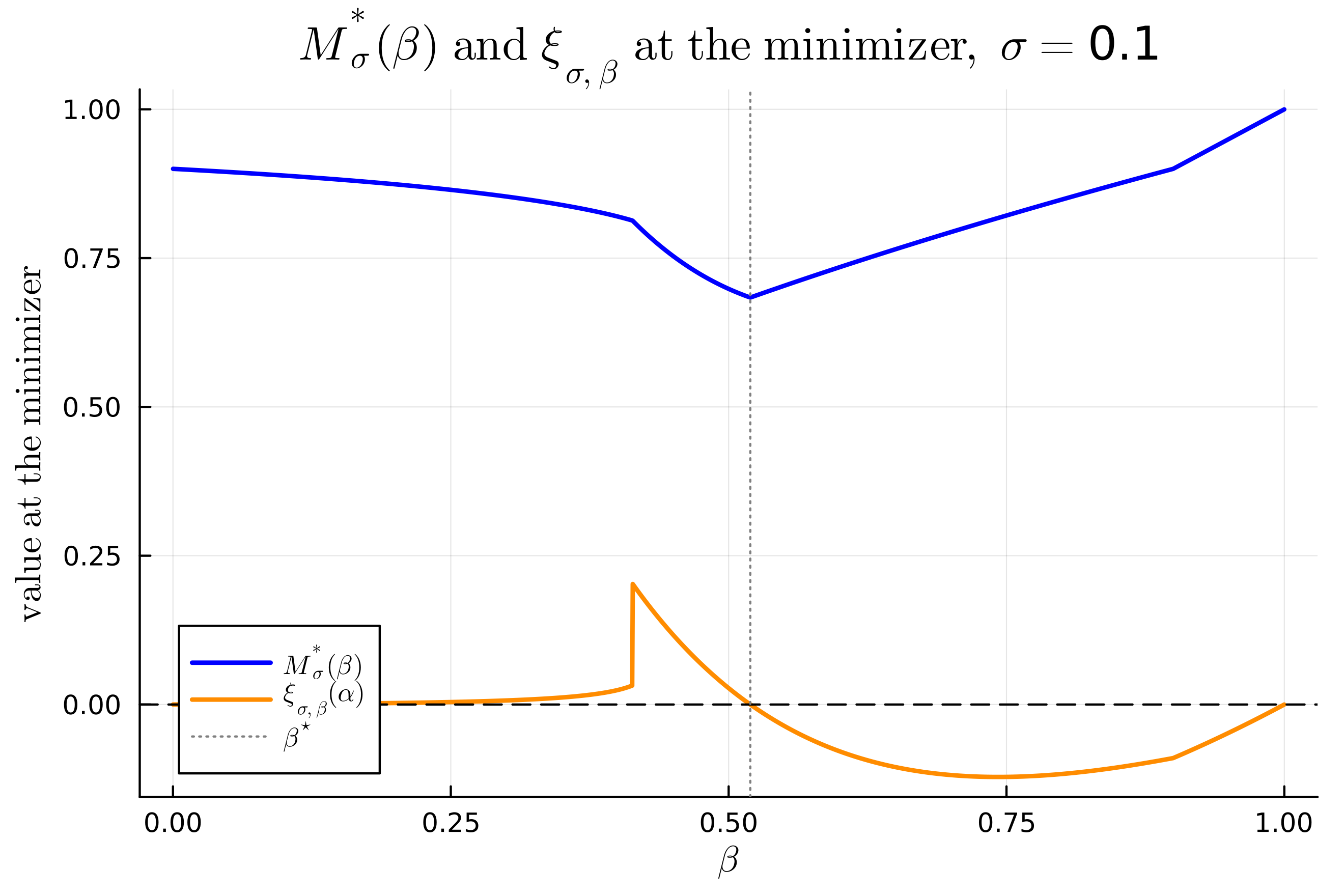}
            }
            \\
            \subfloat[minimizer $\alpha$]{
                \label{fig:mstar-argmin:c}\includegraphics[width=0.47\textwidth]{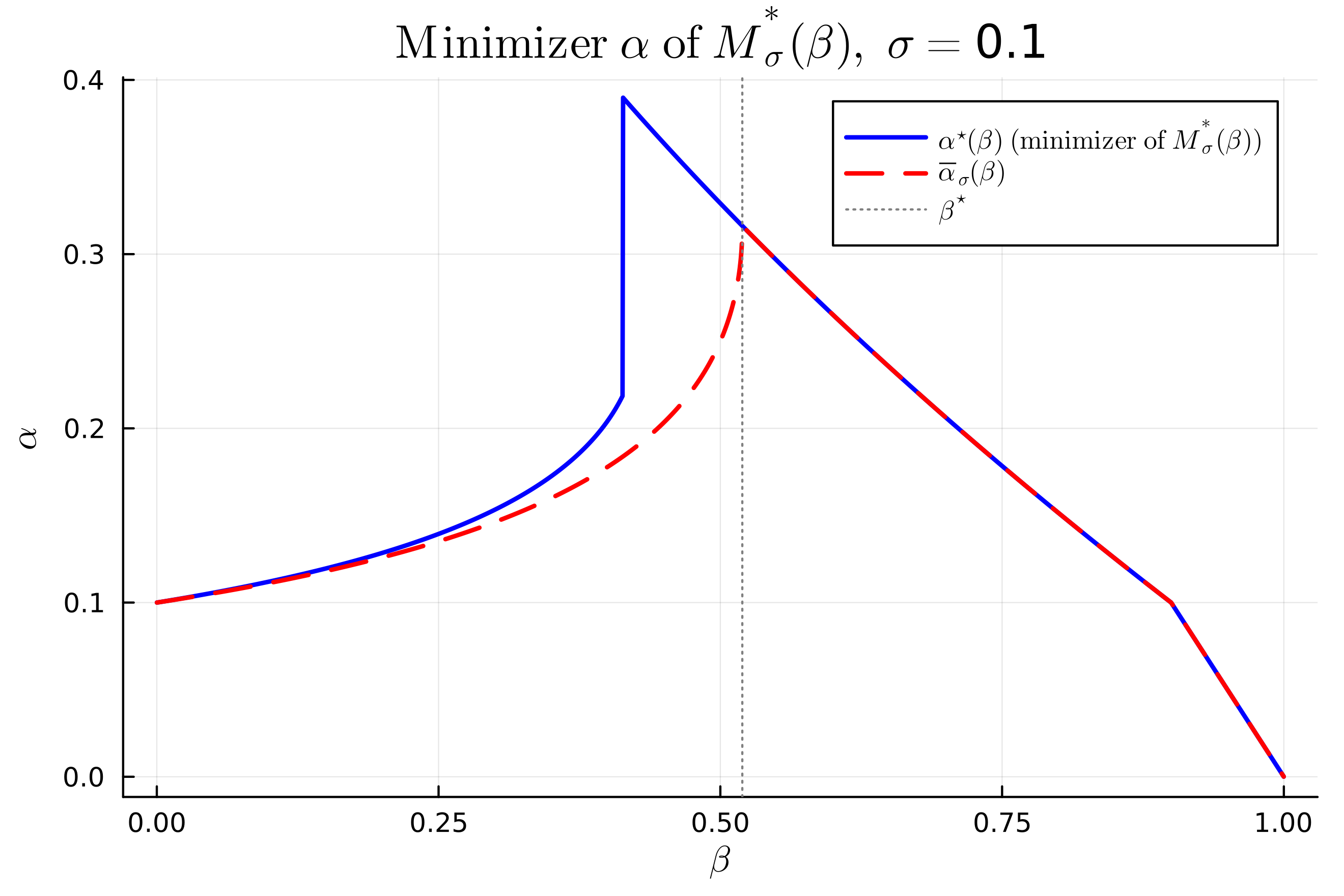}
            }
            \hfill
            \subfloat[minimizer $q$]{
                \label{fig:mstar-argmin:d}\includegraphics[width=0.47\textwidth]{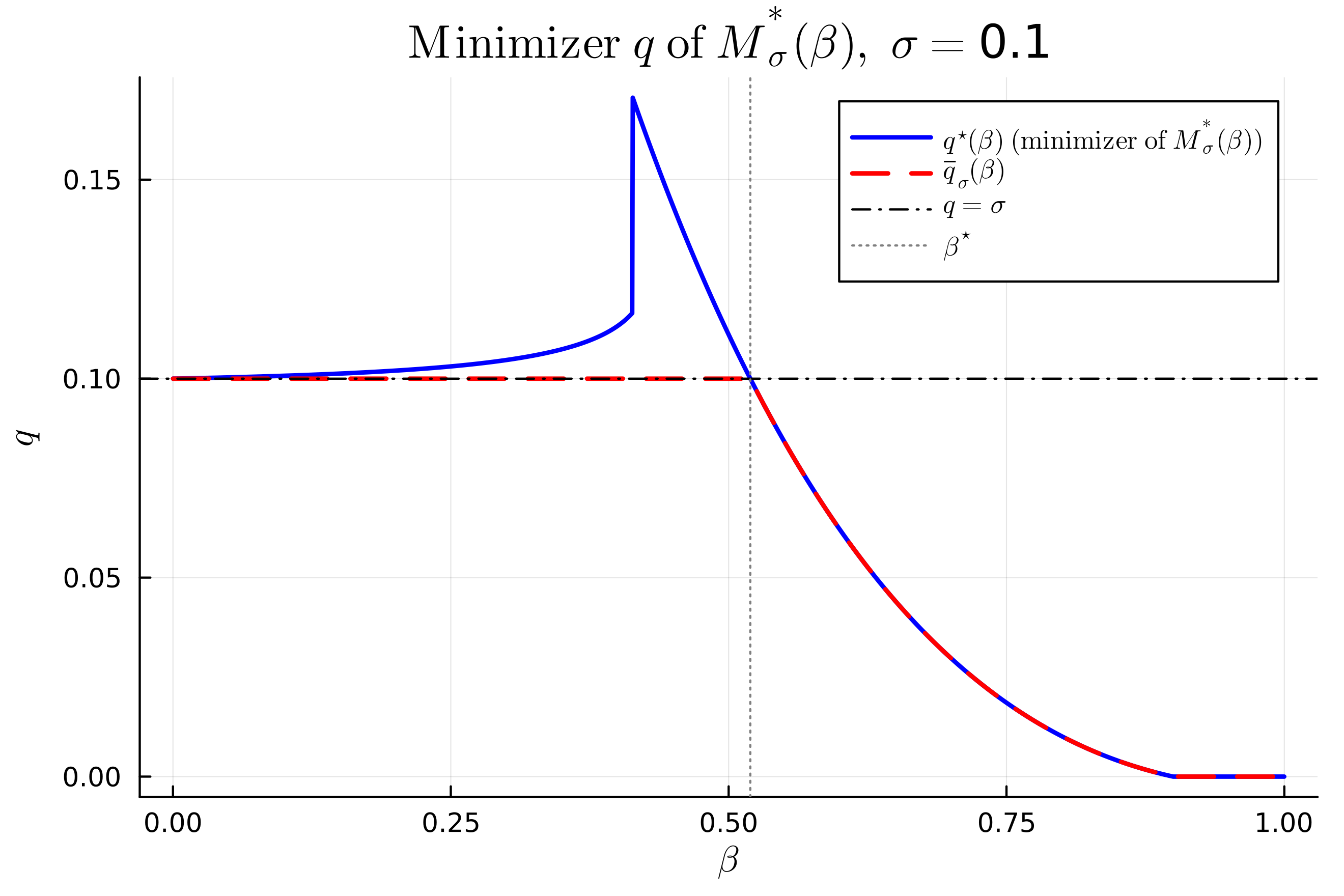}
            }
            \caption{
                (a) $\overline M_\sigma(\beta) - M_\sigma^*(\beta)$.
                (b) $\xi_{\sigma,\beta}$ and the $\max(1-\alpha,C)$ gap, both at the minimizer.
                (c) minimizer $\alpha^\star$ and $\overline \alpha$.
                (d) minimizer $q^\star$ and $\overline q$.
                All at $\sigma = 0.1$.
            }
            \label{fig:mstar-argmin}
        \end{figure}

        Figure \ref{fig:mstar-argmin} illustrates the gap between $\overline M_\sigma$ and $M_\sigma^*$ for values of $\beta < \beta^\star$. 
        It aims to convey the fact that the true minimizer of $M_\sigma^*$ makes $\xi_{\sigma, \beta}$ larger than zero. 
        \par
        Figure \ref{fig:mstar-argmin:a} plots the gap $\overline M_\sigma(\beta) - M_\sigma^*(\beta)$.
        It is zero for all $\beta \ge \beta^\star$.
        Hence $\overline M_\sigma(\beta) = M_\sigma^*(\beta)$.
        However, it is positive when $\beta < \beta^\star$.
        \par
        Figure \ref{fig:mstar-argmin:b} examines the gap further by plotting $\xi_{\sigma,\beta}$ together with $M_\sigma^*(\beta)$.
        It shows that the gap in Figure \ref{fig:mstar-argmin:a} is nonzero when $\xi_{\sigma,\beta}$ is positive.
        This is not assumed in Proposition \ref{prop:mstar-beta-le-bstar} when deriving $\overline M_\sigma$ for the sake of simpler algebra\footnote{
            The minimizer $\alpha^\star$ of $M_\sigma^*$ can be solved analytically and it is the root of a cubic equation because $\xi_{\sigma, \beta}$ is involved.
            We choose not to present it and settle for the upper bound and these figures instead.
        }.
        Therefore $\xi_{\sigma, \beta}$ plays a role in minimizing $M_{\sigma,\beta}$ when $\beta < \beta^\star$.
        \par
        Figures \ref{fig:mstar-argmin:c} and \ref{fig:mstar-argmin:d} compare $(\alpha^\star, q^\star)$ of \eqref{eqn:mstar-with-alpha-only} with $(\overline\alpha, \overline q)$ given by Propositions \ref{prop:mstar-beta-le-bstar} and \ref{prop:mstar-beta-ge-betastar}.
        The curves illustrate that $\alpha^\star \ge \overline \alpha$ and $q^\star \ge \overline q$ when $\beta \in (0, \beta^\star)$.
        These two figures raise the suspicion that the trajectory of $\alpha^\star, q^\star$  is not a continuous function with respect to $\beta$.

\section{PEP Verification of Our Convergence Results}\label{sec:pep-verification}
    Our main result Theorem \ref{thm:lin-cnvg-constn-mmntm} is numerically confirmed via the performance estimation problem (PEP) of Drori and Teboulle \cite{drori_performance_2014} and Taylor et al. \cite{taylor_smooth_2016}.
    The method of PEP computes the exact worst case value of a performance metric over $(F, f, g, L, \mu)$ satisfying Assumption \ref{ass:global-setup} as a semidefinite program.
    For $L = 1$ the choice $\mu = \sigma$ realizes \eqref{def:sigma}, and for each $(\sigma, \beta, k)$ we solve for
    \begin{align*}
        \tau(\sigma, \beta, k) := \max\left\{ F(x_k^+) - F(\bar x) \,:\, \Vert \bar x - \hat x \Vert^2 \le 1 \right\}.
    \end{align*}
    Here $(x_k)_{k \in \N}, (x_k^+)_{k \in \N}$ satisfy Definition \ref{def:const-momentum-alg} with parameter $\beta$ and initial conditions $x_{-1}^+ = x_0^+ := T_L^{(f,g)}(\hat x)$ which is the hypothesis of Theorem \ref{thm:lin-cnvg-constn-mmntm}.
    \par
    Tables \ref{table:pep-verification-b1}, \ref{table:pep-verification-b2}, and \ref{table:pep-verification-b3} report $\tau(\sigma, \beta, k)$ against the theorem's bound $\overline M_\sigma(\beta)^{k} \cdot \tfrac{L}{2}$, one table per case of $\overline M_\sigma$ in \eqref{def:overline-m}.
    Each table lists the theorem's bound $\overline M_\sigma(\beta)^{k} \cdot \tfrac{L}{2}$ in the fourth column and the ratio $\tfrac{\tau(\sigma, \beta, k)}{\overline M_\sigma(\beta)^{k} \cdot L/2}$ in the sixth column.
    Each table also reports a numerical estimate of $M_\sigma^*(\beta)$ in the seventh column and the ratio $\tfrac{\tau(\sigma, \beta, k)}{M_\sigma^*(\beta)^{k} \cdot L/2}$ in the eighth column.
    We estimate $M_\sigma^*(\beta)$ in \eqref{eqn:mstar-with-alpha-only} using the method of grid search.
    Every table showcases $\sigma \in \{0.1, 0.05, 0.025\}$ and $k \in \{8, 16, 32\}$.
    Each case is sampled at two values of $\beta$ evenly spaced across the domain.
    Table \ref{table:pep-verification-b2} additionally includes the case boundary $\beta = 1 - \sigma$.
    \par
    There are two separate numerical issues affecting these tables.
    First, our numerical experiments exclude the exact value $\beta = \beta^\star$ from every table because the semidefinite program always fails to converge within a five minute solver budget using both CLARABEL\footnote{\href{https://clarabel.org/stable/}{https://clarabel.org/stable/}} and SCS\footnote{\href{https://www.cvxgrp.org/scs/}{https://www.cvxgrp.org/scs/}}.
    Second, within the sampled grid three instances in Table \ref{table:pep-verification-b2} at $\sigma = 0.1$ and $k = 32$ return no value.
    CLARABEL reports that the solve failed and the fallback SCS exceeds the five minute time limit.
    These three entries appear as NaN in the table.
    \par
    Every instance of Tables \ref{table:pep-verification-b1}, \ref{table:pep-verification-b2}, and \ref{table:pep-verification-b3} with a value of $\tau(\sigma, \beta, k)$ satisfies $\tau(\sigma, \beta, k) \le \overline M_\sigma(\beta)^{k} \cdot \tfrac{L}{2}$.
    These results are consistent with Theorem \ref{thm:lin-cnvg-constn-mmntm}.
    Each $\tau(\sigma, \beta, k)$ is the optimum solved by CLARABEL implemented using PEPit \cite{pepit2024}.
    CLARABEL falls back to SCS if a solver failure occurs.
    {
        \footnotesize
        \captionsetup{font=normalsize}
        \begin{longtable}{|c|c|c|c|c|c|c|c|}
            \caption{PEP verification for $\beta \le \beta^\star$}\label{table:pep-verification-b1}\\
            \hline
            \rule{0pt}{2.8ex}$\sigma$ & $\beta$ & $k$ & $\overline M_\sigma(\beta)^{k} \cdot \tfrac{L}{2}$ & $\tau(\sigma,\beta,k)$ & $\tfrac{\tau(\sigma, \beta, k)}{\overline M_\sigma(\beta)^{k} \cdot L/2}$ & $M_\sigma^*(\beta)$ & $\tfrac{\tau(\sigma, \beta, k)}{M_\sigma^*(\beta)^{k} \cdot L/2}$ \\
            \hline\hline
            \endfirsthead
            \hline
            \rule{0pt}{2.8ex}$\sigma$ & $\beta$ & $k$ & $\overline M_\sigma(\beta)^{k} \cdot \tfrac{L}{2}$ & $\tau(\sigma,\beta,k)$ & $\tfrac{\tau(\sigma, \beta, k)}{\overline M_\sigma(\beta)^{k} \cdot L/2}$ & $M_\sigma^*(\beta)$ & $\tfrac{\tau(\sigma, \beta, k)}{M_\sigma^*(\beta)^{k} \cdot L/2}$ \\
            \hline\hline
            \endhead
            \hline
            \endfoot
            \hline
            \endlastfoot
            $0.100$ & $0.1732$ & $8$  & $1.773 \times 10^{-1}$ & $6.064 \times 10^{-3}$ & $0.0342$ & $0.8784$ & $0.0342$ \\
            $0.100$ & $0.1732$ & $16$ & $6.284 \times 10^{-2}$ & $6.891 \times 10^{-4}$ & $0.0110$ & $0.8784$ & $0.0110$ \\
            $0.100$ & $0.1732$ & $32$ & $7.898 \times 10^{-3}$ & $1.095 \times 10^{-5}$ & $0.0014$ & $0.8784$ & $0.0014$ \\
            $0.100$ & $0.3463$ & $8$  & $1.253 \times 10^{-1}$ & $3.439 \times 10^{-3}$ & $0.0274$ & $0.8405$ & $0.0276$ \\
            $0.100$ & $0.3463$ & $16$ & $3.139 \times 10^{-2}$ & $2.029 \times 10^{-4}$ & $0.0065$ & $0.8405$ & $0.0065$ \\
            $0.100$ & $0.3463$ & $32$ & $1.971 \times 10^{-3}$ & $3.469 \times 10^{-6}$ & $0.0018$ & $0.8405$ & $0.0018$ \\
            \hline
            $0.050$ & $0.2115$ & $8$  & $2.954 \times 10^{-1}$ & $1.215 \times 10^{-2}$ & $0.0411$ & $0.9363$ & $0.0411$ \\
            $0.050$ & $0.2115$ & $16$ & $1.745 \times 10^{-1}$ & $3.223 \times 10^{-3}$ & $0.0185$ & $0.9363$ & $0.0185$ \\
            $0.050$ & $0.2115$ & $32$ & $6.093 \times 10^{-2}$ & $3.529 \times 10^{-4}$ & $0.0058$ & $0.9363$ & $0.0058$ \\
            $0.050$ & $0.4230$ & $8$  & $2.362 \times 10^{-1}$ & $7.746 \times 10^{-3}$ & $0.0328$ & $0.9103$ & $0.0328$ \\
            $0.050$ & $0.4230$ & $16$ & $1.115 \times 10^{-1}$ & $1.393 \times 10^{-3}$ & $0.0125$ & $0.9103$ & $0.0125$ \\
            $0.050$ & $0.4230$ & $32$ & $2.488 \times 10^{-2}$ & $6.592 \times 10^{-5}$ & $0.0026$ & $0.9103$ & $0.0027$ \\
            \hline
            $0.025$ & $0.2423$ & $8$  & $3.820 \times 10^{-1}$ & $1.640 \times 10^{-2}$ & $0.0429$ & $0.9669$ & $0.0429$ \\
            $0.025$ & $0.2423$ & $16$ & $2.919 \times 10^{-1}$ & $6.191 \times 10^{-3}$ & $0.0212$ & $0.9669$ & $0.0212$ \\
            $0.025$ & $0.2423$ & $32$ & $1.704 \times 10^{-1}$ & $1.590 \times 10^{-3}$ & $0.0093$ & $0.9669$ & $0.0093$ \\
            $0.025$ & $0.4846$ & $8$  & $3.325 \times 10^{-1}$ & $1.096 \times 10^{-2}$ & $0.0330$ & $0.9502$ & $0.0330$ \\
            $0.025$ & $0.4846$ & $16$ & $2.211 \times 10^{-1}$ & $3.257 \times 10^{-3}$ & $0.0147$ & $0.9502$ & $0.0147$ \\
            $0.025$ & $0.4846$ & $32$ & $9.777 \times 10^{-2}$ & $5.267 \times 10^{-4}$ & $0.0054$ & $0.9502$ & $0.0054$ \\
        \end{longtable}
    }
    {
        \footnotesize
        \captionsetup{font=normalsize}
        \setlength{\tabcolsep}{4pt}
        \begin{longtable}{|c|c|c|c|c|c|c|c|}
            \caption{PEP verification for $\beta^\star < \beta \le 1 - \sigma$}\label{table:pep-verification-b2}\\
            \hline
            \rule{0pt}{2.8ex}$\sigma$ & $\beta$ & $k$ & $\overline M_\sigma(\beta)^{k} \cdot \tfrac{L}{2}$ & $\tau(\sigma,\beta,k)$ & $\tfrac{\tau(\sigma,\beta,k)}{\overline M_\sigma(\beta)^{k} \cdot L/2}$ & $M_\sigma^*(\beta)$ & $\tfrac{\tau(\sigma,\beta,k)}{M_\sigma^*(\beta)^{k} \cdot L/2}$ \\
            \hline\hline
            \endfirsthead
            \hline
            \rule{0pt}{2.8ex}$\sigma$ & $\beta$ & $k$ & $\overline M_\sigma(\beta)^{k} \cdot \tfrac{L}{2}$ & $\tau(\sigma,\beta,k)$ & $\tfrac{\tau(\sigma,\beta,k)}{\overline M_\sigma(\beta)^{k} \cdot L/2}$ & $M_\sigma^*(\beta)$ & $\tfrac{\tau(\sigma,\beta,k)}{M_\sigma^*(\beta)^{k} \cdot L/2}$ \\
            \hline\hline
            \endhead
            \hline
            \endfoot
            \hline
            \endlastfoot
            $0.100$ & $0.6463$ & $8$  & $5.724 \times 10^{-2}$ & $4.002 \times 10^{-3}$ & $0.0699$ & $0.7627$ & $0.0699$ \\
            $0.100$ & $0.6463$ & $16$ & $6.552 \times 10^{-3}$ & $8.986 \times 10^{-5}$ & $0.0137$ & $0.7627$ & $0.0137$ \\
            $0.100$ & $0.6463$ & $32$ & $8.586 \times 10^{-5}$ & NaN & NaN & $0.7627$ & NaN \\
            $0.100$ & $0.7732$ & $8$  & $1.172 \times 10^{-1}$ & $9.171 \times 10^{-3}$ & $0.0782$ & $0.8342$ & $0.0782$ \\
            $0.100$ & $0.7732$ & $16$ & $2.749 \times 10^{-2}$ & $5.858 \times 10^{-4}$ & $0.0213$ & $0.8342$ & $0.0213$ \\
            $0.100$ & $0.7732$ & $32$ & $1.512 \times 10^{-3}$ & NaN & NaN & $0.8342$ & NaN \\
            $0.100$ & $0.9000 = 1-\sigma$ & $8$  & $2.152 \times 10^{-1}$ & $2.251 \times 10^{-2}$ & $0.1046$ & $0.9000$ & $0.1046$ \\
            $0.100$ & $0.9000 = 1-\sigma$ & $16$ & $9.265 \times 10^{-2}$ & $3.905 \times 10^{-3}$ & $0.0422$ & $0.9000$ & $0.0422$ \\
            $0.100$ & $0.9000 = 1-\sigma$ & $32$ & $1.717 \times 10^{-2}$ & NaN & NaN & $0.9000$ & NaN \\
            \hline
            $0.050$ & $0.7397$ & $8$  & $1.219 \times 10^{-1}$ & $1.621 \times 10^{-2}$ & $0.1330$ & $0.8383$ & $0.1330$ \\
            $0.050$ & $0.7397$ & $16$ & $2.973 \times 10^{-2}$ & $1.227 \times 10^{-3}$ & $0.0413$ & $0.8383$ & $0.0413$ \\
            $0.050$ & $0.7397$ & $32$ & $1.768 \times 10^{-3}$ & $1.238 \times 10^{-5}$ & $0.0070$ & $0.8383$ & $0.0070$ \\
            $0.050$ & $0.8448$ & $8$  & $2.074 \times 10^{-1}$ & $3.197 \times 10^{-2}$ & $0.1541$ & $0.8959$ & $0.1541$ \\
            $0.050$ & $0.8448$ & $16$ & $8.606 \times 10^{-2}$ & $5.082 \times 10^{-3}$ & $0.0591$ & $0.8959$ & $0.0591$ \\
            $0.050$ & $0.8448$ & $32$ & $1.481 \times 10^{-2}$ & $2.108 \times 10^{-4}$ & $0.0142$ & $0.8959$ & $0.0142$ \\
            $0.050$ & $0.9500 = 1-\sigma$ & $8$  & $3.317 \times 10^{-1}$ & $6.609 \times 10^{-2}$ & $0.1992$ & $0.9500$ & $0.1992$ \\
            $0.050$ & $0.9500 = 1-\sigma$ & $16$ & $2.201 \times 10^{-1}$ & $2.238 \times 10^{-2}$ & $0.1017$ & $0.9500$ & $0.1017$ \\
            $0.050$ & $0.9500 = 1-\sigma$ & $32$ & $9.686 \times 10^{-2}$ & $4.379 \times 10^{-3}$ & $0.0452$ & $0.9500$ & $0.0452$ \\
            \hline
            $0.025$ & $0.8096$ & $8$  & $1.941 \times 10^{-1}$ & $3.499 \times 10^{-2}$ & $0.1802$ & $0.8885$ & $0.1802$ \\
            $0.025$ & $0.8096$ & $16$ & $7.537 \times 10^{-2}$ & $6.782 \times 10^{-3}$ & $0.0900$ & $0.8885$ & $0.0900$ \\
            $0.025$ & $0.8096$ & $32$ & $1.136 \times 10^{-2}$ & $2.979 \times 10^{-4}$ & $0.0262$ & $0.8885$ & $0.0262$ \\
            $0.025$ & $0.8923$ & $8$  & $2.864 \times 10^{-1}$ & $5.865 \times 10^{-2}$ & $0.2048$ & $0.9327$ & $0.2048$ \\
            $0.025$ & $0.8923$ & $16$ & $1.641 \times 10^{-1}$ & $1.797 \times 10^{-2}$ & $0.1095$ & $0.9327$ & $0.1095$ \\
            $0.025$ & $0.8923$ & $32$ & $5.385 \times 10^{-2}$ & $2.438 \times 10^{-3}$ & $0.0453$ & $0.9327$ & $0.0453$ \\
            $0.025$ & $0.9750 = 1-\sigma$ & $8$  & $4.083 \times 10^{-1}$ & $1.001 \times 10^{-1}$ & $0.2452$ & $0.9750$ & $0.2452$ \\
            $0.025$ & $0.9750 = 1-\sigma$ & $16$ & $3.335 \times 10^{-1}$ & $5.361 \times 10^{-2}$ & $0.1608$ & $0.9750$ & $0.1608$ \\
            $0.025$ & $0.9750 = 1-\sigma$ & $32$ & $2.224 \times 10^{-1}$ & $2.445 \times 10^{-2}$ & $0.1099$ & $0.9750$ & $0.1099$ \\
        \end{longtable}
    }
    {
        \footnotesize
        \captionsetup{font=normalsize}
        \setlength{\tabcolsep}{4pt}
        \begin{longtable}{|c|c|c|c|c|c|c|c|}
            \caption{PEP verification for $\beta > 1 - \sigma$}\label{table:pep-verification-b3}\\
            \hline
            \rule{0pt}{2.8ex}$\sigma$ & $\beta$ & $k$ & $\overline M_\sigma(\beta)^{k} \cdot \tfrac{L}{2}$ & $\tau(\sigma,\beta,k)$ & $\tfrac{\tau(\sigma,\beta,k)}{\overline M_\sigma(\beta)^{k} \cdot L/2}$ & $M_\sigma^*(\beta)$ & $\tfrac{\tau(\sigma,\beta,k)}{M_\sigma^*(\beta)^{k} \cdot L/2}$ \\
            \hline\hline
            \endfirsthead
            \hline
            \rule{0pt}{2.8ex}$\sigma$ & $\beta$ & $k$ & $\overline M_\sigma(\beta)^{k} \cdot \tfrac{L}{2}$ & $\tau(\sigma,\beta,k)$ & $\tfrac{\tau(\sigma,\beta,k)}{\overline M_\sigma(\beta)^{k} \cdot L/2}$ & $M_\sigma^*(\beta)$ & $\tfrac{\tau(\sigma,\beta,k)}{M_\sigma^*(\beta)^{k} \cdot L/2}$ \\
            \hline\hline
            \endhead
            \hline
            \endfoot
            \hline
            \endlastfoot
            $0.100$ & $0.9333$ & $8$  & $2.878 \times 10^{-1}$ & $2.868 \times 10^{-2}$ & $0.0997$ & $0.9333$ & $0.0997$ \\
            $0.100$ & $0.9333$ & $16$ & $1.657 \times 10^{-1}$ & $6.414 \times 10^{-3}$ & $0.0387$ & $0.9333$ & $0.0387$ \\
            $0.100$ & $0.9333$ & $32$ & $5.491 \times 10^{-2}$ & $3.919 \times 10^{-4}$ & $0.0071$ & $0.9333$ & $0.0071$ \\
            $0.100$ & $0.9667$ & $8$  & $3.813 \times 10^{-1}$ & $3.663 \times 10^{-2}$ & $0.0961$ & $0.9667$ & $0.0961$ \\
            $0.100$ & $0.9667$ & $16$ & $2.908 \times 10^{-1}$ & $1.059 \times 10^{-2}$ & $0.0364$ & $0.9667$ & $0.0364$ \\
            $0.100$ & $0.9667$ & $32$ & $1.692 \times 10^{-1}$ & $1.062 \times 10^{-3}$ & $0.0063$ & $0.9667$ & $0.0063$ \\
            \hline
            $0.050$ & $0.9667$ & $8$  & $3.813 \times 10^{-1}$ & $7.451 \times 10^{-2}$ & $0.1954$ & $0.9667$ & $0.1954$ \\
            $0.050$ & $0.9667$ & $16$ & $2.908 \times 10^{-1}$ & $2.834 \times 10^{-2}$ & $0.0974$ & $0.9667$ & $0.0974$ \\
            $0.050$ & $0.9667$ & $32$ & $1.692 \times 10^{-1}$ & $7.163 \times 10^{-3}$ & $0.0423$ & $0.9667$ & $0.0423$ \\
            $0.050$ & $0.9833$ & $8$  & $4.370 \times 10^{-1}$ & $8.409 \times 10^{-2}$ & $0.1924$ & $0.9833$ & $0.1924$ \\
            $0.050$ & $0.9833$ & $16$ & $3.819 \times 10^{-1}$ & $3.581 \times 10^{-2}$ & $0.0938$ & $0.9833$ & $0.0938$ \\
            $0.050$ & $0.9833$ & $32$ & $2.917 \times 10^{-1}$ & $1.168 \times 10^{-2}$ & $0.0400$ & $0.9833$ & $0.0400$ \\
            \hline
            $0.025$ & $0.9833$ & $8$  & $4.370 \times 10^{-1}$ & $1.057 \times 10^{-1}$ & $0.2420$ & $0.9833$ & $0.2420$ \\
            $0.025$ & $0.9833$ & $16$ & $3.819 \times 10^{-1}$ & $6.005 \times 10^{-2}$ & $0.1572$ & $0.9833$ & $0.1572$ \\
            $0.025$ & $0.9833$ & $32$ & $2.917 \times 10^{-1}$ & $3.125 \times 10^{-2}$ & $0.1071$ & $0.9833$ & $0.1071$ \\
            $0.025$ & $0.9917$ & $8$  & $4.677 \times 10^{-1}$ & $1.118 \times 10^{-1}$ & $0.2390$ & $0.9917$ & $0.2390$ \\
            $0.025$ & $0.9917$ & $16$ & $4.376 \times 10^{-1}$ & $6.740 \times 10^{-2}$ & $0.1540$ & $0.9917$ & $0.1540$ \\
            $0.025$ & $0.9917$ & $32$ & $3.829 \times 10^{-1}$ & $4.009 \times 10^{-2}$ & $0.1047$ & $0.9917$ & $0.1047$ \\
        \end{longtable}
    }

\section{Conclusion and Future Works}\label{sec:conclusion}
    For Nesterov's constant momentum accelerated proximal gradient method we resolve the literature's gap on whether convergence is possible away from the classic momentum extrapolation $\tfrac{1 - \sqrt{\sigma}}{1 + \sqrt{\sigma}}$ proposed by Nesterov and Beck \cite{nesterov_lectures_2018,beck_first-order_2017}.
    We prove that the method converges linearly at every momentum parameter in the interval $(0, 1)$.
    We also give an explicit closed form bound on the convergence rate that holds throughout this range.
    \par
    For the heavy ball method, Goujaud et al. \cite{goujaud_provable_2025} prove that every parameter choice either fails to accelerate or fails to converge.
    Our results prove that Nesterov's method converges linearly throughout the range of extrapolation parameters in $(0, 1)$.
    They therefore complement the existing literature.
    We corroborate this convergence rate numerically via performance estimation.
    \par
    Aujol et al. \cite[Theorem 1]{aujol_heavy_2025} prove linear convergence of V-FISTA under the weaker quadratic growth condition without assuming strong convexity.
    Their momentum coefficient is given by a formula that requires knowing the modulus of strong convexity and the Lipschitz smoothness constant.
    Whether a similar result holds under this weaker condition for the range of momentum parameters we study is an open question.

\section*{Acknowledgements}
    Hongda Li and Xianfu Wang are supported by NSERC Discovery Grants.

\section*{Declaration of generative AI use}
    The authors used Claude Code Sonnet 5.5 for copyediting, language refinement, and coding assistance.
    The authors checked the terms of use, and confirm the originality, accuracy of the content.
    We take full responsibility for its integrity, including references.

\appendix

\section{Algebraic Identities for the Abstract Extrapolation Form}\label{sec:appendix-abstract-momentum-form}
    \subsection{Proof of Lemma \ref{lemma:abstract-dvt}}\label{ssec:abstract-dvt-pf}
        \begin{proof}
            \todoinline{
                Lemma \ref{lemma:abstract-dvt} proof depends on:
                \begin{itemize}[nosep]
                    \item Assumption \ref{ass:abst-nes-extrp}\ref{ass:abst-nes-extrp-lambdak}.
                    \item Assumption \ref{ass:abst-nes-extrp}\ref{ass:abst-nes-extrp-x}.
                    \item Assumption \ref{ass:abst-nes-extrp}\ref{ass:abst-nes-extrp-xd}.
                \end{itemize}
            }
            Assumption \ref{ass:abst-nes-extrp}\ref{ass:abst-nes-extrp-lambdak} and Assumption \ref{ass:abst-nes-extrp}\ref{ass:abst-nes-extrp-x} give $(1 - q_k)x_k = (\alpha_k - q_k)x_{k - 1}^\dagger + (1 - \alpha_k)x_{k - 1}^+$.
            Rearranging it gives $\alpha_k x_{k - 1}^\dagger = (1 - q_k)x_k + q_k x_{k - 1}^\dagger - (1 - \alpha_k)x_{k - 1}^+$ so
            \begin{align*}
                \alpha_k(x_{k - 1}^\dagger - \bar x)
                &\underset{(1)}{=} (1 - q_k)x_k + q_k x_{k - 1}^\dagger - (1 - \alpha_k)x_{k - 1}^+ - \alpha_k\bar x
                \\
                &= x_k - \underbrace{\left(\alpha_k \bar x + (1 - \alpha_k)x_{k - 1}^+\right)}_{z_k} + q_k(x_{k - 1}^\dagger - x_k)
                \\
                &= x_k - z_k + q_k(x_{k - 1}^\dagger - x_k).
            \end{align*}
            At (1) we substituted $\alpha_k x_{k - 1}^\dagger$ using the rearranged recursion above.
            \par
            Similarly substituting $x_k^\dagger = x_{k - 1}^+ + \alpha_k^{-1}(x_k^+ - x_{k - 1}^+)$ from Assumption \ref{ass:abst-nes-extrp}\ref{ass:abst-nes-extrp-xd} gives
            \begin{align*}
                \alpha_k(x_k^\dagger - \bar x)
                &= \alpha_k x_{k - 1}^+ + x_k^+ - x_{k - 1}^+ - \alpha_k \bar x
                \\
                &= x_k^+ - \left(\alpha_k \bar x + (1 - \alpha_k)x_{k - 1}^+\right)
                \\
                &= x_k^+ - z_k.
            \end{align*}
        \end{proof}
    \subsection{Proof of Lemma \ref{lemma:abstract-dist-splt}}\label{ssec:abstract-dist-splt-pf}
        \begin{proof}
            \todoinline{
                Lemma \ref{lemma:abstract-dist-splt} proof depends on:
                \begin{itemize}[nosep]
                    \item Assumption \ref{ass:abst-nes-extrp}\ref{ass:abst-nes-extrp-lambdak}.
                    \item Assumption \ref{ass:abst-nes-extrp}\ref{ass:abst-nes-extrp-x}.
                    \item Lemma \ref{lemma:ecld-id-extd}.
                \end{itemize}
            }
            By Assumption \ref{ass:abst-nes-extrp}\ref{ass:abst-nes-extrp-x}, $x_k = \lambda_k x_{k - 1}^\dagger + (1 - \lambda_k)x_{k - 1}^+$ where $\lambda_k$ is given by Assumption \ref{ass:abst-nes-extrp}\ref{ass:abst-nes-extrp-lambdak}. It follows that
            \begin{align*}
                \Vert \bar x - x_k\Vert^2
                &= \Vert \bar x - \lambda_k x_{k - 1}^\dagger - (1 - \lambda_k)x_{k - 1}^+\Vert^2
                \\
                &\underset{(1)}{=} (1 - \lambda_k)\Vert \bar x - x_{k - 1}^+\Vert^2
                + \lambda_k\Vert \bar x - x_{k - 1}^\dagger\Vert^2
                - \lambda_k(1 - \lambda_k)\Vert x_{k - 1}^\dagger - x_{k - 1}^+\Vert^2.
            \end{align*}
            At (1) we applied Lemma \ref{lemma:ecld-id-extd} with $(w, x, \bar x, \lambda) = (\bar x, x_{k - 1}^+, x_{k - 1}^\dagger, \lambda_k)$.
        \end{proof}
    \subsection{Proof of Lemma \ref{lemma:abstract-gp-splt}}\label{ssec:abstract-gp-splt-pf}
        \begin{proof}
            \todoinline{
                Lemma \ref{lemma:abstract-gp-splt} proof depends on:
                \begin{itemize}[nosep]
                    \item Assumption \ref{ass:abst-nes-extrp}\ref{ass:abst-nes-extrp-lambdak}.
                    \item Assumption \ref{ass:abst-nes-extrp}\ref{ass:abst-nes-extrp-x}.
                    \item Lemma \ref{lemma:abstract-dist-splt}.
                \end{itemize}
            }
            We have
            \begin{align*}
                \frac{1}{2}\Vert x_k - z_k\Vert^2
                &= \frac{1}{2}\Vert \alpha_k(x_{k - 1}^\dagger - \bar x) - q_k(x_{k - 1}^\dagger - x_k) \Vert^2
                \\
                &= \frac{\alpha_k^2}{2}\Vert x_{k - 1}^\dagger - \bar x\Vert^2
                + \frac{q_k^2}{2}\Vert x_{k - 1}^\dagger - x_k\Vert^2
                - \alpha_k q_k\langle x_{k - 1}^\dagger - \bar x, x_{k - 1}^\dagger - x_k \rangle
                \\
                &\underset{(1)}{=}
                \frac{\alpha_k^2 - \alpha_k q_k}{2}\Vert x_{k - 1}^\dagger - \bar x\Vert^2
                + \frac{q_k^2 - \alpha_k q_k}{2}\Vert x_{k - 1}^\dagger - x_k\Vert^2
                + \frac{\alpha_k q_k}{2}\Vert x_k - \bar x \Vert^2
                \\
                &\underset{(2)}{=}
                \frac{\alpha_k\lambda_k(1 - q_k)}{2}\Vert x_{k - 1}^\dagger - \bar x\Vert^2
                - \frac{q_k\lambda_k(1 - q_k)(1 - \lambda_k)^2}{2}\Vert x_{k - 1}^\dagger - x_{k - 1}^+\Vert^2
                \\&\quad
                + \frac{\alpha_kq_k}{2}\Vert x_k - \bar x \Vert^2.
            \end{align*}
            At (1) we used $\langle u, v\rangle = \tfrac{1}{2}\left(\Vert u\Vert^2 + \Vert v\Vert^2 - \Vert u - v\Vert^2\right)$ with $u = x_{k - 1}^\dagger - \bar x$ and $v = x_{k - 1}^\dagger - x_k$ so that $u - v = x_k - \bar x$.
            At (2) we substituted $\alpha_k - q_k = \lambda_k(1 - q_k)$ from Assumption \ref{ass:abst-nes-extrp}\ref{ass:abst-nes-extrp-lambdak} and used the relation $x_k = \lambda_k x_{k - 1}^\dagger + (1 - \lambda_k)x_{k - 1}^+$ of Assumption \ref{ass:abst-nes-extrp}\ref{ass:abst-nes-extrp-x}.
            This gives $x_{k - 1}^\dagger - x_k = (1 - \lambda_k)(x_{k - 1}^\dagger - x_{k - 1}^+)$.
            \par
            It remains to expand $\Vert x_k - \bar x\Vert^2$. By Lemma \ref{lemma:abstract-dist-splt} scaled by $\tfrac{\alpha_k q_k}{2}$
            \begin{align*}
                \frac{\alpha_kq_k}{2}\Vert x_k - \bar x \Vert^2
                &= \frac{\alpha_kq_k(1 - \lambda_k)}{2}\left\Vert \bar x - x_{k - 1}^+ \right\Vert^2
                + \frac{\alpha_kq_k\lambda_k}{2}\Vert \bar x - x_{k - 1}^\dagger \Vert^2
                \\&\quad
                - \frac{\alpha_kq_k\lambda_k(1 - \lambda_k)}{2}\left\Vert x_{k - 1}^\dagger - x_{k - 1}^+ \right\Vert^2.
            \end{align*}
            Substituting this back to equality at (2) the coefficient of $\Vert x_{k - 1}^\dagger - \bar x\Vert^2$ becomes
            \begin{align*}
                \frac{\alpha_k\lambda_k(1 - q_k)}{2} + \frac{\alpha_kq_k\lambda_k}{2}
                = \frac{\alpha_k\lambda_k}{2}\big[(1 - q_k) + q_k\big]
                = \frac{\alpha_k\lambda_k}{2},
            \end{align*}
            and the coefficient of $\Vert x_{k - 1}^\dagger - x_{k - 1}^+\Vert^2$ collapses via $(1 - q_k)(1 - \lambda_k) + \alpha_k = (1 - \alpha_k) + \alpha_k = 1$ where the first equality is Assumption \ref{ass:abst-nes-extrp}\ref{ass:abst-nes-extrp-lambdak}:
            \begin{align*}
                - \frac{q_k\lambda_k(1 - q_k)(1 - \lambda_k)^2}{2} - \frac{\alpha_kq_k\lambda_k(1 - \lambda_k)}{2}
                &= -\frac{q_k\lambda_k(1 - \lambda_k)}{2}\big[(1 - q_k)(1 - \lambda_k) + \alpha_k\big]
                \\
                &= -\frac{q_k\lambda_k(1 - \lambda_k)}{2}.
            \end{align*}
            The term $\Vert \bar x - x_{k - 1}^+\Vert^2$ arises only from the expansion above with coefficient $\tfrac{\alpha_kq_k(1 - \lambda_k)}{2}$.
        \end{proof}

\section{\texorpdfstring{Auxiliary Results for the Proof of Theorem \ref{thm:lin-cnvg-constn-mmntm}}{Auxiliary Results for the Proof of the Linear Convergence Theorem}}\label{sec:appendix-lin-cnvg-aux}
    \subsection{Proof of Lemma \ref{lemma:dbeta-equiv}}\label{ssec:dbeta-equiv-pf}
        \begin{proof}
            \todoinline{
                Lemma \ref{lemma:dbeta-equiv} proof depends on:
                \begin{itemize}[nosep]
                    \item \eqref{def:dbeta-set}.
                \end{itemize}
            }
            Fix any $\beta \in (0, 1)$.
            We first show the $\subseteq$ relation.
            Let $(\alpha, q) \in \mathcal D_\beta$.
            By \eqref{def:dbeta-set},
            \begin{align*}
                \beta &= \left(
                    \tfrac{\alpha - q}{1 - q}
                \right)
                \tfrac{1 - \alpha}{\alpha}
                \\
                &=
                \left(
                    1 - \tfrac{1 - \alpha}{1 - q}
                \right)
                \tfrac{1 - \alpha}{\alpha}
                \\
                &\le
                \left(
                    1 - \inf_{q \in [0, \alpha)} \tfrac{1 - \alpha}{1 - q}
                \right)
                \tfrac{1 - \alpha}{\alpha}
                \\
                &=\alpha \cdot \tfrac{1 - \alpha}{\alpha}  = 1 - \alpha.
            \end{align*}
            Therefore $\beta \le 1 - \alpha$, that is, $0 < \alpha \le 1 - \beta$.
            Next we verify the relationship on $q$.
            Using $\alpha \in (0, 1 - \beta]$ and the equality we presented above we solve it for $q$ in terms of $\alpha, \beta$:
            \begin{align*}
                (\alpha, q) \in \mathcal D_\beta
                \implies&
                \beta = \left(
                    \tfrac{\alpha - q}{1 - q}
                \right)
                \tfrac{1 - \alpha}{\alpha}
                \\
                \iff&
                \alpha(1 - q)\beta = (\alpha - q)(1 - \alpha)
                \\
                \iff&
                \alpha \beta - \alpha \beta q = \alpha - q - \alpha^2 + \alpha q
                \\
                \underset{(1)}\iff&
                \alpha\beta - \alpha + \alpha^2 = q(\alpha\beta + \alpha - 1)
                \\
                \iff&
                q = \tfrac{\alpha\beta - \alpha + \alpha^2}{\alpha\beta + \alpha - 1}
                = \alpha \cdot \tfrac{1 - \alpha - \beta}{1 - \alpha - \alpha\beta}.
            \end{align*}
            At (1) observe that $\alpha\beta + \alpha - 1 \le \alpha\beta - \beta = \beta(\alpha - 1) < 0$ because for all $\beta \in (0, 1)$ it follows that $\alpha \le 1 - \beta < 1$.
            Hence we can divide by $\alpha\beta + \alpha - 1$.
            At this point we have shown the $\subseteq$ relationship.
            \par
            Next we show the $\supseteq$ relationship.
            Given any $\alpha \in (0, 1 - \beta]$ we will show that $q \in [0, \alpha]\setminus \{1\}$ as a consequence.
            We claim $q < 1$ because:
            \begin{align*}
                0 \le q &= \alpha\cdot\tfrac{1 - \alpha - \beta}{1 - \alpha(1 + \beta)} \le \alpha \le 1 - \beta < 1.
            \end{align*}
            The last inequality is due to the fact that $\beta > 0$.
            The second last inequality is due to $1 - \alpha - \beta \le 1 - \alpha - \alpha \beta$ because $\beta \in (0, 1)$.
            The edge case $\beta = 0$ forces either $q = \alpha$ or $\alpha = 1$ giving $\mathcal D_0 = \{(\alpha, \alpha): \alpha \in (0, 1)\}\cup \{1\}\times [0, 1)$.
        \end{proof}
    \subsection{Proof of Lemma \ref{lemma:m-in-alpha-beta}}\label{ssec:m-in-alpha-beta-pf}
        \begin{proof}
            \todoinline{
                Lemma \ref{lemma:m-in-alpha-beta} proof depends on:
                \begin{itemize}[nosep]
                    \item \eqref{def:dbeta-set}.
                    \item \eqref{def:c-sigma-alpha-q}.
                    \item Lemma \ref{lemma:dbeta-equiv}.
                \end{itemize}
                The proof implicitly depends on:
                \begin{itemize}[nosep]
                    \item \eqref{def:xi-sigma-alpha-q}.
                \end{itemize}
            }
            Fix any $\beta \in (0, 1)$ and $\sigma \in (0, 1)$.
            For all $(\alpha, q) \in \mathcal D_\beta$, \eqref{def:dbeta-set} asserts that $\beta = \left(\tfrac{\alpha - q}{1 - q}\right)\tfrac{1 - \alpha}{\alpha}$ and therefore substituting it into \eqref{def:c-sigma-alpha-q} yields:
            \begin{align*}
                C_\sigma(\alpha, q) = \frac{(1 - \sigma)(\alpha - q)}{\alpha(1 - q)}
                = \frac{1 - \sigma}{\alpha}\cdot\frac{\alpha - q}{1 - q}
                = \frac{(1 - \sigma)\beta}{1 - \alpha}.
            \end{align*}
            Using Lemma \ref{lemma:dbeta-equiv} we have
            \begin{align*}
                1 - q &= \frac{(1 - \alpha)^2}{1 - \alpha(1 + \beta)},
                \\
                \alpha - q &= \frac{\alpha\beta(1 - \alpha)}{1 - \alpha(1 + \beta)},
                \\
                q - \sigma &= \frac{(1 - \alpha)(\alpha - \sigma) - \alpha\beta(1 - \sigma)}{1 - \alpha(1 + \beta)}
            \end{align*}
            so that $\tfrac{q - \sigma}{1 - q}= \tfrac{(1 - \alpha)(\alpha - \sigma) - \alpha\beta(1 - \sigma)}{(1 - \alpha)^2}$.
            This enables an expression for $\xi_\sigma(\alpha, q)$ in terms of $\beta, \alpha$:
            \begin{align*}
                \xi_\sigma(\alpha, q) =
                \frac{\alpha(1 - \alpha)(q - \sigma)}{\sigma(1 - q)}
                &= \frac{\alpha(1 - \alpha)}{\sigma}\cdot\frac{(1 - \alpha)(\alpha - \sigma) - \alpha\beta(1 - \sigma)}{(1 - \alpha)^2}
                \\
                &= \frac{\alpha}{\sigma}\left(
                    (\alpha - \sigma) - \alpha\cdot\frac{(1 - \sigma)\beta}{1 - \alpha}
                \right).
            \end{align*}
        \end{proof}
    \subsection{Proof of Lemma \ref{lemma:xi-sigma-beta-nonpos}}\label{ssec:xi-sigma-beta-nonpos-pf}
        \begin{proof}
            \todoinline{
                Lemma \ref{lemma:xi-sigma-beta-nonpos} proof depends on:
                \begin{itemize}[nosep]
                    \item \eqref{def:xi-sigma-beta-in-alpha}.
                    \item Lemma \ref{lemma:m-in-alpha-beta}.
                \end{itemize}
            }
            Fix any $\beta \in (0, 1)$.
            We solve the inequality to produce the result.
            For $\alpha \in (0, 1 - \beta]$, \eqref{def:xi-sigma-beta-in-alpha} and Lemma \ref{lemma:m-in-alpha-beta} give:
            \begin{align}\label{eqn:xi-sigma-beta-nonpos-pitem1}\begin{split}
                0 &\ge  \tfrac{\sigma}{\alpha}(1 - \alpha)\xi_{\sigma, \beta}(\alpha)
                \\
                & =
                (1 - \alpha)(\alpha - \sigma) - \alpha\beta(1 - \sigma)
                \\
                &= \alpha - \sigma - \alpha^2 + \alpha \sigma - \alpha \beta + \alpha\beta\sigma
                \\
                &= - \alpha^2 + \alpha(1 + \sigma - \beta + \beta\sigma) - \sigma
                \\
                &= - \alpha^2 + \alpha(1 + \sigma - \beta(1 - \sigma)) - \sigma.
            \end{split}
            \end{align}
            The last expression is a negative quadratic function of $\alpha$.
            It will be less than zero whenever its discriminant is negative.
            Its discriminant is given by:
            \begin{align*}
                0&\ge  (1 + \sigma - \beta(1 - \sigma))^2 - 4\sigma
                \\
                &=
                \left(1 + \sigma - \beta(1 - \sigma) + 2\sqrt{\sigma}\right)
                \left(1 + \sigma - \beta(1 - \sigma) - 2\sqrt{\sigma}\right)
                \\
                &=
                \left((1 + \sqrt{\sigma})^2 - \beta(1 - \sigma)\right)
                \left((1 - \sqrt{\sigma})^2 - \beta(1 - \sigma)\right)
                \\
                &=
                (1 + \sqrt{\sigma})(1 - \sqrt{\sigma})
                \left((1 + \sqrt{\sigma}) - \beta(1 - \sqrt{\sigma})\right)
                \left((1 - \sqrt{\sigma}) - \beta(1 + \sqrt{\sigma})\right).
            \end{align*}
            Observe that only the last term in the product can be negative.
            The first two terms are nonnegative because $\sigma \in (0, 1)$.
            By the hypothesis that $\beta \in (0, 1)$ the third term is nonnegative: $1 + \sqrt{\sigma} - \beta(1 - \sqrt{\sigma}) \ge 2\sqrt{\sigma} > 0$.
            Therefore the last term introduces a sign constraint on $\beta$ which is:
            \begin{align*}
                & 0 \ge  (1 - \sqrt{\sigma}) - \beta(1 + \sqrt{\sigma})
                \iff \beta \ge \tfrac{1 - \sqrt{\sigma}}{1 + \sqrt{\sigma}}.
            \end{align*}
        \end{proof}

\bibliographystyle{siam}
\bibliography{references/beta_momentum,references/refs,references/software_tools}

@book{nesterov_lectures_2018,
	series = {Springer {Optimization} and {Its} {Applications}},
	title = {Lectures on {Convex} {Optimization}},
	volume = {137},
	isbn = {978-3-319-91577-7 978-3-319-91578-4},
	url = {http://link.springer.com/10.1007/978-3-319-91578-4},
	urldate = {2023-10-11},
	publisher = {Springer International Publishing},
	author = {Nesterov, Yurii},
	year = {2018},
}

@book{beck_first-order_2017,
	series = {{MOS}-{SIAM} {Series} in {Optimization}},
	title = {First-order {Methods} in {Optimization}},
	isbn = {978-1-61197-498-0},
	url = {https://epubs.siam.org/doi/book/10.1137/1.9781611974997},
	language = {en},
	urldate = {2023-10-19},
	publisher = {SIAM},
	author = {Beck, Amir},
	year = {2017},
}

@book{bauschke_convex_2017,
	address = {Cham},
	series = {{CMS} {Books} in {Mathematics}},
	title = {Convex {Analysis} and {Monotone} {Operator} {Theory} in {Hilbert} {Spaces}},
	isbn = {978-3-319-48310-8 978-3-319-48311-5},
	url = {https://link.springer.com/10.1007/978-3-319-48311-5},
	language = {en},
	urldate = {2023-11-29},
	publisher = {Springer International Publishing},
	author = {Bauschke, Heinz H. and Combettes, Patrick L.},
	year = {2017},
}

@article{aujol_optimal_2019,
	title = {Optimal convergence rates for nesterov acceleration},
	volume = {29},
	url = {https://hal.science/hal-01786117},
	doi = {10.1137/18M1186757},
	number = {4},
	urldate = {2024-06-25},
	journal = {SIAM Journal on Optimization},
	author = {Aujol, Jean François and Dossal, Charles H and Rondepierre, Aude},
	month = dec,
	year = {2019},
	pages = {3131--3153},
}

@misc{alamo_restart_2019,
	title = {Restart {FISTA} with global linear convergence},
	url = {http://arxiv.org/abs/1906.09126},
	doi = {10.48550/arXiv.1906.09126},
	urldate = {2023-10-18},
	publisher = {arXiv},
	author = {Alamo, Teodoro and Krupa, Pablo and Limon, Daniel},
	month = dec,
	year = {2019},
}

@inproceedings{alamo_gradient_2019,
	title = {Gradient based restart {FISTA}},
	url = {https://ieeexplore.ieee.org/document/9029983},
	doi = {10.1109/CDC40024.2019.9029983},
	urldate = {2023-10-18},
	booktitle = {2019 {IEEE} 58th {Conference} on {Decision} and {Control} ({CDC})},
	author = {Alamo, Teodoro and Krupa, Pablo and Limon, Daniel},
	month = dec,
	year = {2019},
	pages = {3936--3941},
}

@article{guler_new_1992,
	title = {New proximal point algorithms for convex minimization},
	volume = {2},
	issn = {1052-6234},
	url = {https://epubs.siam.org/doi/10.1137/0802032},
	doi = {10.1137/0802032},
	number = {4},
	urldate = {2023-11-30},
	journal = {SIAM Journal on Optimization},
	author = {Güler, Osman},
	month = nov,
	year = {1992},
	pages = {649--664},
}

@article{su_differential_2016,
	title = {A differential equation for modeling nesterov's accelerated gradient method: {Theory} and {Insights}},
	volume = {17},
	issn = {1533-7928},
	shorttitle = {A {Differential} {Equation} for {Modeling} {Nesterov}'s {Accelerated} {Gradient} {Method}},
	url = {http://jmlr.org/papers/v17/15-084.html},
	language = {en},
	number = {153},
	urldate = {2023-10-09},
	journal = {Journal of Machine Learning Research},
	author = {Su, Weijie and Boyd, Stephen and Candes, Emmanuel J.},
	year = {2016},
	pages = {1--43},
}

@article{aujol_parameter-free_2024,
	title = {Parameter-{Free} {FISTA} by adaptive restart and backtracking},
	url = {https://epubs.siam.org/doi/10.1137/23M158961X},
	language = {en},
	urldate = {2023-10-09},
	journal = {SIAM Journal on Optimization},
	author = {Aujol, Jean-François and Calatroni, Luca and Dossal, Charles and Labarrière, Hippolyte and Rondepierre, Aude},
	month = dec,
	year = {2024},
	pages = {3259--3285},
}

@article{attouch_first-order_2022,
	title = {First-order optimization algorithms via inertial systems with {Hessian} driven damping},
	volume = {193},
	issn = {1436-4646},
	url = {https://doi.org/10.1007/s10107-020-01591-1},
	doi = {10.1007/s10107-020-01591-1},
	language = {en},
	number = {1},
	urldate = {2024-12-29},
	journal = {Mathematical Programming},
	author = {Attouch, Hedy and Chbani, Zaki and Fadili, Jalal and Riahi, Hassan},
	month = may,
	year = {2022},
	pages = {113--155},
}

@article{necoara_linear_2019,
	title = {Linear convergence of first order methods for non-strongly convex optimization},
	volume = {175},
	issn = {1436-4646},
	url = {https://doi.org/10.1007/s10107-018-1232-1},
	doi = {10.1007/s10107-018-1232-1},
	language = {en},
	number = {1},
	urldate = {2023-10-11},
	journal = {Mathematical Programming},
	author = {Necoara, I. and Nesterov, Yu. and Glineur, F.},
	month = may,
	year = {2019},
	pages = {69--107},
}

@inproceedings{lee_geometric_2021,
	title = {A {Geometric} structure of acceleration and its role in making gradients small fast},
	volume = {34},
	url = {https://papers.nips.cc/paper/2021/hash/647c722bf90a49140184672e0d3723e3-Abstract.html},
	urldate = {2024-06-25},
	booktitle = {Advances in {Neural} {Information} {Processing} {Systems}},
	author = {Lee, Jongmin and Park, Chanwoo and Ryu, Ernest},
	year = {2021},
	pages = {11999--12012},
}

@article{lin_catalyst_2018,
	title = {Catalyst acceleration for first-order convex optimization: from theory to practice},
	volume = {18},
	url = {http://jmlr.org/papers/v18/17-748.html},
	number = {212},
	journal = {Journal of Machine Learning Research},
	author = {Lin, Hongzhou and Mairal, Julien and Harchaoui, Zaid},
	year = {2018},
	pages = {1--54},
}

@inproceedings{schmidt_convergence_2011,
	title = {Convergence rates of inexact proximal-gradient methods for convex optimization},
	volume = {24},
	url = {https://proceedings.neurips.cc/paper_files/paper/2011/hash/8f7d807e1f53eff5f9efbe5cb81090fb-Abstract.html},
	urldate = {2024-10-23},
	booktitle = {Advances in {Neural} {Information} {Processing} {Systems}},
	publisher = {Curran Associates, Inc.},
	author = {Schmidt, Mark and Roux, Nicolas and Bach, Francis},
	year = {2011},
}

@article{drori_performance_2014,
	title = {Performance of first-order methods for smooth convex minimization: a novel approach},
	volume = {145},
	issn = {0025-5610, 1436-4646},
	shorttitle = {Performance of first-order methods for smooth convex minimization},
	url = {https://doi.org/10.1007/s10107-013-0653-0},
	doi = {10.1007/s10107-013-0653-0},
	language = {en},
	number = {1-2},
	urldate = {2025-02-03},
	journal = {Mathematical Programming},
	publisher = {Springer},
	author = {Drori, Yoel and Teboulle, Marc},
	month = jun,
	year = {2014},
	pages = {451--482},
}

@article{odonoghue_adaptive_2015,
	title = {Adaptive restart for accelerated gradient schemes},
	volume = {15},
	issn = {1615-3383},
	url = {https://doi.org/10.1007/s10208-013-9150-3},
	doi = {10.1007/s10208-013-9150-3},
	language = {en},
	number = {3},
	urldate = {2025-04-06},
	journal = {Foundations of Computational Mathematics},
	author = {O’Donoghue, Brendan and Candès, Emmanuel},
	month = jun,
	year = {2015},
	pages = {715--732},
}

@article{villa_accelerated_2013,
	title = {Accelerated and inexact forward-backward algorithms},
	volume = {23},
	issn = {1052-6234},
	url = {https://epubs.siam.org/doi/10.1137/110844805},
	doi = {10.1137/110844805},
	number = {3},
	urldate = {2025-09-09},
	journal = {SIAM Journal on Optimization},
	author = {Villa, Silvia and Salzo, Saverio and Baldassarre, Luca and Verri, Alessandro},
	month = jan,
	year = {2013},
	pages = {1607--1633},
}

@article{bello-cruz_inexact_2020,
	title = {On inexact accelerated proximal gradient methods with relative error rules},
	url = {https://www.semanticscholar.org/paper/On-Inexact-Accelerated-Proximal-Gradient-Methods-Bello-Cruz-Gonccalves/b1a9b0502c3bdbb8a96516a9e29689649b505657},
	urldate = {2026-03-27},
	journal = {arXiv: Optimization and Control},
	author = {Bello-Cruz, Yunier and Gonccalves, Max L. N. and Krislock, Nathan},
	month = may,
	year = {2020},
}

@misc{aujol_heavy_2025,
	title = {Heavy ball momentum for non-strongly convex optimization},
	url = {http://arxiv.org/abs/2403.06930},
	doi = {10.48550/arXiv.2403.06930},
	urldate = {2026-04-27},
	publisher = {arXiv},
	author = {Aujol, Jean-François and Dossal, Charles and Labarrière, Hippolyte and Rondepierre, Aude},
	month = jan,
	year = {2025},
	note = {10.48550/arXiv.2403.06930
Prepint only.},
}

@inproceedings{kim_unifying_2023,
	title = {Unifying {Nesterov}’s accelerated gradient methods for convex and strongly convex objective functions},
	issn = {2640-3498},
	url = {https://proceedings.mlr.press/v202/kim23y.html},
	language = {en},
	urldate = {2026-04-27},
	booktitle = {Proceedings of the 40th {International} {Conference} on {Machine} {Learning}},
	publisher = {PMLR},
	author = {Kim, Jungbin and Yang, Insoon},
	month = jul,
	year = {2023},
	pages = {16897--16954},
}

@article{goujaud_provable_2025,
	title = {Provable non-accelerations of the heavy-ball method},
	issn = {1436-4646},
	url = {https://doi.org/10.1007/s10107-025-02269-2},
	doi = {10.1007/s10107-025-02269-2},
	language = {en},
	urldate = {2026-08-04},
	journal = {Mathematical Programming},
	author = {Goujaud, Baptiste and Taylor, Adrien and Dieuleveut, Aymeric},
	month = oct,
	year = {2025},
	note = {Online first; arXiv:2307.11291. https://doi.org/10.1007/s10107-025-02269-2},
}

@article{lessard_analysis_2016,
	title = {Analysis and design of optimization algorithms via integral quadratic constraints},
	volume = {26},
	issn = {1052-6234},
	url = {https://epubs.siam.org/doi/abs/10.1137/15M1009597},
	doi = {10.1137/15M1009597},
	number = {1},
	urldate = {2026-08-04},
	journal = {SIAM Journal on Optimization},
	publisher = {Society for Industrial and Applied Mathematics},
	author = {Lessard, Laurent and Recht, Benjamin and Packard, Andrew},
	month = jan,
	year = {2016},
	pages = {57--95},
}

@misc{taylor_smooth_2016,
	title = {Smooth strongly convex interpolation and exact worst-case performance of first-order methods},
	doi = {10.48550/arXiv.1502.05666},
	urldate = {2026-08-04},
	publisher = {arXiv},
	author = {Taylor, Adrien B. and Hendrickx, Julien M. and Glineur, François},
	month = oct,
	year = {2016},
	note = {arXiv:1502.05666 [math]},
}

@article{polyak_methods_1964,
	title = {Some methods of speeding up the convergence of iteration methods},
	volume = {4},
	issn = {0041-5553},
	url = {https://www.sciencedirect.com/science/article/pii/0041555364901375},
	doi = {10.1016/0041-5553(64)90137-5},
	number = {5},
	urldate = {2026-08-04},
	journal = {USSR Computational Mathematics and Mathematical Physics},
	author = {Polyak, B. T.},
	month = jan,
	year = {1964},
	pages = {1--17},
}

@article{goujaud_counter-examples_2023,
	title = {Counter-examples in first-order optimization: a constructive approach},
	volume = {7},
	issn = {2475-1456},
	shorttitle = {Counter-{Examples} in {First}-{Order} {Optimization}},
	url = {https://ieeexplore.ieee.org/document/10153415},
	doi = {10.1109/LCSYS.2023.3286277},
	urldate = {2026-08-07},
	journal = {IEEE Control Systems Letters},
	author = {Goujaud, Baptiste and Dieuleveut, Aymeric and Taylor, Adrien},
	year = {2023},
	pages = {2485--2490},
}

@inproceedings{zhou_direct_2019,
	title = {Direct acceleration of {SAGA} using sampled negative momentum},
	issn = {2640-3498},
	url = {https://proceedings.mlr.press/v89/zhou19c.html},
	language = {en},
	urldate = {2026-09-28},
	booktitle = {Proceedings of the {Twenty}-{Second} {International} {Conference} on {Artificial} {Intelligence} and {Statistics}},
	publisher = {PMLR},
	author = {Zhou, Kaiwen and Ding, Qinghua and Shang, Fanhua and Cheng, James and Li, Danli and Luo, Zhi-Quan},
	month = apr,
	year = {2019},
	pages = {1602--1610},
}

@inproceedings{shalev-shwartz_accelerated_2014,
	title = {Accelerated proximal stochastic dual coordinate ascent for regularized loss minimization},
	issn = {1938-7228},
	url = {https://proceedings.mlr.press/v32/shalev-shwartz14.html},
	language = {en},
	urldate = {2026-09-28},
	booktitle = {Proceedings of the 31st {International} {Conference} on {Machine} {Learning}},
	publisher = {PMLR},
	author = {Shalev-Shwartz, Shai and Zhang, Tong},
	month = jan,
	year = {2014},
	pages = {64--72},
}

@inproceedings{zhou_simple_2018,
	title = {A simple stochastic variance reduced algorithm with fast convergence rates},
	issn = {2640-3498},
	url = {https://proceedings.mlr.press/v80/zhou18c.html},
	language = {en},
	urldate = {2026-09-28},
	booktitle = {Proceedings of the 35th {International} {Conference} on {Machine} {Learning}},
	publisher = {PMLR},
	author = {Zhou, Kaiwen and Shang, Fanhua and Cheng, James},
	month = jul,
	year = {2018},
	pages = {5980--5989},
}

@article{allen-zhu_katyusha_2018,
	title = {Katyusha: {The} {First} {Direct} {Acceleration} of {Stochastic} {Gradient} {Methods}},
	volume = {18},
	issn = {1533-7928},
	shorttitle = {Katyusha},
	url = {http://jmlr.org/papers/v18/16-410.html},
	number = {221},
	urldate = {2026-09-30},
	journal = {Journal of Machine Learning Research},
	author = {Allen-Zhu, Zeyuan},
	year = {2018},
	pages = {1--51},
}

@article{nesterov_method_1983,
	title = {A method for solving the convex programming problem with convergence rate {O}(1/k{\textasciicircum}2)},
	url = {https://www.semanticscholar.org/paper/A-method-for-solving-the-convex-programming-problem-Nesterov/8d3a318b62d2e970122da35b2a2e70a5d12cc16f},
	urldate = {2024-10-10},
	journal = {Proceedings of the USSR Academy of Sciences},
	author = {Nesterov, Y.},
	year = {1983},
	pages = {543--547},
}

@article{beck_fast_2009,
	title = {A fast iterative shrinkage-thresholding algorithm for linear inverse problems},
	volume = {2},
	issn = {1936-4954},
	url = {http://epubs.siam.org/doi/10.1137/080716542},
	doi = {10.1137/080716542},
	language = {en},
	number = {1},
	urldate = {2023-11-16},
	journal = {SIAM Journal on Imaging Sciences},
	author = {Beck, Amir and Teboulle, Marc},
	month = jan,
	year = {2009},
	pages = {183--202},
}

@article{pepit2024,
  title={{PEPit}: computer-assisted worst-case analyses of first-order optimization methods in {P}ython},
  author={Goujaud, Baptiste and Moucer, C\'eline and Glineur, Fran\c{c}ois and Hendrickx, Julien and Taylor, Adrien and Dieuleveut, Aymeric},
  journal={Math.~Prog.~Comp.},
  volume={16},
  pages={337–367},
  year={2024},
  publisher={Springer},
  doi={https://doi.org/10.1007/s12532-024-00259-7}
}

\end{document}